\documentclass[12pt, a4paper]{article}
\usepackage{mathptmx}  \usepackage[T1]{fontenc}
\usepackage{graphicx, theorem, color} \definecolor{webblue}{rgb}{0,0,0.6}
\usepackage[colorlinks=true, linkcolor=webblue, citecolor=webblue,
 urlcolor=webblue, pdfstartview=FitBH,
 pdftitle={Core entropy and biaccessibility of quadratic polynomials},
 pdfauthor={Wolf Jung},
 bookmarksopen=true, bookmarksnumbered=true, pdfpagemode=UseNone]{hyperref}
\makeatletter 
\long\def\@makecaption#1#2{%
  \vskip\abovecaptionskip
  \sbox\@tempboxa{{\small\textbf{#1}: #2}}%
  \ifdim \wd\@tempboxa >\hsize
    {\small\textbf{#1}: #2}\par
  \else
    \global \@minipagefalse
    \hb@xt@\hsize{\hfil\box\@tempboxa\hfil}%
  \fi
  \vskip\belowcaptionskip}
\makeatother
\theoremstyle{break}
\theorembodyfont{\itshape}
   \newtheorem{thm}{Theorem}[section]
\newtheorem{dfn}[thm]{Definition} \newtheorem{prop}[thm]{Proposition}
\newtheorem{lem}[thm]{Lemma} \newtheorem{conj}[thm]{Conjecture}
\theorembodyfont{\upshape}
\newtheorem{rmk}[thm]{Remark}   \newtheorem{xmp}[thm]{Example}
\newcommand{\be}{\begin{equation}} \newcommand{\ee}{\end{equation}}
\renewcommand{\hat}{\widehat} \renewcommand{\tilde}{\widetilde}
\newcommand{\mybox}{\hspace*{\fill}\rule{2mm}{2mm}}
\DeclareSymbolFont{AMSb}{U}{msb}{m}{n}
\DeclareSymbolFontAlphabet{\mathbb}{AMSb}
\newcommand{\C}{\mathbb{C}}
\newcommand{\Q}{\mathbb{Q}}  \newcommand{\R}{\mathbb{R}}
\newcommand{\Z}{\mathbb{Z}}
  \newcommand{\disk}{\mathbb{D}}

\newcommand{\e}[1]{\displaystyle {\rm e}^{\displaystyle #1}}
\let\inodot\i  \renewcommand{\i}{\mathrm{i}}
\newcommand\brevei{\u\inodot}   
\newcommand{\eps}{\varepsilon}  
\renewcommand{\phi}{\varphi}
\renewcommand{\r}{\mathcal{R}}  \renewcommand{\P}{\mathcal{P}}
  \newcommand{\K}{\mathcal{K}}
\newcommand{\M}{\mathcal{M}}  \newcommand{\sM}{{\scriptscriptstyle M}}
\renewcommand{\O}{\mathcal{O}}  
\renewcommand{\top}{{\scriptscriptstyle\mathrm{top}}}
\newcommand{\comb}{{\scriptscriptstyle\mathrm{comb}}}
\date{}
\author{Wolf Jung\\
{\small Gesamtschule Aachen-Brand}\\
{\small Rombachstrasse 99, 52078 Aachen, Germany.}\\
{\small E-mail: \href{mailto:jung@mndynamics.com}{jung@mndynamics.com}}}
\title{Core entropy and biaccessibility\\of quadratic polynomials}
\begin{document} \maketitle
\begin{abstract}
\noindent For complex quadratic polynomials,  the topology of the
Julia set and the dynamics are understood from another perspective by
considering the Hausdorff dimension of biaccessing angles and the core
entropy:  the topological entropy on the Hubbard tree.  These quantities
are related according to Thurston.  Tiozzo \cite{tiob} has shown
continuity on principal veins of the Mandelbrot set $\M$.  This result
is extended to all veins here,  and it is shown that continuity with
respect to the external angle $\theta$ will imply continuity in the
parameter $c$.  Level sets of the biaccessibility dimension are
described,  which are related to renormalization.  H\"older asymptotics
at rational angles are found,  confirming the H\"older exponent given
by Bruin--Schleicher \cite{bshb}. 
Partial results towards local maxima at dyadic angles are obtained as
well,  and a possible self-similarity of the dimension as a function of
the external angle is suggested.
\end{abstract}

\section{Introduction} 
For a real unimodal map $f(x)$, the topological entropy $h=\log\lambda$
is quantifying the complexity of iteration: e.g., the number of monotonic
branches of $f^n(x)$ grows like $\lambda^n$. Moreover, $f(x)$ is
semi-conjugate to a tent map of slope $\pm\lambda$; so $\lambda$ is an
averaged rate of expansion, which is a topological invariant
\cite{mtimi, miszepm, alm}. Consider a complex quadratic polynomial
$f_c(z)=z^2+c$ with its filled Julia set $\K_c$\,, which is defined in
\textbf{Section~\ref{2}}:
\begin{itemize}
\item At least in the postcritically finite case, the interesting dynamics
happens on the Hubbard tree $T_c\subset\K_c$\,: other arcs are iterated
homeomorphically to $T_c$\,, which is folded over itself, producing chaotic
behavior. The core entropy $h(c)$ is the topological entropy of $f_c(z)$
on $T_c$ \cite{alfa, taoli, tabstr}.
\item On the other hand, for $c\neq-2$ the external angles of these arcs have
measure 0, and the endpoints of $\K_c$ correspond to angles of full measure
in the circle. This phenomenon is quantified by the Hausdorff dimension
$B_\top(c)$ of biaccessing angles: the biaccessibility dimension
\cite{zbqj, zdun, ssms, bks, mesh, bshb}.
\end{itemize}
According to Thurston, these quantities are related by
$h(c)=\log2\cdot B_\top(c)$, which allows to combine tools from different
approaches: e.g., $B_\top(c)$ is easily defined for every parameter $c$ in the
Mandelbrot set $\M$, but hard to compute explicitly. $h(c)$ is easy to compute
and to analyze when $f_c(z)$ is postcritically finite. But for some parameters
$c$, the core may become too large by taking the closure of the connected hull
of the critical orbit in $\K_c$\,.

The postcritically finite case is discussed in \textbf{Section~\ref{3}}.
A Markov matrix $A$ is associated to the Hubbard tree, describing transitions
between the edges. Its largest eigenvalue $\lambda$ gives the core entropy
$h(c)=\log\lambda$. An alternative matrix is due to Thurston
\cite{gaotl, gao}. Or $\lambda$ is obtained from matching conditions
for a piecewise-linear model with constant expansion rate. This approach is
more convenient for computing specific examples, but the matrix $A$ may be
easier to analyze. Known results are extended to all renormalizable maps
$f_c(z)$. Here $A$ is reducible or imprimitive, and its blocks are compared:
pure satellite renormalization gives a rescaling
$h(c_p\ast\hat c)=\frac1p\,h(\hat c)$\,. Lower bounds of $h(c)$ for
$\beta$-type Misiurewicz points and for primitive centers show that in the
primitive renormalizable case, the dynamics on the small Julia sets is
negligible in terms of entropy, so $h(c)$ is constant on maximal-primitive
small Mandelbrot sets $\M_p\subset\M$. Moreover, it is strictly monotonic
between them. An alternative proof of $h(c)=\log2\cdot B_\top(c)$ shows that
the external angles of $T_c$ have finite positive Hausdorff measure.

In \textbf{Section~\ref{4}}, the biaccessibility dimension is defined
combinatorially for every angle $\theta\in S^1$ and topologically for every
parameter $c\in\M$. When $c$ belongs to the impression of the parameter ray
$\r_\sM(\theta)$, we have $B_\top(c)=B_\comb(\theta)$ \cite{bks}. This
relation means that non-landing dynamic rays have angles of negligible
Hausdorff dimension, but a discussion of non-local connectivity is avoided
here by generalizing results from the postcritically finite case: the
biaccessibility dimension is constant on maximal-primitive Mandelbrot sets,
and strictly monotonic between them. Components of a level set of positive
$B_\top(c)$ are maximal-primitive Mandelbrot sets or points. Examples of
accumulation of point components are discussed. In \cite{bshb, tiob} the
Thurston relation $h(c)=\log2\cdot B_\top(c)$ was obtained for all parameters
$c$, such that the core $T_c$ is topologically finite. The proof extends to
compact trees with infinitely many endpoints.

In an email of March 2012 quoted in \cite{gaotl}, Thurston announced proofs of
continuity for $B_\comb(\theta)$ by Hubbard, Bruin--Schleicher, and himself.
In May 2012, a proof with symbolic dynamics was given in \cite{bshb}, but it
is currently under revision. In the present paper, it is shown that continuity
of $B_\comb(\theta)$ on $S^1$ will imply continuity of $B_\top(c)$ on $\M$.
Again, a discussion of non-local connectivity can be avoided, since the
biaccessibility dimension is constant on primitive Mandelbrot sets. See
version~2 of \cite{bshb} for an alternative argument.
Tiozzo \cite{tiob} has shown that $h(c)$ and $B_\top(c)$ are continuous
on principal veins of $\M$; this result is extended to all veins here.

\begin{figure}[h!t!b!]
\unitlength 0.001\textwidth 
\begin{picture}(990, 1055)
\put(10, 740){\includegraphics[width=0.42\textwidth]{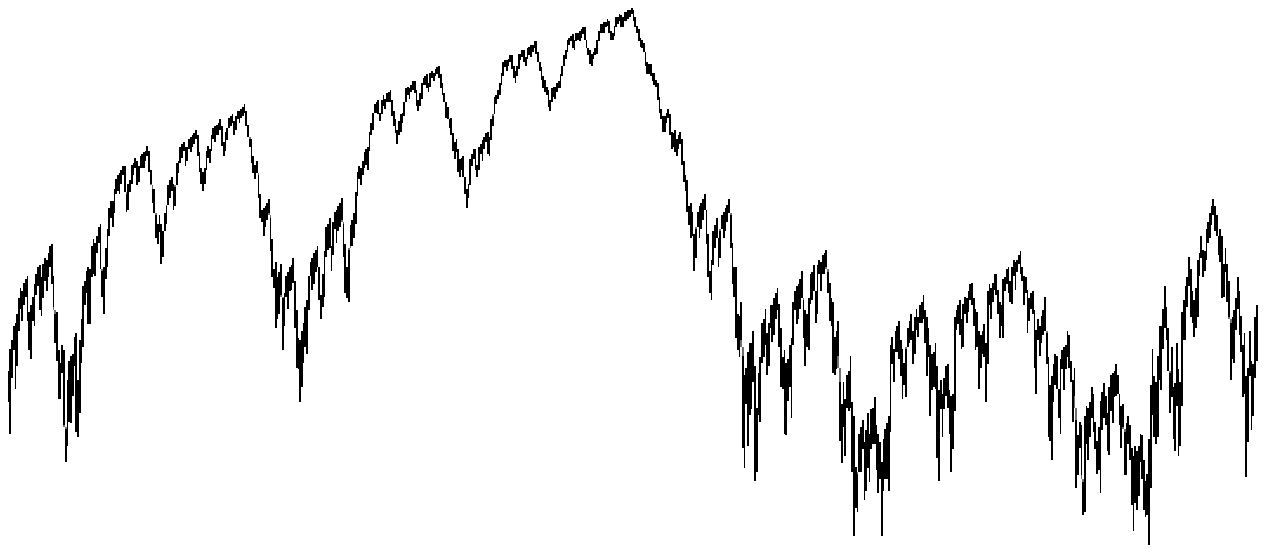}}
\put(570, 740){\includegraphics[width=0.42\textwidth]{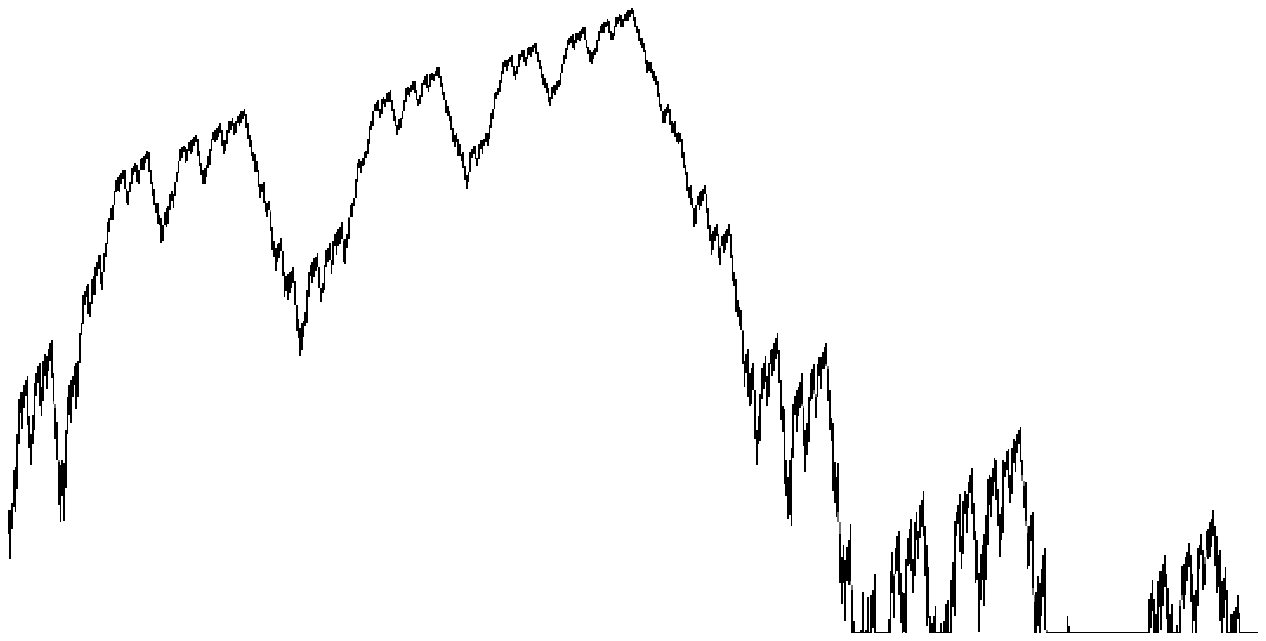}}
\put(10, 370){\includegraphics[width=0.42\textwidth]{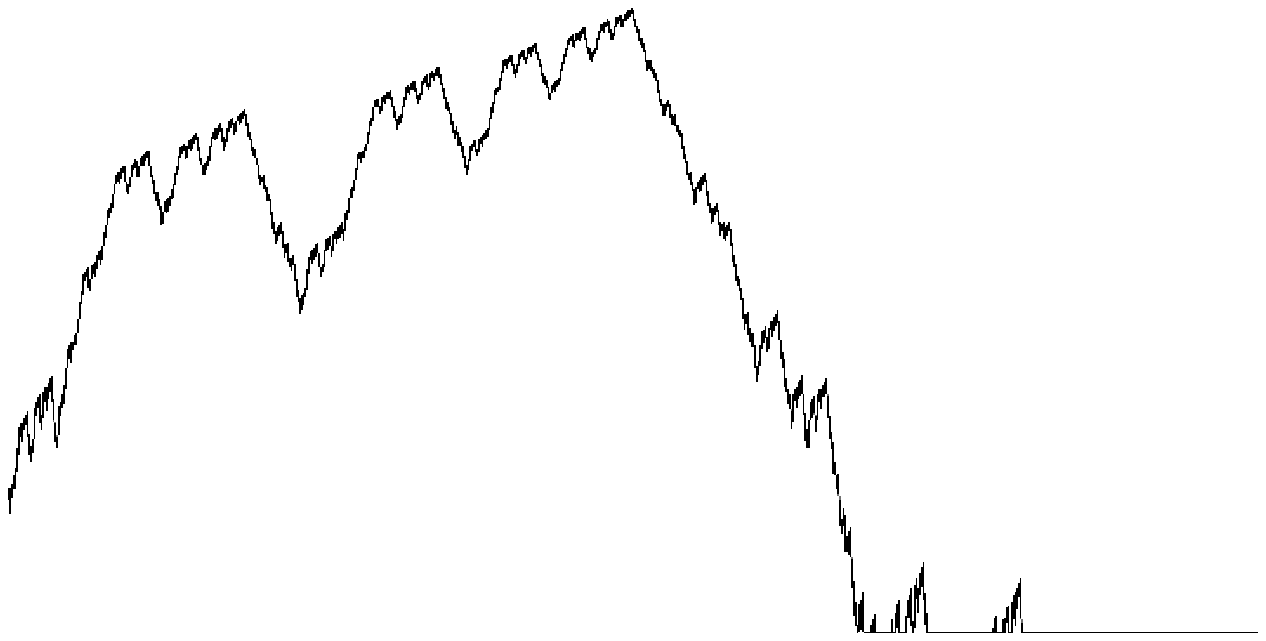}}
\put(570, 370){\includegraphics[width=0.42\textwidth]{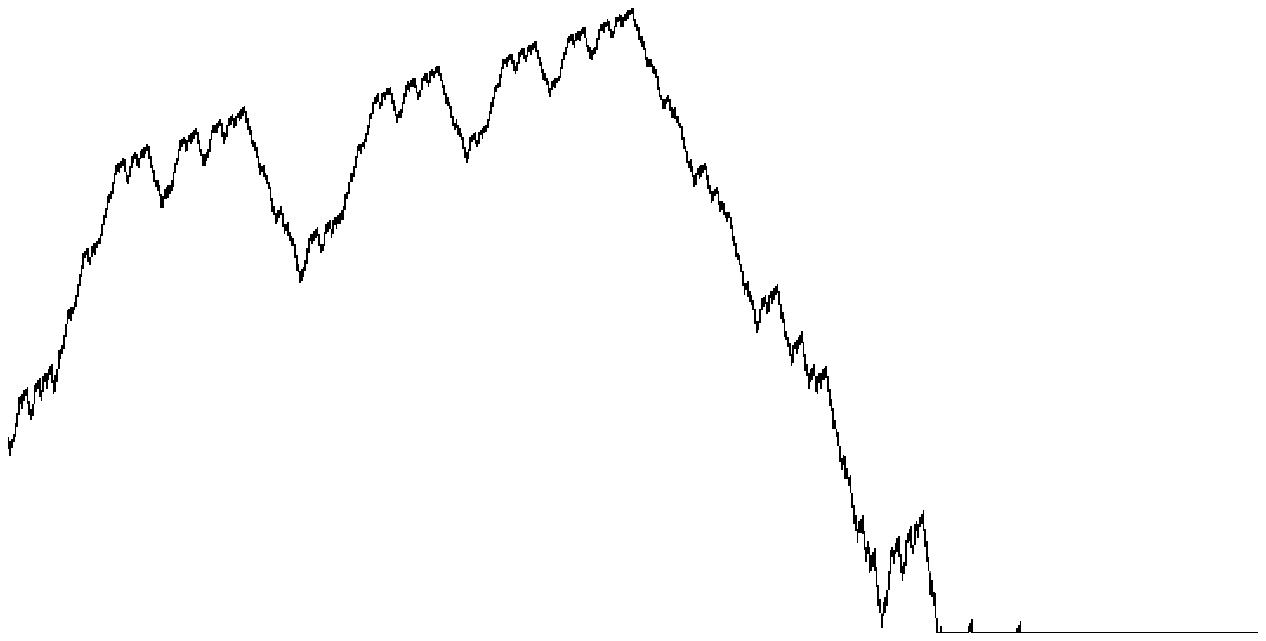}}
\put(10, 0){\includegraphics[width=0.42\textwidth]{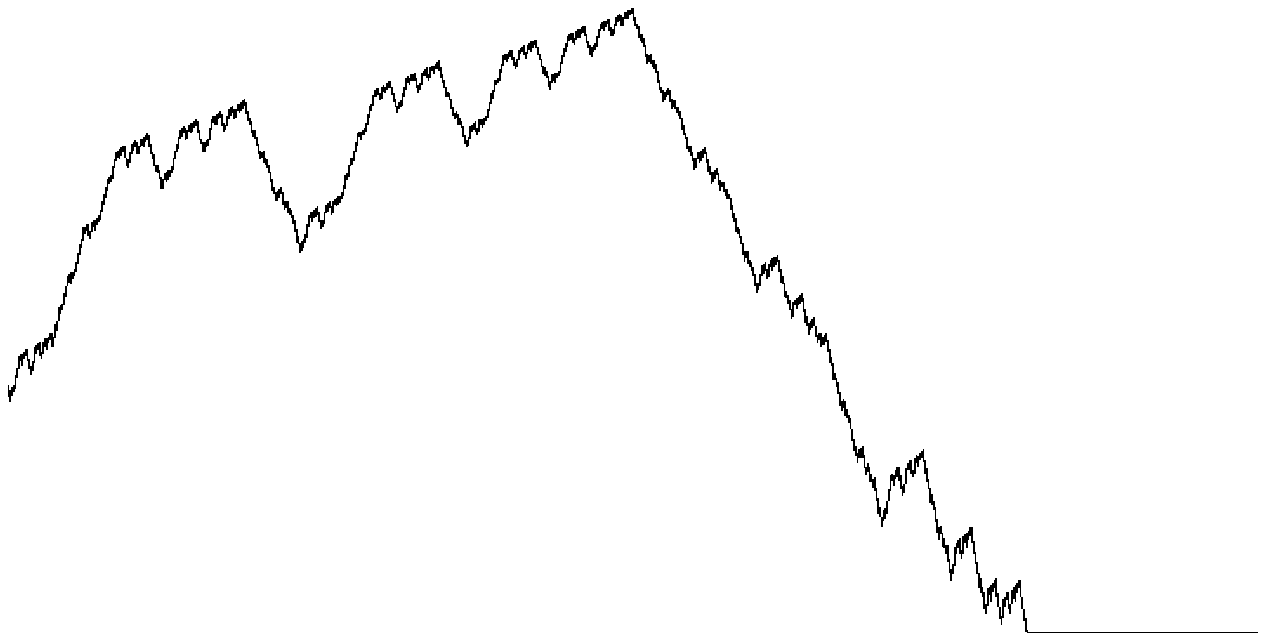}}
\put(570, 0){\includegraphics[width=0.42\textwidth]{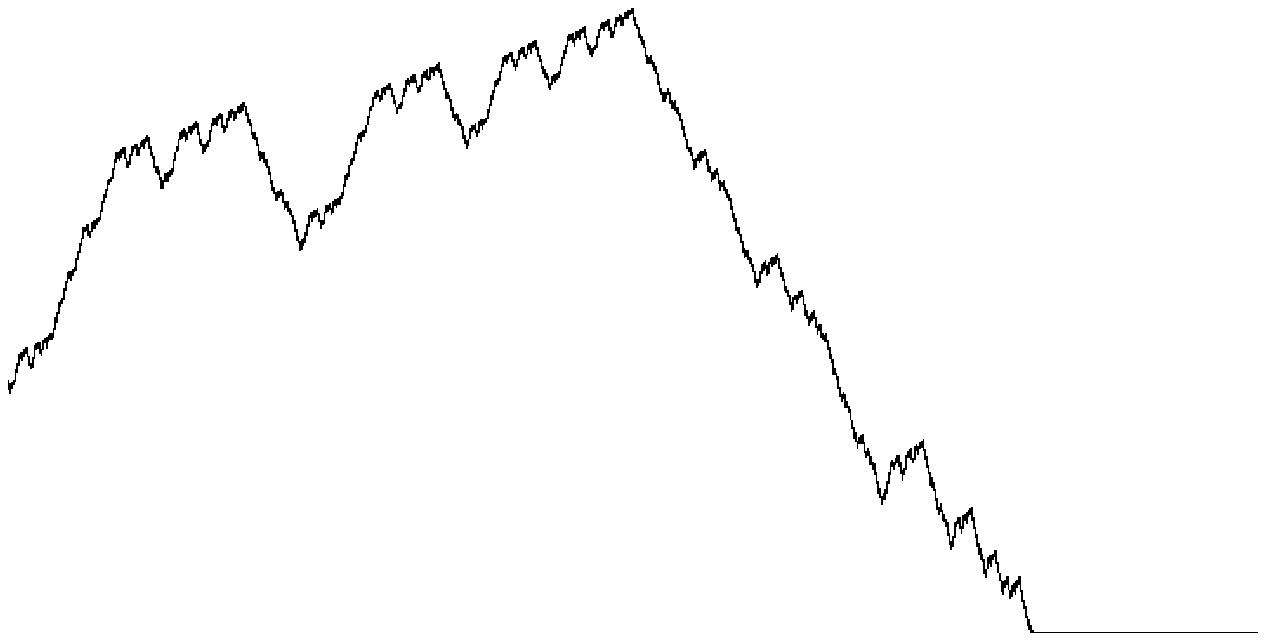}}
\thinlines
\multiput(10, 740)(560, 0){2}{\line(1, 0){420}}
\multiput(10, 1055)(560, 0){2}{\line(1, 0){420}}
\multiput(10, 740)(420, 0){2}{\line(0, 1){315}}
\multiput(570, 740)(420, 0){2}{\line(0, 1){315}}
\multiput(10, 370)(560, 0){2}{\line(1, 0){420}}
\multiput(10, 685)(560, 0){2}{\line(1, 0){420}}
\multiput(10, 370)(420, 0){2}{\line(0, 1){315}}
\multiput(570, 370)(420, 0){2}{\line(0, 1){315}}
\multiput(10, 0)(560, 0){2}{\line(1, 0){420}}
\multiput(10, 315)(560, 0){2}{\line(1, 0){420}}
\multiput(10, 0)(420, 0){2}{\line(0, 1){315}}
\multiput(570, 0)(420, 0){2}{\line(0, 1){315}}
\multiput(500, 157)(0, 370){3}{\makebox(0, 0)[cc]{$\rightarrow$}}
\multiput(500, 342)(0, 370){2}{\makebox(0, 0)[cc]{$\swarrow$}}
\put(35, 1035){\makebox(0, 0)[lt]{$n=0$}}
\put(595, 1035){\makebox(0, 0)[lt]{$n=1$}}
\put(35, 665){\makebox(0, 0)[lt]{$n=2$}}
\put(595, 665){\makebox(0, 0)[lt]{$n=3$}}
\put(35, 295){\makebox(0, 0)[lt]{$n=4$}}
\put(595, 295){\makebox(0, 0)[lt]{$n=5$}}
\end{picture} \caption[]{\label{Fssb}
The biaccessibility dimension is related to the growth factor $\lambda$
by $B_\comb(\theta)=\log(\lambda(\theta))/\log2$.
Consider zooms of $\lambda(\theta)$ centered at $\theta_0=1/4$ with
$\lambda_0=1.69562077$. The width is $0.201\times2^{-n}$ and the height
is $1.258\times\lambda_0^{-n}$\,. There seems to be a local maximum at
$\theta_0$ and a kind of self-similarity with respect to the combined
scaling by $2$ and by $\lambda_0$\,. See Example~\protect\ref{X14}
for the asymptotics of $B_\comb(\theta)$ on specific sequences of
angles $\theta_n\to\theta_0$\,.}
\end{figure}

In \textbf{Section~\ref{5}}, statements of Bruin--Schleicher, Zakeri, and
Tiozzo \cite{bks, bshb, zerm, zbpm, tiob} on the biaccessibility dimension of
$\M$ are generalized to arbitrary pieces. In \textbf{Section~\ref{6}},
Markov matrices are used again to show a geometric scaling behavior of the
core entropy for specific sequences of angles,
which converge to rational angles; for these examples the H\"older
exponent of $B_\comb(\theta)$ given in \cite{bshb} is optimal. The
asymptotics of sequences suggests the question, whether the graph of
$B_\comb(\theta)$ is self-similar; cf.~the example in Figure~\ref{Fssb}.
Partial results towards the Tiozzo Conjecture \cite{tiob} are obtained as
well, which is concerning local maxima of $B_\comb(\theta)$ at dyadic angles.
Some computations of characteristic polynomials are sketched in
\textbf{Appendix~\ref{Asm}}. For the real case, statements on
piecewise-linear models \cite{mtimi} and on the distribution of Galois
conjugates \cite{tjack, tedo, tiom} are reported in
\textbf{Appendix~\ref{Amt}} to round off the discussion. See \cite{bks, bshb}
for the approach with symbolic dynamics and \cite{gaotl, gao} for the
structure of critical portraits.

The present paper aims at a systematic exposition of algebraic and
analytic aspects; so it contains a mixture of well-known, extended, and new
results. I have tried to give proper credits and references to previous or
independent work, and I apologize for possible omissions.
Several people are working on reconstructing and extending Bill Thurston's
results. I have been inspired by hints from or discussions with
Henk Bruin, Gao Yan, Sarah Koch, Michael Mertens,
Dierk Schleicher, Tan Lei, and Giulio Tiozzo.

\section{Background}\label{2}
A short introduction to the complex dynamics of quadratic polynomials is
given to fix some notations. The definitions of topological entropy and of
Hausdorff dimension are recalled, and concepts for non-negative matrices are
discussed.

\subsection{Quadratic dynamics} \label{2qd}
Quadratic polynomials are parametrized conveniently as $f_c(z)=z^2+c$. The
filled Julia set $\K_c$ contains all points $z$ with a bounded orbit under
the iteration, and the Mandelbrot set $\M$ contains those parameters $c$,
such that $\K_c$ is connected. Dynamic rays $\r_c(\phi)$ are curves
approaching $\partial\K_c$ from the exterior, having the angle $2\pi\phi$ at
$\infty$, such that $f_c(\r_c(\phi))=\r_c(2\phi)$.
They are defined as preimages of straight rays under the B\"ottcher map
$\Phi_c\,:\,\hat\C\setminus\K_c\to\hat\C\setminus\overline\disk$.
Parameter rays $\r_\sM(\theta)$ approach $\partial\M$ \cite{mer, ser}; they
are defined in terms of the Douady map
$\Phi_\sM\,:\,\hat\C\setminus\M\to\hat\C\setminus\overline\disk$ with
$\Phi_\sM(c):=\Phi_c(c)$.
The landing point is denoted by $z=\gamma_c(\phi)$ or $c=\gamma_\sM(\theta)$,
respectively, but the rays need not land for irrational angles, see
Figure~\ref{FFB}. There are two cases of postcritically finite dynamics:
\begin{itemize}
\item When the parameter $c$ is a Misiurewicz point, the critical value
$z=c$ is strictly preperiodic. Both $c\in\partial\M$ and $c\in\partial\K_c$
have the same external angles, which are preperiodic under doubling.
\item When $c$ is the center of a hyperbolic component, the critical orbit
is periodic and contained in superattracting basins. The external angles
of the root of the component coincide with the characteristic angles of the
Fatou basin around $z=c$\,; the characteristic point may have more periodic
angles in the satellite case.
\end{itemize}
In both cases, the Hubbard tree \cite{dcc, hss} is obtained by connecting the
critical orbit with regulated arcs, which are traveling through Fatou basins
along internal rays. Fixing a characteristic angle $\theta$ of $c$, the
circle $S^1=\R/\Z$ is partitioned by the diameter joining $\theta/2$ and 
$(\theta+1)/2$ and the orbit of an angle $\phi\in S^1$ under doubling is
encoded by a sequence of symbols $A,\,B,\,\ast$ or $1,\,0,\,\ast$. There is a
corresponding partition of the filled Julia set, so points $z\in\K_c$ are
described by  symbolic dynamics as well. The kneading sequence is the
itinerary of $\theta$ or of $c$, see \cite{hss, ser, bks}.

$\M$ consists of the closed main cardioid and its
limbs, which are labeled by the rotation number at the fixed point
$\alpha_c$\,; the other fixed point $\beta_c$ is an endpoint of $\K_c$\,.
A partial order on $\M$ is defined such that $c\prec c'$ when $c'$ is
disconnected from $0$ in $\M\setminus\{c\}$. See Sections~\ref{3ren} 
and~\ref{4ctb} for the notion of renormalization
\cite{dhp, daa, hy3, mlcj, sf3, wje}, which is explaining small Julia sets
within Julia sets and small Mandelbrot sets within the Mandelbrot set.
Primitive and satellite renormalization may be nested; a primitive small
Mandelbrot set will be called maximal-primitive, if it is not contained in
another primitive one. A pure satellite is attached to the main cardioid
by a series of satellite bifurcations, so it is not contained in a
primitive Mandelbrot set.

\subsection{Topological entropy} \label{2te}
Suppose $X$ is a compact metric space and $f:X\to X$ is continuous. The
topological entropy is measuring the complexity of iteration from the
growth rate of the number of distinguishable orbits.
The first definition assumes an open cover $U$ and considers the minimal
cardinality $V_U(n)$ of a subcover, such that all points in a set of the
subcover have the same itinerary with respect to $U$ for $n$ steps. See, e.g.,
\cite{dte, taoli}. The second definition is using the minimal number
$V_\eps(n)$ of points, such that every orbit is $\eps$-shadowed by one of
these points for $n$ steps \cite{bow1, mevst}.
We have
\be h_\top(f,\,X) :=
 \sup_U \, \lim_{n\to\infty} \,\frac1n\,\log V_U(n) =
 \lim_{\eps\to0} \, \limsup_{n\to\infty} \,\frac1n\,\log V_\eps(n) \ . \ee
For a continuous, piecewise-monotonic interval map, the growth rate of
monotonic branches (laps) may be used instead, or the maximal growth rate of
preimages \cite{miszepm, dte}. The same result applies to endomorphisms of
a finite tree \cite{alm, taoli}. Moreover, $f$ is semi-conjugate to a
piecewise-linear model of constant expansion rate $\lambda$ when
$h_\top(f,\,X)=\log\lambda>0$; this is shown in \cite{mtimi, dte, alm}
for interval maps and in \cite{bcplm} for tree maps. See Section~\ref{Amt}
for the relation to the kneading determinant and Section~\ref{4ebv}
for continuity results \cite{mtimi, miszepm, mishpc, dte, alm}. If
$\pi:X\to Y$ is a surjective semi-conjugation from $f:X\to X$ to
$g:Y\to Y$, then $h_\top(g,\,Y)\le h_\top(f,\,X)$. Equality follows when every
fiber is finite, but fiber cardinality need not be bounded globally
\cite{bow1, mevst}.

\subsection{Hausdorff dimension} \label{2hd}
The $d$-dimensional Hausdorff measure is a Borel outer measure. It is
defined as follows for a bounded subset $X\subset\R$ or a subset
$X\subset S^1=\R/\Z$\,:
\be\label{eqhdm} \mu_d(X):=\lim_{\eps\to0}\,\inf_U\, \sum_i |U_i|^d \ee
Here the cover $U$ of $X$ is a countable family of intervals $U_i$ of length
$|U_i|\le\eps$. They may be assumed to be open or closed, aligned to nested
grids or not, but the important point is that they may be of different size.
When an interval is replaced with two subintervals, the sum may grow in fact
when $d<1$. In general, the Hausdorff measure of $X$ may be easy to bound
from above by using intervals of the same size, but it will be hard to bound
from below, since this requires to find an optimal cover with intervals of
different sizes.

The Hausdorff dimension $\dim(X)$ is the unique number in $[0,\,1]$, such
that $\mu_d(X)=\infty$ for $0\le d<\dim(X)$ and $\mu_d(X)=0$ for
$\dim(X)<d\le1$. For $d=\dim(X)$, the Hausdorff measure $\mu_d(X)$ may be
$0$, positive and finite, or $\infty$. The Hausdorff dimension of a countable
set is $0$ and the dimension of a countable union is the supremum of the
dimensions. Again, $\dim(X)$ may be easy to bound from above by the box
dimension, which corresponds to equidistant covers, but it is harder to bound
from below. Sometimes this is achieved by constructing a suitable mass
distribution according to the Frostman Lemma \cite{fafr, zerm}.
When $X\subset S^1$ is closed and invariant under doubling $F(\phi)=2\phi$,
the Hausdorff dimension $\dim(X)$ is equal to the box dimension according to
Furstenberg \cite{fu}.

\subsection{Perron--Frobenius theory}\label{2pf}
The Perron theory of matrices with positive elements was extended by Frobenius
to non-negative matrices, see \cite{hojo}. We shall need the following
features of a square matrix $A\ge0$: 
\begin{itemize}
\item There is a non-negative eigenvalue $\lambda$ with non-negative
eigenvector, such that all eigenvalues of $A$ are $\le\lambda$
in modulus. $\lambda$ is bounded above by the maximal sum of rows or columns,
and bounded below by the minimal sum.
\item $A\ge0$ is called reducible, if it is conjugate to a block-triangular
matrix by a permutation. It is irreducible (ergodic) if the corresponding
directed graph is strongly connected.
Then $\lambda$ is positive and it is an algebraically simple eigenvalue.
The eigenvector of $\lambda$ is positive, and other eigenvectors are not
non-negative.
\item An irreducible $A\ge0$ is called primitive (mixing), if the other
eigenvalues have modulus $<\lambda$. Equivalently, $A^n>0$ element-wise
for some $n$.
\item If $A$ is irreducible but imprimitive with $p>1$ eigenvalues of modulus
$\lambda$, its characteristic polynomial is of the form $x^kP(x^p)$. The
Frobenius normal form shows that there are $p$ subspaces mapped cyclically
by $A$.
\item If $B\ge A$ element-wise with $B\neq A$ and $A$ is primitive,
then $\lambda_B>\lambda_A$\,. This is proved by choosing $n$ with $A^n>0$,
noting $B^n\ge A^n$, so $B^{n+1}$ is strictly larger than $A^{n+1}$ in
at least one row and one column, and finally $B^{2n+1}>A^{2n+1}$. Now fix
$\eps>0$ with $B^{2n+1}\ge(1+\eps)A^{2n+1}$ and consider higher powers
to show $\lambda_B\ge\sqrt[2n+1]{1+\eps}\,\lambda_A$\,.
\end{itemize}
To obtain $\lambda$ numerically from $A\ge0$, we do not need to
determine the characteristic polynomial and its roots: for some positive
vector $v_0$ compute $v_n=A^nv_0$ recursively, then
$\lambda=\lim\sqrt[n]{\|v_n\|}$ converges slowly. If $A$ is
irreducible and primitive, $\lambda=\lim{\|v_{n+1}\|}\,/\,{\|v_n\|}$
will converge exponentially fast.

\section{Postcritically finite polynomials and core entropy} \label{3}
Suppose $f_c(z)$ is postcritically finite and consider the Hubbard tree
$T_c$\,. Since the collection of vertices is forward invariant, each edge is
mapped to one edge or to several adjacent edges. Thus the edges form a Markov
partition (strictly speaking, a tessellation). By numbering the edges, the map
is described by a non-negative matrix $A$ with entries $0$ and $1$, such that
the $j$-th column is showing where the $j$-th edge of $T_c$ is mapped by
$f_c(z)$. The Markov matrix $A$ is the transition matrix of the Markov
partition and the adjacency matrix of the Markov diagram. Often the transposed
matrix is used instead. In the preperiodic case, no postcritical point is
mapped to the critical point $z=0$. So we still have a Markov partition when
the two edges at $0$ are considered as one edge, but mapping this edge will
cover the edge before $z=c$ twice, resulting in an entry of $2$ in $A$.

\begin{dfn}[Markov matrix and core entropy]\label{DAlh}
For a postcritically finite quadratic polynomial $f_c(z)$, the Markov matrix
$A$ is the transition matrix for the edges of the Hubbard tree $T_c$\,. Its
highest eigenvalue $\lambda$ gives the core entropy $h(c):=\log\lambda$.
Equivalently, $h(c):=h_\top(f_c\,,\,T_c)$ is the topological entropy
on $T_c$\,.
\end{dfn}

The $i$-th row of $A$ says which edges are mapped to an arc covering the
$i$-th edge. Since $f_c\,:\,T_c\to T_c$ is surjective and at most 2:1,
the sum of each row of $A$ is $1$ or $2$, so the highest eigenvalue of $A$
satisfies $1\le\lambda\le2$. The largest entries of $A^n$ are growing as
$\asymp\lambda^n$ when $\lambda>1$. (Not as $n^k\lambda^n$\,: according to
Section~\ref{3ren}, $A$ may be reducible but $\lambda>1$ corresponds to a
unique irreducible block, which need not be primitive.)
Since the entries of $A^n$ give the number of preimages of edges, the same
estimate applies to the maximal cardinality of preimages $f_c^{-n}(z_0)$
in $T_c$\,. There are various ways to show that $\log\lambda$ is the
topological entropy of $f_c(z)$ on the Hubbard tree $T_c$\,:
\begin{itemize}
\item By expansivity, every edge is iterated to an edge at $z=0$, so the
preimages of $0$ are growing by $\lambda^n$ even if the edges at $0$
correspond to an irreducible block of lower eigenvalues. So the number of
monotonic branches of $f_c^n(z)$ is growing by the same rate, which
determines the topological entropy according to \cite{alm}. 
\item This criterion is obtained from the definition according to
Section~\ref{2te} by showing that the open cover may be replaced with
closed arcs having common endpoints, and that the maximal growth rate is
attained already for the given edges \cite{dte, taoli}.
\item According to \cite{bcplm}, $f_c(z)$ on $T_c$ is semi-conjugate
to a piecewise-linear model with constant expansion rate $\lambda$ when
the topological entropy is $\log\lambda>0$. See also \cite{mtimi, dte} for
real parameters and \cite{tiob} for parameters on veins. Now the highest
eigenvector of the transposed matrix $A'$ is assigning a Markov length to the
edges, such that $f_c(z)$ corresponds to multiplying the length with
$\lambda$. Note that according to the decomposition (\ref{renmat}) in
Section~\ref{3ren}, the edges in a primitive small Julia set have Markov
length $0$ and will be squeezed to points by the semi-conjugation.
\end{itemize}

\subsection{Computing the core entropy} \label{3cce}
Let us start with four examples of postcritically finite parameters
$c$ in the $1/3$-limb of $\M$. The external angle
$\theta=3/15$ gives a primitive center of period $4$.
The other examples are preperiodic and the edges at $z=0$ are united.
$\theta=1/4$ defines a $\beta$-type Misiurewicz point of preperiod $2$ and
$\theta=9/56$ gives an $\alpha$-type Misiurewicz point of preperiod $3$
and ray period $3$; it is satellite renormalizable of period $3$ and the
Markov matrix $A$ is imprimitive of index $3$. And $\theta=1/6$ gives
$c=\i$, a Misiurewicz point with preperiod $1$ and period $2$.

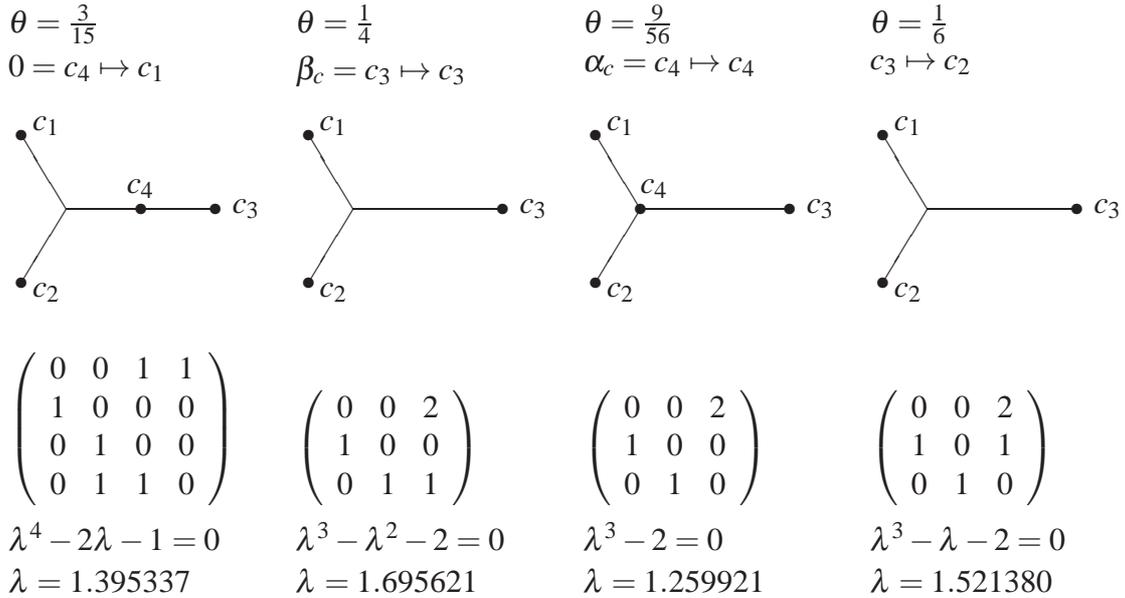
\begin{figure}[h!t!b!]
\unitlength 0.001\textwidth 
\begin{picture}(990, 500) 
\multiput(10, 250)(250, 0){4}{\begin{picture}(230, 180)
\put(50, 90){\line(-3, 5){39}}
\put(50, 90){\line(-3, -5){39}}
\put(50, 90){\line(1, 0){130}}
\put(11, 25){\circle*{10}} \put(21, 25){\makebox(0, 0)[lt]{$c_2$}}
\put(11, 155){\circle*{10}} \put(21, 155){\makebox(0, 0)[lb]{$c_1$}}
\put(180, 90){\circle*{10}} \put(195, 90){\makebox(0, 0)[lc]{$c_3$}}
\end{picture}}
\put(125, 340){\circle*{10}} \put(125, 350){\makebox(0, 0)[cb]{$c_4$}} 
\put(560, 340){\circle*{10}} \put(560, 350){\makebox(0, 0)[lb]{$c_4$}} 
\put(10, 520){\makebox(0, 0)[lt]{$\theta=\frac3{15}$}}
\put(10, 475){\makebox(0, 0)[lt]{$0=c_4\mapsto c_1$}}
\put(10, 25){\makebox(0, 0)[lt]{$\lambda=1.395337$}}
\put(10, 65){\makebox(0, 0)[lt]{$\lambda^4-2\lambda-1=0$}}
\put(10, 80){\makebox(0, 0)[lb]{$\left(\begin{array}{cccc}
0 & 0 & 1 & 1\\
1 & 0 & 0 & 0\\
0 & 1 & 0 & 0\\
0 & 1 & 1 & 0
\end{array}\right)$}}
\put(260, 520){\makebox(0, 0)[lt]{$\theta=\frac14$}}
\put(260, 475){\makebox(0, 0)[lt]{$\beta_c=c_3\mapsto c_3$}}
\put(260, 25){\makebox(0, 0)[lt]{$\lambda=1.695621$}}
\put(260, 65){\makebox(0, 0)[lt]{$\lambda^3-\lambda^2-2=0$}}
\put(260, 80){\makebox(0, 0)[lb]{$\left(\begin{array}{ccc}
0 & 0 & 2\\
1 & 0 & 0\\
0 & 1 & 1
\end{array}\right)$}}
\put(510, 520){\makebox(0, 0)[lt]{$\theta=\frac9{56}$}}
\put(510, 475){\makebox(0, 0)[lt]{$\alpha_c=c_4\mapsto c_4$}}
\put(510, 25){\makebox(0, 0)[lt]{$\lambda=1.259921$}}
\put(510, 65){\makebox(0, 0)[lt]{$\lambda^3-2=0$}}
\put(510, 80){\makebox(0, 0)[lb]{$\left(\begin{array}{ccc}
0 & 0 & 2\\
1 & 0 & 0\\
0 & 1 & 0
\end{array}\right)$}}
\put(760, 520){\makebox(0, 0)[lt]{$\theta=\frac16$}}
\put(760, 475){\makebox(0, 0)[lt]{$c_3\mapsto c_2$}}
\put(760, 25){\makebox(0, 0)[lt]{$\lambda=1.521380$}}
\put(760, 65){\makebox(0, 0)[lt]{$\lambda^3-\lambda-2=0$}}
\put(760, 80){\makebox(0, 0)[lb]{$\left(\begin{array}{ccc}
0 & 0 & 2\\
1 & 0 & 1\\
0 & 1 & 0
\end{array}\right)$}}
\end{picture} \caption[]{\label{F4pcf13}
Examples of Hubbard trees and Markov matrices defined by an
angle $\theta$. Here $f_c(z)$ maps $c=c_1\mapsto c_2\mapsto c_3\mapsto c_4$
and the edges are numbered such that the first edge is before $c_1=c$.}
\end{figure}

Instead of computing the characteristic polynomial of the Markov matrix $A$,
we may use a piecewise-linear model with expansion constant $\lambda>1$ to be
matched: in Figure~\ref{F4pcf13}, the edges or arcs from $\alpha_c$ to
$c_1$\,, $c_2$\,, $c_3$ have length $\lambda$, $\lambda^2$, $\lambda^3$ when
$[0,\,\pm\alpha_c]$ has length $1$.\\
For $\theta=3/15$, we have $[-\alpha_c\,,\,c_3]\mapsto[\alpha_c\,,\,c_4]$,
so $\lambda(\lambda^3-2)=1$.\\
For $\theta=1/4$, we have $[-\alpha_c\,,\,c_3]\mapsto[\alpha_c\,,\,c_3]$,
so $\lambda(\lambda^3-2)=\lambda^3$.\\
For $\theta=9/56$, we have $c_3=-\alpha_c$\,,
so $\lambda^3=2$.\\
For $\theta=1/6$, we have $[-\alpha_c\,,\,c_3]\mapsto[\alpha_c\,,\,c_2]$,
so $\lambda(\lambda^3-2)=\lambda^2$.

\begin{xmp}[Lowest periods and preperiods in limbs]\label{Xprinc}
The $p/q$-limb of $\M$ contains those parameters $c$, such that the fixed
point $\alpha_c$ has the rotation number $p/q$. The principal vein is the arc
from $0$ to the $\beta$-type Misiurewicz point of preperiod $q-1$. The
examples in Figure~\ref{F4pcf13} can be generalized by considering sequences
of Markov matrices or by replacing $\lambda^3-2$ with $\lambda^q-2$ for the
piecewise-linear models. This gives the following polynomials for $\lambda$:\\
$\beta$-type Misiurewicz point of preperiod $q-1$\,: $\,x^q-x^{q-1}-2=0\;$
\cite{alfa}\\ 
Primitive center of period $q+1$\,: $\,x^{q+1}-2x-1=0$\\
$\alpha$-type Misiurewicz point of preperiod $q$\,: $\,x^q=2$\\
These equations show that $h(c)$ is not H\"older continuous with respect to
the external angle $\theta$ as $\theta\to0$, see Section~\ref{4cc}.
\end{xmp}

\begin{xmp}[Sequences on principal veins]\label{Xpriser}
On the principal vein of the $p/q$-limb, the $\beta$-type Misiurewicz point
is approached by a sequence of centers $c_n$ and $\alpha$-type Misiurewicz
points $a_n$ of increasing periods and preperiods. The corresponding
polynomials for $\lambda_n$ are obtained from piecewise-linear models again,
or by considering a sequence of matrices.
The polynomial is simplified by summing a geometric series and multiplying
with $\lambda-1$:\\
Center $c_n$ of period $n\ge q+1$\,: $\,x^{n+1}-x^n-2x^{n+1-q}+x+1=0$\\
$\alpha$-type Misiurewicz point $a_n$ of preperiod $n\ge q$\,:
$\,x^{n+1}-x^n-2x^{n+1-q}+2=0$\\
These polynomials imply monotonicity of $\lambda_n$ and give geometric
asymptotics by writing
$x^{n+1}-x^n-2x^{n+1-q}=x^{n+1-q}\cdot(x^q-x^{q-1}-2)$\,;
note that the largest root $\lambda_0$ of the latter polynomial corresponds
to the endpoint of the vein according to Example~\ref{Xprinc}:\\
For $c_n$ we have $\lambda_n\sim\lambda_0-K_c\cdot\lambda_0^{-n}$
with $K_c=\frac{\lambda_0+1}{q-(q-1)/\lambda_0}>0$.\\
For $a_n$ we have $\lambda_n\sim\lambda_0-K_a\cdot\lambda_0^{-n}$
with $K_a=\frac2{q-(q-1)/\lambda_0}>0$.\\
See Proposition~\ref{Pasympbeta} and Appendix~\ref{Asm} for a detailed
computation and Remark~\ref{Rash} for the relation to H\"older continuity.
\end{xmp}

\begin{xmp}[Sequences on the real axis]\label{Xstsa}
The real axis is the principal vein of the $1/2$-limb; setting $q=2$
in Example~\ref{Xpriser}, dividing by $\lambda+1$ and noting $\lambda_0=2$
at the endpoint $c=-2$ gives:\\
$c_n$ of period $n\ge3$\,: $\,x^n-2x^{n-1}+1=0$\,,
$\lambda_n\sim2-2\cdot2^{-n}$\\
$a_n$ of preperiod $n\ge 2$\,:
$\,x^{n+1}-x^n-2x^{n-1}+2=0$\,,
$\lambda_n\sim2-\frac43\cdot2^{-n}$\\
Now consider the $\alpha$-type Misiurewicz point $a_2$ with the external angle
$\theta=5/12$ and with $\lambda_a=\sqrt2$, which is the tip of the satellite
Mandelbrot set of period $2$. It is approached from above (with
respect to $\prec$, i.e., from the left) by centers $c_n'$ of periods
$n=3,\,5,\,7,\,\dots$ related to the \v{S}harkovski\brevei{} ordering.
The entropy was computed by \v{S}tefan \cite{sts}, and the centers of even
periods before $a_2$ (to the right) are treated analogously. Again, geometric
asymptotics are obtained from the sequence of polynomials for $\lambda_n'$\,:\\
$c_n'$ of period $n=3,\,5,\,7,\,\dots$\,: $\,x^n-2x^{n-2}-1=0$\,,
$\lambda_n'\sim\lambda_a+\frac1{\sqrt2}\cdot\lambda_a^{-n}$\\
$c_n'$ of period $n=4,\,6,\,8,\,\dots$\,: $\,x^n-2x^{n-2}+1=0$\,,
$\lambda_n'\sim\lambda_a-\frac1{\sqrt2}\cdot\lambda_a^{-n}$\\
Note that for even periods $n$, the polynomial for $\lambda_n'$ is
imprimitive of index $2$, which is related to the satellite renormalization
according to Section~\ref{3ren}.
\end{xmp}

Constructing the Hubbard tree $T_c$ from an external angle $\theta$ is quite
involved \cite{bks}. In \cite{gaotl, gao} an alternative matrix $F$ by
Thurston is described, which is obtained from an external angle without
employing the Hubbard tree. Actually, only the kneading sequence of the angle
is required to determine $F$ from the parts of $T_c\setminus\{0\}$\,:

\begin{prop}[Alternative matrix by Thurston and Gao]\label{Path}
From a rational angle or from a $\ast$-periodic or preperiodic kneading
sequence, construct a transition matrix $F$. The basic vectors represent
non-oriented arcs between postcritical points $c_j=f_c^j(0)$, $j\ge1$,
and $[c_j\,,\,c_k]$ is mapped to $[c_{j+1}\,,\,c_{k+1}]$ by $F$ unless
its endpoints are in different parts of $\K_c\setminus\{0\}$\,:
then it is mapped to $[c_1\,,\,c_{j+1}]+[c_1\,,\,c_{k+1}]$.\\
Now the largest eigenvalues of the Thurston matrix $F$ and the
Markov matrix $A$ coincide.
\end{prop}

This combinatorial definition corresponds to the fact that an arc covering
$z=0$ is mapped 2:1 to an arc at $z=c$ by $f_c(z)$. Arcs at $c_0$ are omitted
in the preperiodic case, because they would generate a diagonal 0-block
anyway.  Note that in general $F$ is considerably larger than $A$; it will
contain large nilpotent blocks and it may contain additional blocks, which
seem to be cyclic. 
In the case of $\beta$-type Misiurewicz points, a small irreducible block of
$F$ is obtained in Proposition~\ref{Pbeta}.3.

\textbf{Proof:} Gao \cite{gao} is using a non-square incidence matrix $C$,
which is mapping each arc to a sum of edges, so $AC=CF$.
Consider the non-negative Frobenius eigenvectors to
obtain equality of the highest eigenvalues: if
$Fy=\lambda_{\scriptscriptstyle F}y$ then $Cy$ is an
eigenvector of $A$ with eigenvalue $\lambda_{\scriptscriptstyle F}$\,. The
transposed matrices satisfy $F'C'=C'A'$ and if
$A'x=\lambda_{\scriptscriptstyle A}x$, then $C'x$ is an eigenvector of
$F'$ with eigenvalue $\lambda_{\scriptscriptstyle A}$\,. (Note that $Cy$ and
$C'x$ are not $0$, because $C$ has a non-zero entry in each row and each
column.)

As an alternative argument, define a topological space
$X_c$ as a union of arcs $[c_j\,,\,c_k]\subset\K_c$\,, which are considered
to be disjoint except for common endpoints. There is a natural projection
$\pi_c\,:\,X_c\to T_c$ and a lift $F_c\,:\,X_c\to X_c$ of $f_c(z)$, such that
$\pi_c$ is a semi-conjugation and $F$ is the transition matrix of $F_c$\,.
Now any $z\in T_c$ has a finite fiber
$\pi_c^{-1}(z)=\{x_i\}\subset X_c$ and we have the disjoint union
$\bigcup F_c^{-n}(x_i)=\pi_c^{-1}(f_c^{-n}(z))$. Choosing $z$ such that the
cardinality of $f_c^{-n}(z)$ is growing by $\lambda_{\scriptscriptstyle A}^n$
shows $\lambda_{\scriptscriptstyle A}\le\lambda_{\scriptscriptstyle F}$\,. And
choosing $z$ such that the cardinality of $F_c^{-n}(x_1)$ is growing by
$\lambda_{\scriptscriptstyle F}^n$ gives
$\lambda_{\scriptscriptstyle F}\le\lambda_{\scriptscriptstyle A}$\,.

Note that $F$ is determined from the kneading sequence of $c_1=c$, which
can be obtained from the external angle $\theta$ as an itinerary.
Alternatively, consider a matrix where the basic vectors represent pairs of
angles $\{2^{j-1}\theta,\,2^{k-1}\theta\}$\,; it will be the same as $F$
except 
in the preperiodic satellite case, where $c_1$ is entering a repelling
$p$-cycle of ray period $rp>p$. (This cycle contains the characteristic
point of a satellite component before the Misiurewicz point $c$.)
The matrix of transitions between pairs of angles will be different from $F$,
but the largest eigenvalue will be the same: pairs of equivalent angles
are permuted cyclically and correspond to eigenvalues of modulus $1$
\cite{gao}. Arcs of $X_c$ with at least one $p$-periodic endpoint are
represented by several pairs of angles, but when a topological space is built
from multiple copies of arcs, it comes with a semi-conjugation to $X_c$ or to
$T_c$ again. \mybox

\subsection{Estimates of the core entropy} \label{3ece}
The edges of the Hubbard tree are connecting the marked points: the critical
orbit $f_c^n(0)$, $n\ge0$, which includes all endpoints, and additional
branch points.

\begin{lem}[Modified Markov matrix]\label{LHtree}
For a postcritically finite $f_c(z)$, the Hubbard tree $T_c$ and the
associated Markov matrix $A$ may be changed as follows, without changing
the highest eigenvalue $\lambda$, and without changing irreducibility or
primitivity $($except for item~$4$ and for $c=-1$ in item~$2)$:

$1$. If $c$ is preperiodic, the two edges at $z=0$ may be considered as one
edge by removing $z=0$ from the marked points; thus an eigenvalue $0$ has
been removed from $A$.

$2$. For $c$ in the $1/2$-limb, the unmarked fixed point $\alpha_c$ may be
marked, splitting one edge in two. This gives an additional eigenvalue
of $-1$.

$3$. An unmarked preimage of a marked periodic or preperiodic point may
be marked.

$4$. Extend the Hubbard tree $T_c$ by attaching edges towards the fixed
point $\beta_c$ and/or some of its preimages. This gives additional
eigenvalues of $1$ and $0$, and it makes $A$ reducible.
\end{lem}

These modifications may be combined, and item~3 can be applied recursively.
The \textbf{proof} is deferred to Appendix~\ref{Asm}. Figure~\ref{F4pcf13}
gives examples of item~1. Items~2 and~3 are applied in Section~\ref{3ren}.
Item~4 shows that for a biaccessible parameter $c$, the Hubbard tree $T_c$
may be extended in a uniform way for all parameters on the vein. E.g., in
the real case we may replace $T_c=[c,\,f_c(c)]$ with $[-\beta_c\,,\,\beta_c]$.
Item~4 is used as well to prove Proposition~\ref{Pmaxlimbs}.
Note that a matrix with an eigenvalue 0, -1, or 1 will be reduced by a
conjugation in $GL(\Q^N)$, but only permutations are considered for a
Frobenius irreducible matrix.

\begin{prop}[Monotonicity of core entropy, Penrose and Tao Li]\label{Pmonopcf}
Core entropy is monotonic: for postcritically finite
$c\prec c'$ we have $h(c)\le h(c')$.
\end{prop}

Not all parameters are comparable with the partial order $\prec$\,.
In particular, many parameters $c'$ are approached by branch points $c$
before them, and parameters $c''$ in different branches are not
comparable to $c'$. --- The proof below employs Hubbard trees.
Tao Li \cite{taoli} used the semi-conjugation from the angle doubling map,
see Section~\ref{4ce}. Penrose had obtained a more general statement for
kneading sequences \cite{chprk}. See also Proposition~\ref{Pmonocomb} for
monotonicity with respect to external angles.

\textbf{Proof:} The periodic and preperiodic points marked in the Hubbard tree
$T_c$ of $f_c(z)$ move holomorphically for parameters in the wake of $c$.
The Hubbard tree $T_{c'}$ contains the characteristic point $z_{c'}$
corresponding to $c$ and the connected hull $T\subset T_{c'}$ of its orbit is
homeomorphic to $T_c$\,, but the dynamics of $f_{c'}(z)$ is different in a
neighborhood of $z=0$, which is mapped behind the characteristic point.
There is a forward-invariant Cantor set $C\subset T$ defined by removing
preimages of that neighborhood. To obtain the lower estimate of $h(c')$,
either note that the preimages of a suitable point in $C$ under $f_{c'}(z)$
correspond to those in $T_c$ under $f_c(z)$, or consider a semi-conjugation
$\pi$ from $C\subset T_{c'}$ to $T_c$\,. If $c$ is a center, hyperbolic
arcs in $T_c$ must be collapsed first. \mybox

A parameter $c\in\M$ is a $\beta$-type Misiurewicz point, if
$f_c^k(c)=\beta_c$\,; the minimal $k\ge1$ is the preperiod. The following
results will be needed for Proposition~\ref{Prenpcf}. The weaker estimate
$h(c)\ge\frac{\log2}{k+1}$ is obtained from
Corollary~E in \cite{alfa}; it does not give Proposition~\ref{Prenpcf}.2.

\begin{prop}[$\beta$-type Misiurewicz points]\label{Pbeta}
Suppose $c$ is a $\beta$-type Misiurewicz point of preperiod $k\ge1$.

$1$. The Markov matrix $A$ is irreducible and primitive.

$2$. The largest eigenvalue satisfies $\lambda^k\ge2$,
so $h(c)\ge\frac{\log2}k$\,, with strict inequality for $k\ge2$.

$3$. For the Thurston matrix $F$ according to Proposition~\emph{\ref{Path}},
the arc $[c,\,\beta_c]$ is generating an irreducible primitive block
$B$ of $F$, which corresponds to the largest eigenvalue $\lambda$.
\end{prop}

\textbf{Proof:} Consider the claim that the union
$\bigcup\,f_c^j([0,\,\beta_c])$, $1\le j\le k$, gives the complete Hubbard
tree $T_c$\,: we have $f_c^1([-\alpha_c\,,\,\beta_c])=[\alpha_c\,,\,\beta_c]$
and $f_c^j([0,\,-\alpha_c])$, $1\le j\le k$, provides arcs from $\alpha_c$ to
all endpoints except $\beta_c$\,. When an arc from $\alpha_c$ is crossing
$z=0$, then either its endpoint is behind $-\alpha_c$\, and the part before
$-\alpha_c$ is chopped away before the next iteration. Or this arc is
branching off between $0$ and $-\alpha_c$\,, and further iterates will be
branching from arcs at $\alpha_c$ constructed earlier. This proves the claim.

2. Now the arc $[0,\,\beta_c]$ is containing preimages of itself, so
$f_c^k([0,\,\beta_c])$ is the Hubbard tree $T_c$\,. Since $f_c^k(z)$ is even,
every arc of the Hubbard tree has at least two preimages under $f_c^k(z)$
on the spine $[-\beta_c\,,\,\beta_c]\subset T_c$\,. Each row of $A^k$
has a sum $\ge2$, so $\lambda^k\ge2$.

1. Every edge $e$ of $T_c$ covers $z=0$ in finitely many
iterations. Then it is iterated to an arc at $\beta_c$\,, so the edge $e$
contains an arc iterated to $[0,\,\beta_c]$. Now $f_c^n(e)=T_c$ for
$n\ge n_e$\,, and $A^n$ is strictly positive for
$n\ge\max\{n_e\,|\,e\subset T_c\}$, thus $A$ is irreducible and primitive.
(Alternatively, this follows from Lemma~\ref{Lrenpcf}.4, since $f_c$ is not
renormalizable.) Finally, if $k\ge2$ and we had $\lambda^k=2$, the
characteristic polynomial of $A$ would contain the irreducible factor
$x^k-2$ and $A$ would be imprimitive.

3. Denote the postcritical points by $c_j=f_c^{j-1}(c)$ again. The iteration
of any arc is giving two image arcs frequently, but at least one of these has
an endpoint with index $j$ growing steadily, so reaching $\beta_c=c_{k+1}$\,.
Further iterations produce the arc $[c,\,\beta_c]$ after the other endpoint
was in part $A$ of $\K_c\setminus\{0\}$, which happens when it becomes
$-\beta_c$ or earlier. So the forward-invariant subspace generated by
$[c,\,\beta_c]$ corresponds to an irreducible block $B$ of $F$. It is
primitive by the same arguments as above, since
$[c,\,\beta_c]\mapsto[c,\,\beta_c]+[c,\,c_2]$. Consider the lift
$F_c\,:\,X_c\to X_c$ from the proof of Proposition~\ref{Path} and the
invariant subset $X_c'$ corresponding to $B$. Now $\pi_c\,:\,X_c'\to T_c$
is surjective, since it maps $[c,\,\beta_c]\subset X_c'$ to
$[c,\,\beta_c]\subset T_c$\,, and we have
$\lambda_{\scriptscriptstyle B}=\lambda_{\scriptscriptstyle A}$ again.
Certainly this block $B$ of $F$ need not be equivalent to $A$\,: for
$\theta=3/16$, $A$ is $7\times7$ and $B$ is $6\times6$. And for
$\theta=3/32$, $A$ is $8\times8$ and $B$ is $9\times9$. \mybox

\subsection{Renormalization} \label{3ren}
A filled Julia set $\K_c$ may contain a copy $\K_c^p$ of another
Julia set $\K_{\hat c}$\,, and the Mandelbrot set $\M$ contains
a corresponding copy $\M_p$ of itself. According to Douady and
Hubbard, this phenomenon is explained by renormalization: in a
suitable neighborhood of $\K_c^p$\,, the iterate $f_c^p(z)$ is
quasi-conformally conjugate to $f_{\hat c}(\hat z)$. See
\cite{dhp, daa, hy3, mlcj, sf3}. Many basic results are hard to find
in the literature; a self-contained exposition will be given in \cite{wje}.

For a primitive or satellite center $c_p$ of period $p$ there is a
corresponding tuning map $\M\to\M_p$\,, $\hat c\mapsto c=c_p\ast\hat c$.
This notation suggests that the centers are acting on $\M$ as a semigroup.
Primitive and satellite renormalization are referred to as simple
renormalization, in contrast to crossed renormalization
\cite{mcr, rscr, mlcj}. A small Mandelbrot set $\M_p$ is maximal, if it is
not contained in another small Mandelbrot set; then it is either primitive
maximal or generated by a satellite of the main cardioid. We shall
need two non-standard notations: a \emph{pure satellite component} is not
contained in a primitive Mandelbrot set, but attached to the main cardioid
with a series of satellite bifurcations. And a \emph{maximal-primitive} small
Mandelbrot set is not contained in another primitive one, but it may be
either maximal and primitive or contained in a pure satellite Mandelbrot set.

If $c_p$ is primitive and $\hat c$ is postcritically finite, the
Hubbard tree for $c=c_p\ast\hat c$ is understood in this way: visualize the
$p$-periodic marked points in the Hubbard tree of $c_p$ as small disks,
which are mapped to the next one 1:1 or 2:1, and replace each disk with
a copy of the Hubbard tree of $\hat c$. Actually, the latter Hubbard tree
is extended to include $\pm\beta_{\hat c}$\,. The Markov matrix is obtained
in block form as follows, where $\mathrm{I}$ denotes an identity matrix:
\be\label{renmat}
A =
\left(\begin{array}{c|c}
B & 0\\ \hline
X & R
\end{array}\right)
\qquad\mbox{with}\qquad
R =
\left(\begin{array}{ccccc}
0 & 0 & 0 & \cdots & \hat A\\
\mathrm{I} & 0 & 0 & \ddots & 0\\
0 & \mathrm{I} & 0 & \ddots & 0\\
\vdots & \ddots & \ddots & \ddots & 0\\
0 & 0 & 0 & \mathrm{I} & 0
\end{array}\right)
\ee

\begin{lem}[Renormalization and Markov matrices]\label{Lrenpcf}
$1$. Suppose $c=c_p\ast\hat c$ for a postcritically finite parameter
$\hat c\neq0$. In particular, $c_p$ is a center of period $p\ge2$ and
a suitable restriction of $f_c^p(z)$ is conjugate to $f_{\hat c}(\hat z)$\,.
Then the edges of the Hubbard tree can be labeled such that the Markov matrix
$A$ has the block form $($\ref{renmat}$)$. Here $R$ is imprimitive and
$\hat A$ is the Markov matrix of $\hat c$. 
Moreover:\\
a$)$ If $c_p$ is a primitive center or $\hat c$ is not a $\beta$-type
Misiurewicz point, then $B$ is the Markov matrix of $c_p$\,.\\
b$)$ If $c_p$ is a satellite of $c_p'$ and $\hat c$ is of $\beta$-type,
then $B$ belongs to $c_p'$\,. It is missing completely for immediate
satellite renormalization with $c_p'=0$.

$2$. In all cases of simple renormalization, we have
$h(c)=\max\Big\{h(c_p)\,,\,\frac1p\,h(\hat c)\Big\}$\,.

$3$. Suppose $f_c(z)$ is crossed renormalizable of immediate type and
a suitable restriction of $f_c^p(z)$ is conjugate to $f_{\hat c}(\hat z)$\,.
Then its Markov matrix $A$ is imprimitive and it can be given the form
$A=R$ from $($\ref{renmat}$)$. Here $\hat A$ is the Markov matrix
of $\hat c$, with $\alpha_{\hat c}$ added to the marked points if
necessary, so $h(c)=\frac1p\,h(\hat c)$\,.

$4$. If the Markov matrix $A$ of $f_c(z)$ is reducible or imprimitive,
then $f_c(z)$ is simply or crossed renormalizable with $\hat c\neq0$,
or $c$ is an immediate satellite center.
\end{lem}

In the preperiodic case, $f_c(z)$ is topologically transitive on the Hubbard
tree $T_c$\,, if $A$ is irreducible, and total transitivity corresponds
to a primitive matrix. In this context, results similar to items~1 and~4
have been obtained in \cite{alfa}. However, it is assumed that the small
Hubbard trees are  disjoint only in the case of primitive renormalization.
But $f_c(z)$ is not transitive in the pure satellite case either, except in
the immediate satellite case with $\hat c$ of $\beta$-type.
--- In real dynamics, the relation of imprimitivity and renormalization is
classic. Item~2 is found in \cite{mtimi, dte} for real unimodal maps, in
\cite{taoli} for postcritically finite complex polynomials, and in
\cite{tiob} for arbitrary parameters on veins.

\textbf{Proof of Lemma~\ref{Lrenpcf}:} Note that for $c=c_p$ and $\hat c=0$,
the map $f_c(z)$ is $p$-renormalizable as well, but the small Hubbard tree
is reduced to a point and $R$ is empty. On the other hand, items~1--3 do
not require $p$ to be the maximal non-trivial renormalization period.
Likewise, the proof of item~4 will not produce that period in general.

1. The Julia set of $f_{c_p}(z)$ contains both a superattracting orbit of
period $p$ and a characteristic point of period dividing $p$. By tuning,
the superattracting basins are replaced with a $p$-cycle of small Julia sets,
which are attached to the characteristic point and its images. These sets
are mapped homeomorphically to the next one; only the small Julia set at
$z=0$ is mapped 2:1 to the set at $z=c$. The $I$-blocks in $R$ represent
the homeomorphic restrictions of $f_c^1(z)$, and the diagonal block
$\hat A$ in $R^p$ corresponds to $f_c^p(z)$ and $f_{\hat c}(\hat z)$.
Except in the $\beta$-type case of $\hat c$, the edges of the Hubbard tree
for $f_{c_p}(z)$ are extended into the
small Julia sets, ending at marked points of the small Hubbard trees, which
are before those corresponding $\pm\beta_{\hat c}$\,. This extension does
not happen in the case of a $\beta$-type Misiurewicz point $\hat c$, because
the $\beta$-fixed point of the small Hubbard tree coincides with the
characteristic point corresponding to $c_p$\,. In the satellite case,
this means that the edges from the characteristic point to the
$p$-periodic points are lost completely; only edges from the Hubbard tree
of $c_p'$ are represented in $B$. Note that in the preperiodic case, we must
join the edges at $z=0$ according to Lemma~\ref{LHtree}.1; otherwise the
block $R$ will be more involved.

2. The characteristic polynomials satisfy
$\chi_A(x)=\chi_B(x)\cdot\chi_R(x)=\chi_B(x)\cdot\chi_{\hat A}(x^p)$. For
the second equality, consider a complex eigenvector in block form to see that
$x$ is an eigenvalue of $R$ if and only if $x^p$ is an eigenvalue of $\hat A$.
Now the largest eigenvalue of $B$ is related to $h(c_p)$ and the largest
eigenvalue of $R$ corresponds to $\frac1p\,h(\hat c)$. When $B$ represents
$c_p'$ instead of $c_p$\,, first apply the same argument to show
$h(c_p')=h(c_p)$. 

3. The $p$-cycle of small Julia sets is mapped $p-1$ times homeomorphically
and once like $f_c^p(z)$ and $f_{\hat c}(\hat z)$ as above. Since they are
crossing at the fixed point $\alpha_c$\,, $\alpha_{\hat c}$ is marked in the
Hubbard tree determining the block $\hat A$. If $\hat c$ belongs to the
$1/2$-limb and is not of $\alpha$-type, this point would not be marked in the
minimal Hubbard tree of $f_{\hat c}(\hat z)$,  but doing so does not change
the largest eigenvalue according to Lemma~\ref{LHtree}.2. There are no further
edges in the Hubbard tree of $f_c(z)$.

4. Since $f_c(z)$ is expanding on the Hubbard tree $T_c$\,, every edge will
cover the edge $e_c$ before $z=c$ under the iteration. So the images
$f_c^j(e_c)$ form an absorbing invariant family of edges, and there is a
corresponding irreducible block $R$ in $A$. If $A$ is reducible, $R$ is not
all of $A$, and the corresponding subset of $T_c$ has $p\ge2$ connected
components: if it was connected, it would be all of $T_c$\,. By surjectivity,
each of these components is mapped onto another one. Each component is
forming a small tree, and collapsing all components to points
gives a $p$-periodic tree.

If $A$ is irreducible and imprimitive of index $p$, it has a Frobenius
normal form with a block structure similar to $R$ in (\ref{renmat}), except
the non-0 blocks need not be quadratic. So the family of edges is a disjoint
union of $E_1\,,\,\dots\,,\,E_p$ with $f_c:E_j\to E_{j+1}\,,\,E_p\to E_1$\,.
There is a corresponding subdivision of the Hubbard tree into $p$ subsets;
these are disjoint except for common vertices, and we do not need to show
that they are connected. Since every subset is mapped into itself under $kp$
iterations only, every edge has return numbers divisible by $p$. Every
characteristic periodic point has a ray period divisible by $p$. If $c$ is
not satellite renormalizable, the preperiodic critical value $c$ or the
primitive characteristic point $z_1$ is approached from below by primitive
characteristic points $x_n$ of minimal periods, which are increasing and
divisible by $p$. Thus for large $n$, the internal address of $x_n$ will be
ending on $\dots$--$\mathrm{Per}(x_{n-1})$--$\mathrm{Per}(x_n)$, so $x_n$ is
$p$-renormalizable according to \cite{bks}. The same applies to $c$,
since that renormalization locus is closed. \mybox

\begin{rmk}[Same entropy for different Hubbard trees] \label{Rsedht}
1. For any postcritically finite parameter not on the real axis, there are
corresponding parameters with a homeomorphic Hubbard tree,
but different rotation numbers at periodic branch points \cite{bks}.
These parameters have the same core entropy: neither the Markov matrix $A$
nor the Thurston matrix $F$ according to Proposition~\ref{Path} depends
on the rotation numbers. The dynamics of these parameters are conjugate
by homeomorphisms, which are not orientation-preserving \cite{bks}.

2. Suppose $p,\,q\ge2$ and $\hat c$ is a postcritically finite parameter
in the $1/q$-limb. Consider the following parameters, whose
$p$-renormalization is conjugate to $f_{\hat c}(\hat z)$\,: $c$ is satellite
renormalizable in the $1/p$-limb, and $c'$ is immediately crossed
renormalizable in the $\frac1{pq}$-limb. Then $h(c)=h(c')$ although the
dynamics are not conjugate globally. If $q\ge3$ and $\hat c$ is of
$\beta$-type, the Markov matrices for the minimal Hubbard trees are
identical in fact.

3. The Misiurewicz parameter $c=\gamma_\sM(3/14)$ has a few remarkable
properties. The characteristic polynomial of $F$ is
$P(x)=x^4-3x-2=(x^2-x-1)\cdot(x^2+x+2)$ and $A$ has $x\cdot P(x)$.
This factorization is non-trivial, i.e., not involving roots of unity.
Moreover, the first factor agrees for the primitive Mandelbrot set of
period 3: the core entropy is $h(c_3)=h(c)$ but the dynamics are unrelated.
$\K_c$ is discussed in \cite{bmplm, ersqs} in a different context: $f_c(z)$
on $\K_c$ is quasi-conformally conjugate to a piecewise-affine map $\pm sx-1$
with $s^3-s^2+s+2=0$. 
\end{rmk}

\begin{prop}[Renormalization and entropy]\label{Prenpcf}
Suppose $c=c_p\ast\hat c$ for a postcritically finite parameter
$\hat c$.

$1$. If the center $c_p$ is of pure satellite type, we have
$h(c_p)=0$ and $h(c)=\frac1p\,h(\hat c)$\,.

$2$. If $c_p$ is a primitive center and $c=c_p\ast\hat c$, then
$h(c)=h(c_p)>\frac{\log2}p$\,.

$3$. Consider a maximal-primitive small Mandelbrot set $\M_p=c_p\ast\M$
as defined above, and denote its tip by $c_p'=c_p\ast(-2)$. There is a
sequence of Misiurewicz points $a_n\to c_p'$ on the vein behind $c_p'$\,,
such that $h(a_n)>h(c_p')=h(c_p)$\,.
\end{prop}

Item~2 means that $h(c)$ is constant on primitive Mandelbrot sets,
see also Theorem~\ref{Tmonoren}. This fact will simplify various proofs
in Section~\ref{4}. Item~3 shows that the monotonicity
from Proposition~\ref{Pmonopcf} is strict, when $c$ and $c'$ are separated
by a maximal-primitive Mandelbrot set. According to \cite{sf2}, two
postcritically finite parameters are separated by a root. Now there is a
maximal-primitive root as well, unless both parameters are of pure satellite
type or belong to the same primitive Mandelbrot set.

\textbf{Proof of Proposition~\ref{Prenpcf}:}
1. Suppose $c$ is an immediate satellite center of the main cardioid, then
the Hubbard tree is a star with an endpoint at $z=0$, so $f_c(z)$ is
injective on it. Moreover, the characteristic polynomial of $A$ is $x^p-1$
when $p\ge3$. Both arguments give $h(c)=0$. For a pure satellite center
$c_p=c_k\ast\dots\ast c_1$\,, apply Lemma~\ref{Lrenpcf}.2
recursively to show $h(c_p)=0$. The same relation gives
$h(c_p\ast\hat c)=\frac1p\,h(\hat c)$\,.

2. Suppose that $c_p$ is not pure satellite renormalizable. There is
an immediate satellite $c_q$ of period $q$ and a $\beta$-type Misiurewicz
point $b_k$ of preperiod $k$, such that $c_p$ is behind the $\alpha$-type
Misiurewicz point $c_q\ast b_k$\,. By the estimate of lowest periods in
decorations \cite{wje} we have $p\ge(k-1)q+(q+1)=kq+1$. Monotonicity and
the estimate from Proposition~\ref{Pbeta}.2 give $h(c_p)\ge h(c_q\ast b_k)
=\frac1q\,h(b_k)\ge\frac1q\,\frac{\log2}k>\frac{\log2}p$\,. More
generally, if $c_p=c_s\ast c_q$ with $c_s$ pure satellite type of period $s$
and $c_q$ of period $q$ primitive and not pure satellite renormalizable,
then $h(c_p)=\frac1s\,h(c_q)>\frac1s\,\frac{\log2}q=\frac{\log2}p$ again.
Finally, for any postcritically finite parameter $\hat c$ we have
$h(c_p\ast\hat c)=\max\Big\{h(c_p)\,,\,\frac1p\,h(\hat c)\Big\}=h(c_p)$\,,
since $\frac1p\,h(\hat c)\le\frac{\log2}p<h(c_p)$\,.

3. Suppose that $c_p$ is primitive and maximal, not pure satellite
renormalizable. For $c=a_n$\,, the critical value is behind the disconnected
small Julia set. After $p$ iterations it will be before the characteristic
point, following its orbit but moving farther away with every $p$-cycle of
iterations. After $pn\pm k$ iterations it will join a repelling cycle
independent of $n$. Misiurewicz points with these properties are constructed
to describe the domains of renormalization \cite{wje}; in the maximal case
we may use $\alpha$-type Misiurewicz points. Fixing any large $n$, set
$c=a_{n+1}$ and $c'=a_n\succ c\succ c_p'$\,. We shall see that $h(c')>h(c)$.

The relevant preperiodic points exist in both Julia sets, as does the
$p$-cycle of disconnected small Julia sets. The edges shall be corresponding
as follows: first, the edges leading to the postcritical points behind the
small Julia sets in $T_c$ and $T_{c'}$ are identified. Second, preperiodic
marked points of $T_c$ between the small Julia sets are marked in $T_{c'}$
in addition to the postcritical points; by Lemma~\ref{LHtree}.3 this does not
change the largest eigenvalue $\lambda'$ of $A'$.
Due to the identifications, the only difference of $A$ and $A'$ is related
to a small edge at the first interior postcritical point: it has one
preimage under $A$ and two preimages under $A'$. Since $A$ is primitive by
Lemma~\ref{Lrenpcf}.4, and $A'\ge A$ in each component, we have
$\lambda'>\lambda$ according to Section~\ref{2pf}. Finally, if $c_p$ is
maximal-primitive but not primitive maximal, the Misiurewicz points
will belong to the same pure satellite Mandelbrot set $c_s\ast\M$, and
the inequalities for $\hat a_n$ are transferred to $a_n=c_s\ast\hat a_n$
according to $h(a_n)=\frac1s\,h(\hat a_n)$. \mybox

Item~2 and Lemma~\ref{Lrenpcf} show that a postcritically finite parameter
$c\neq0$ has the following property, if and only if it is neither pure
satellite renormalizable nor immediately crossed renormalizable:
the largest eigenvalue $\lambda$ of $A$ is simple, and there is no other
eigenvalue of the same modulus. --- Suppose $c_p$ has the internal
address $1$--$\dots$--$s$--$q$--$\dots$--$p$, where $s$ is of pure satellite
type and $q$ is primitive, then item~2 is strengthened by monotonicity to
$h(c_p)>\frac{\log2}q$\,. The estimate $h(c_p)\ge\frac{\log2}q$ is due to
Bruin--Schleicher \cite{bshb}; 
equality can be ruled out by Lemma~\ref{Lrenpcf}.4 as well, since $c_q$
is not $q$-renormalizable with $\hat c\neq0$.

\subsection{Biaccessibility dimension of postcritically finite maps}
\label{3biacpcf}
When $c$ is a Misiurewicz point, the filled Julia set $\K_c$ is a countable
union of arcs between periodic and preperiodic points, plus an 
uncountable union of accumulation points. The interior points of the arcs are
pinching points, so they are biaccessible: there are at least two dynamic rays
landing. Each of these points is iterated to the Hubbard tree in a finite
number of steps. When $c$ is a center, the arcs of the Hubbard tree and
their preimages are meeting a dense family of Fatou basins each. For
a pure satellite center, these basins have a countable family of common
boundary points and accumulation points on an arc, so the biaccessible points
in $\partial\K_c$ are countable. When $c$ is a primitive center, or a
satellite center within a primitive Mandelbrot set, the primitive small Julia
sets do not have common boundary points, and there is an uncountable family of
biaccessible accumulation points on every arc. In any case, a countable union
of linear preimages shows that the Hausdorff dimension of biaccessing angles
of rays is the same for the Hubbard tree $T_c$ and for $\partial\K_c$\,. It is
related to the core entropy $h(c)$ as follows:

\begin{prop}[Hausdorff dimension and Hausdorff measure]\label{Pbhpcf}
$1$. Suppose $f_c(z)$ is postcritically finite and consider the subset
$\K_c'\subset\K_c$ of biaccessible points. The Hausdorff dimension of
biaccessing angles is related to the core entropy by
$b:=\dim\gamma_c^{-1}(\K_c')=\dim\gamma_c^{-1}(T_c)=h(c)/\log2$
according to Thurston \emph{\cite{mesh, gaotl, bshb, tiob}}.

$2$. When the biaccessibility dimension is $0<b\le1$, the Hausdorff
measure of the external angles of the Hubbard tree satisfies
$0<\mu_b(\gamma_c^{-1}(T_c))<\infty$. Taking all biaccessing angles of
$\K_c$\,, we have $\mu_b(\gamma_c^{-1}(\K_c'))=\infty$ when $0<b<1$.
\end{prop}

The Thurston relation $h(c)=b\cdot\log2$ will be proved in Section~\ref{4ce}
under more general assumptions, using the fact that
$\gamma_c\,:\,\gamma_c^{-1}(T_c)\subset S^1\to T_c \cap\partial\K_c$ is a
semi-conjugation from the angle doubling map $F(\phi)$ to $f_c(z)$\,;
this proof from \cite{bshb, tiob} has no essential simplification in the
postcritically finite case. The following alternative argument is based on
Hausdorff measure, and it gives a partial answer to a
question in \cite[Remark~1.1]{mesh}.

\textbf{Proof:} There is a finite number of intervals in $S^1$, such that
$\gamma_c^{-1}(T_c)$ is obtained by removing these intervals and their
preimages recursively \cite{dte, taoli, tiob}. Since $f_c(z)$ is even and
$f_c:T_c\to T_c$ is
surjective but not $2:1$ globally (unless $c=-2$), these intervals correspond
to the components of $(-T_c)\setminus T_c$\,: if a component is attached at
a pinching point in $\partial\K_c$\,, the interval of angles is obvious.
If the cut point of a component is a superattracting periodic point, the
interval is determined from the repelling periodic point on the boundary of
the Fatou basin. E.g., for the primitive center $c$ of period $4$ in the
$1/3$-limb, the initial intervals are $(9/14,\,11/14)$ and $(12/15,\,1/15)$
in cyclic order. 

Removing preimages up to order $n$ can be described in terms of an
equidistant subdivision: there are $x_n$ closed intervals $W_j$ of length
$C2^{-n}$ left, where $C$ is related to the common denominator of the
initial intervals. We shall see that $x_n$ is bounded above and below as
$x_n\asymp\lambda^n$ when $h(c)=\log\lambda>0$, and obtain the Hausdorff
dimension from $\lambda$. The basic idea is to consider preimages of the cut
points in $T_c$\,, since preimages in other parts of $\K_c$ are removed
together with another branch. After a few initial steps, the preimages are
preperiodic and their number equals the number of preimages of certain edges.
The decomposition of $A$ from Lemma~\ref{Lrenpcf} shows that this number is
growing by $\lambda$; although certain edges have fewer preimages in the case
of pure satellite renormalization, these do not contain initial cut points.
Now the remaining intervals have total length $x_nC2^{-n}$ given by the tails
of geometric series (or their derivatives) in the halved eigenvalues of $A$.
So $x_n\sim K\lambda^n$ asymptotically when $\lambda$ corresponds to a
primitive block of $A$ and $x_n\asymp\lambda^n$ in general, since $\lambda>1$
when pure satellite centers are excluded.

The equidistant covering of $\gamma_c^{-1}(T_c)$ by $x_n$ intervals $W_j$ of
length $C2^{-n}$ gives the box dimension $b=\log\lambda/\log2=h(c)/\log2$ and
an upper estimate $\mu_b(\gamma_c^{-1}(T_c))\le\limsup x_n(C2^{-n})^b<\infty$
when $b>0$ according to the definition (\ref{eqhdm}). Since
$\gamma_c^{-1}(T_c)$ is closed and invariant under doubling $F(\phi)=2\phi$,
the Hausdorff dimension $\dim\gamma_c^{-1}(T_c)=b$ is obtained from
Furstenberg \cite[Proposition~III.1]{fu}. We shall use similar arguments to
show $\mu_b(\gamma_c^{-1}(T_c))>0$, following the exposition by Gao
\cite[Section~5.2.1]{gao}. Using an equidistant subdivision $W_j'$ of length
$2^{-n}$ still gives $\sum|W_j'|^b\asymp\lambda^n2^{-bn}\asymp1$. Assuming
$\mu_b(\gamma_c^{-1}(T_c))=0$, there is a finite cover by closed intervals
with $|U_i|=2^{-n_i}$, aligned to the dyadic grids of different levels
$n_i$\,, such that $\sum|U_i|^b<1$. Define puzzle pieces by taking preimages
under $F^{n_i}$ restricted to $U_i$\,: by $F$-invariance there are nested
puzzle pieces around the points of $\gamma_c^{-1}(T_c)$. The pieces of depth
$k$ contribute $(\sum|U_i|^b)^k$ to (\ref{eqhdm}). For large
$n$, each $W_j'$ will be contained in a puzzle piece of large depth $k$, but
our estimates do not give a contradiction yet: there may be several $W_j'$
in the same piece. So, set $N=\max n_i$ and choose a puzzle piece of maximal
level $\sum n_i$ for each $W_j'$\,. By induction, the level will be $\ge n-N$,
since images under $F^{n_i}$ still have maximal levels. Now there are at most
$2^N$ of the $W_j'$ contained in one piece, the depth is $k\ge n/N-1$, and
\be 1 \asymp \sum_j|W_j'|^b \le
 2^N\sum_{k\sim n/N}^\infty\Big(\sum_i|U_i|^b\Big)^k \to 0
 \quad\mbox{as}\quad n\to\infty \ . \ee
This is a contradiction disproving the assumption
$\mu_b(\gamma_c^{-1}(T_c))=0$.
When $b<1$, an arc in $(-T_c)\setminus T_c\neq\emptyset$ has disjoint
preimages. In each step, their number is doubled and their individual measure
is divided by $2^b<2$, giving $\mu_b(\gamma_c^{-1}(\K_c'))=\infty$. \mybox

For the principal $\alpha$-type Misiurewicz point $c$ in the $p/q$-limb,
the angles of an edge of $T_c$ are obtained explicitly: on each side, we
have a standard Cantor set, where an interval is divided into $2^q$ pieces
recursively and the outer two intervals are kept. So for $b=1/q$, the sum
in (\ref{eqhdm}) is independent of the levels of subdivision.

\section{The biaccessibility dimension} \label{4}
\begin{figure}[h!t!b!]
\unitlength 0.001\textwidth 
\begin{picture}(990, 220)
\put(6, 0){\includegraphics[width=0.22\textwidth]{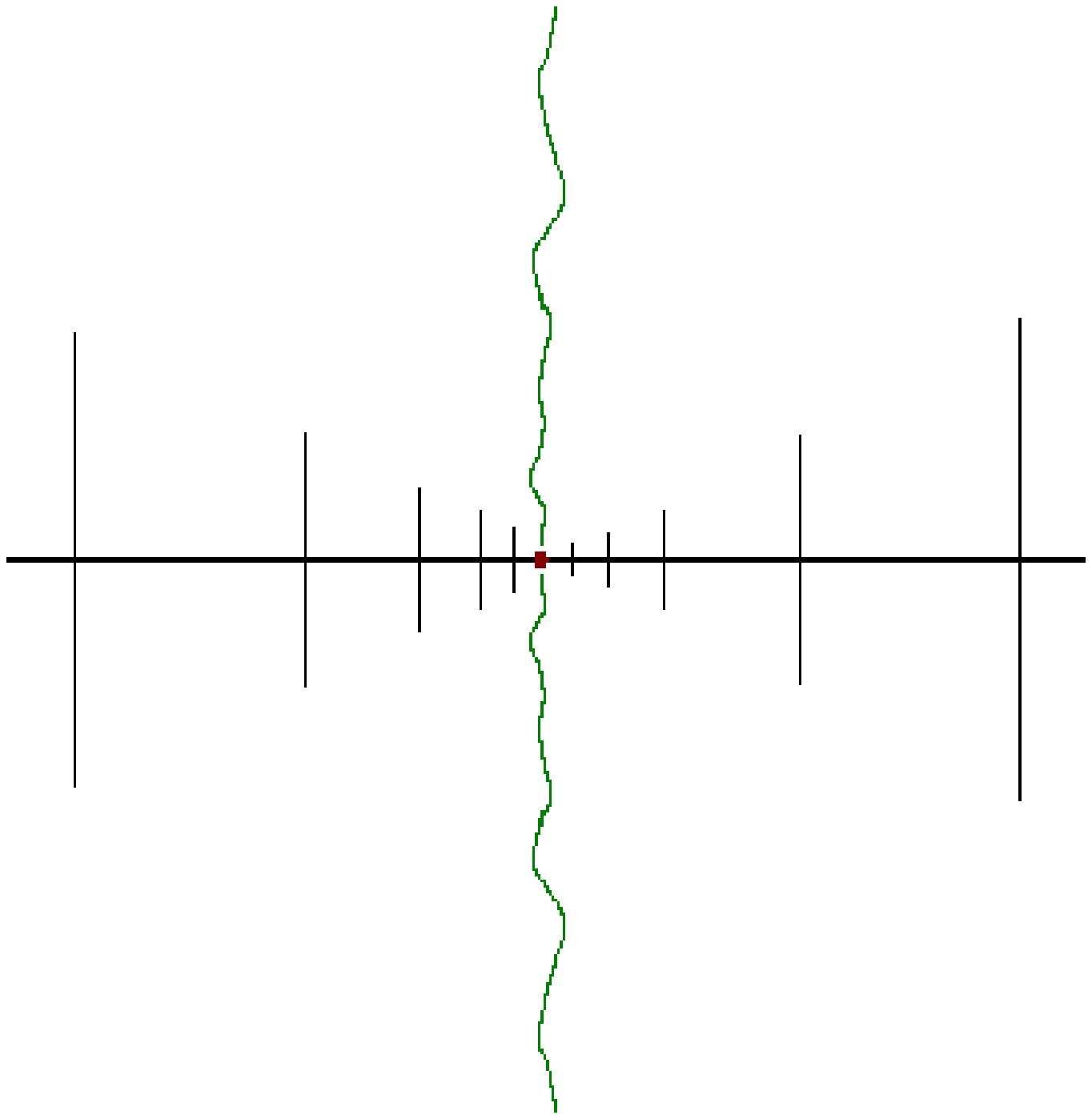}}
\put(262, 0){\includegraphics[width=0.22\textwidth]{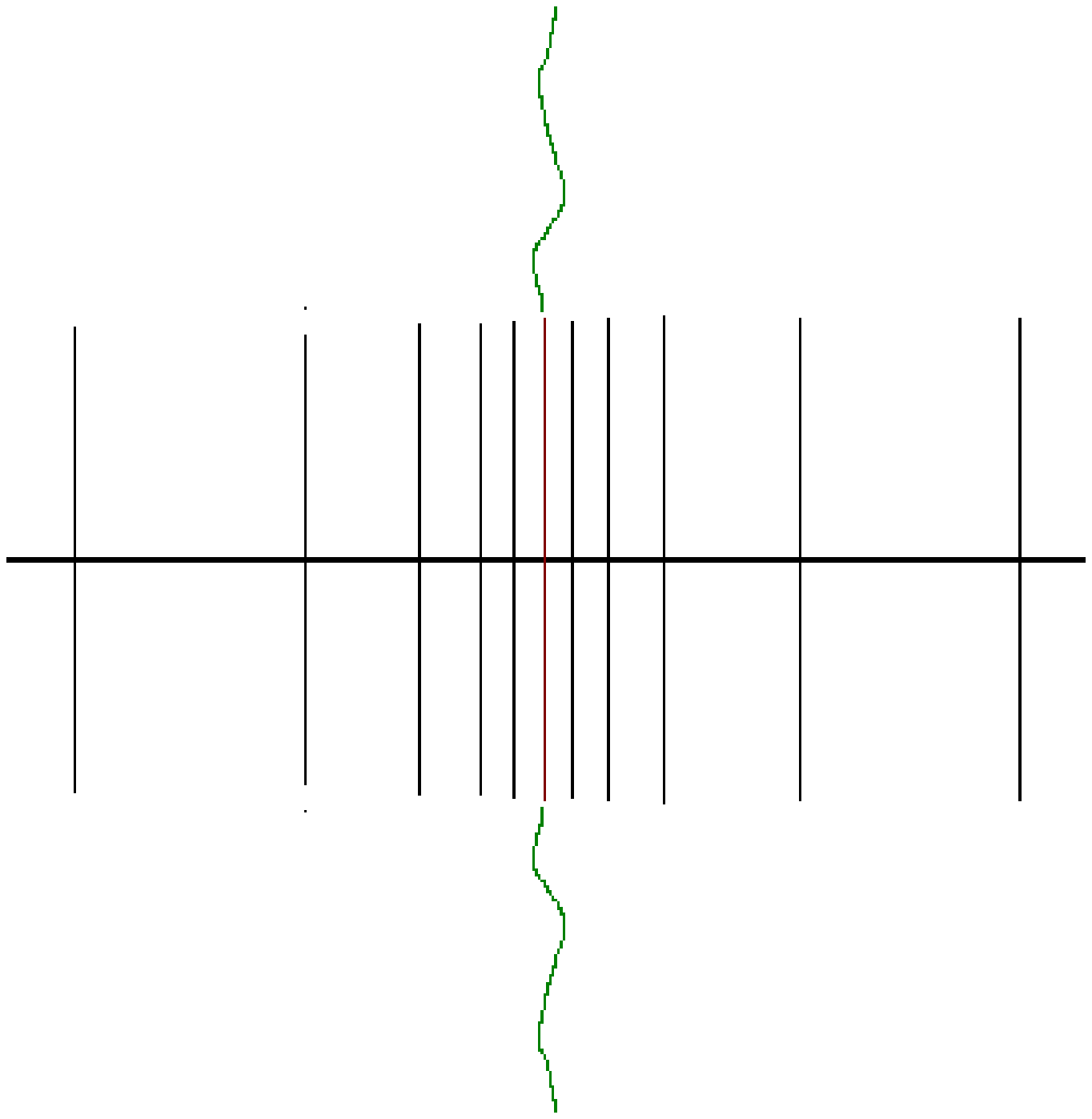}}
\put(518, 0){\includegraphics[width=0.22\textwidth]{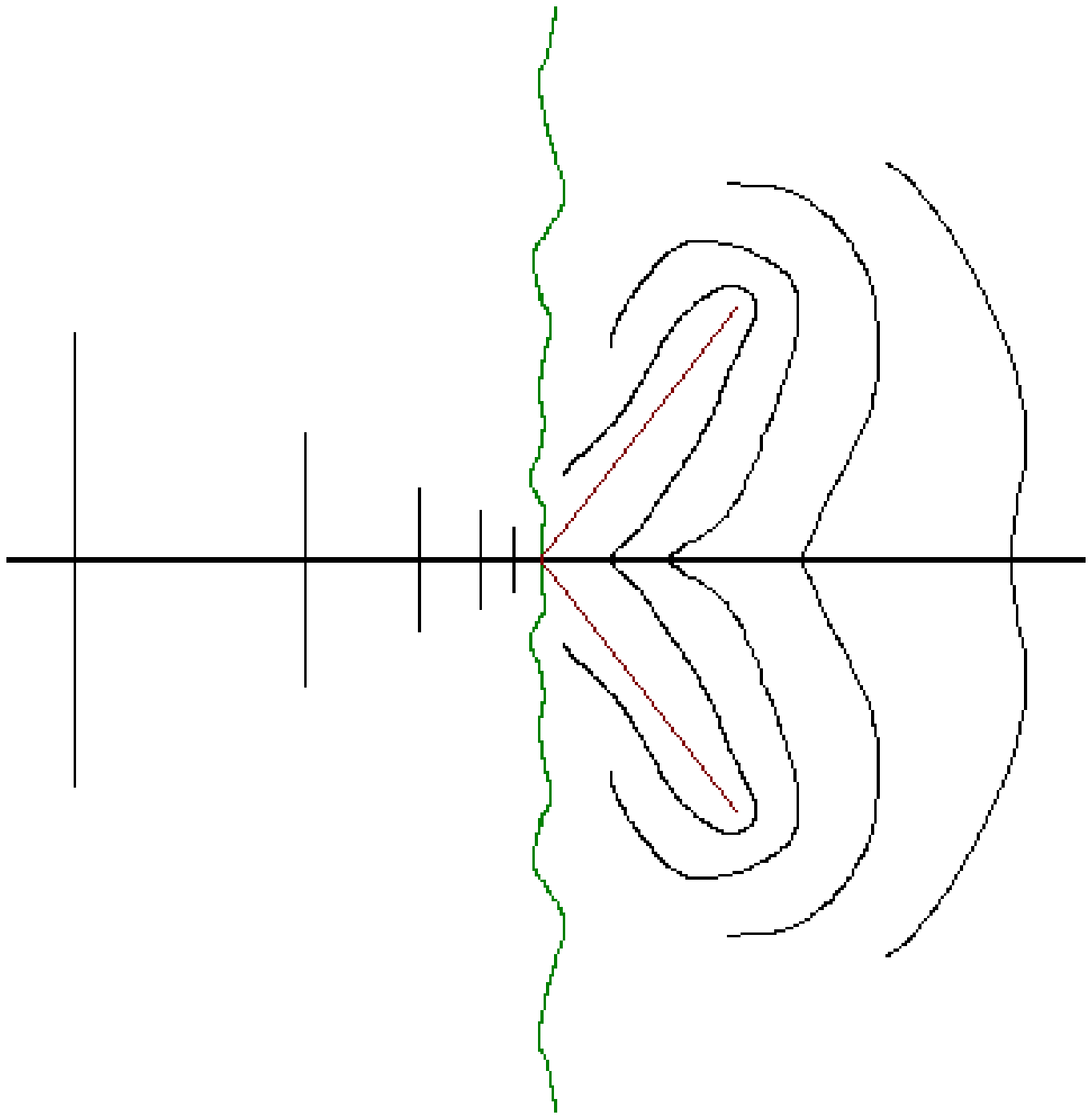}}
\put(774, 0){\includegraphics[width=0.22\textwidth]{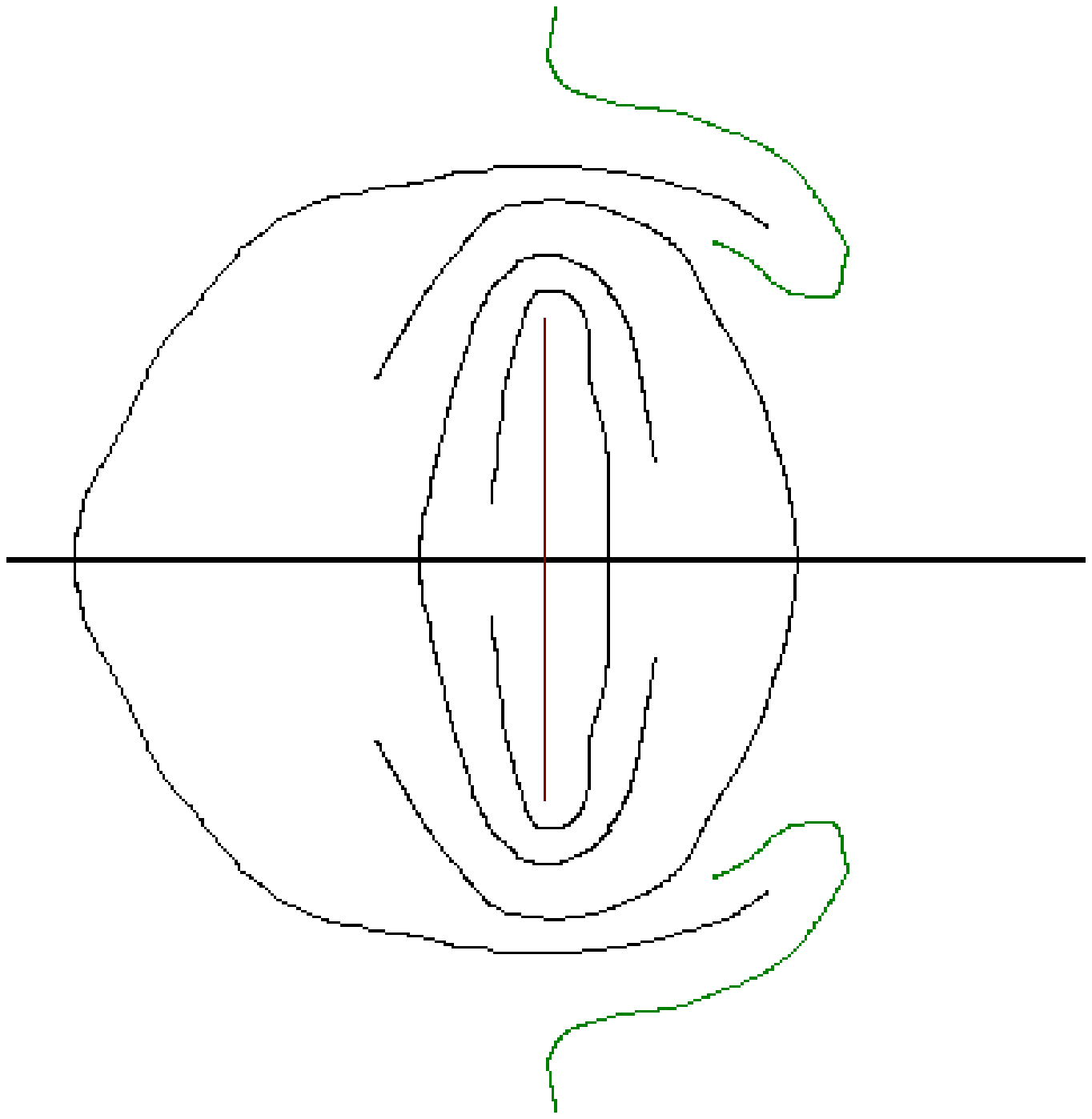}}
\put(7, 219){\makebox(0, 0)[lt]{a)}}
\put(263, 219){\makebox(0, 0)[lt]{b)}}
\put(519, 219){\makebox(0, 0)[lt]{c)}}
\put(775, 219){\makebox(0, 0)[lt]{d)}}
\end{picture} \caption[]{\label{FFB}
Possible scenarios at the Feigenbaum parameter
$c_{\mbox{\tiny F}}\in\M\cap\R$: in a) $\M$ is openly locally connected at
$c_{\mbox{\tiny F}}$\,, and the fiber of $c_{\mbox{\tiny F}}$ is
non-trivial in the other cases \cite{dcc, sf2}. It is
formed from the impressions of two rays. In case d) the rays do not land;
they have non-trivial accumulation sets. This example could be modified, such
that the impressions are larger than the accumulation sets. 
Note that in cases b) and d), $\M\setminus\{c_{\mbox{\tiny F}}\}$ is not
disconnected.\\
Similar situations arise for combinatorially biaccessing dynamic rays of a
filled Julia set $\K_c$\,, but they need not be symmetric. ---
Non-trivial fibers may correspond to a single external ray as well.}
\end{figure}

A filled Julia set $\K_c$\,, or the Mandelbrot set $\M$, is locally connected
if and only if every external ray lands and the landing point depends
continuously on the angle (Carath\'eodory). Some possible scenarios are shown
in Figure~\ref{FFB}. It may happen that a pair of rays is approximated by
pairs of rays landing together, without doing the same. Therefore we must
distinguish between combinatorial biaccessibility and topological
biaccessibility. However, every pinching point $z\in\K_c$\,, i.e.,
$\K_c\setminus\{z\}$ is disconnected, must be the landing point of at least
two rays \cite[p.~85]{mcr}. Non-local connectivity of $\M$ would mean that
$\M$ contains non-trivial fibers; these compact connected subsets are
preimages of points under the projection $\pi_\sM:\partial\M\to S^1\!/\!\sim$
of \cite{dcc} described in Section~\ref{4cc}. As an alternative
characterization, these subsets cannot be disconnected by pinching points
with rational angles \cite{sf2} (and they are disjoint from closed
hyperbolic components). There is an analogous concept for Julia sets, but
rational angles do not suffice to describe fibers of Julia sets for
Siegel and Cremer parameters \cite{sf3}.

\subsection{Combinatorial and topological biaccessibility} \label{4ctb}
For a given angle $\theta$, the itinerary of $\phi$ is a sequence of symbols
$A,\,B,\,\ast$ or $1,\,0,\,\ast$ describing the orbit of $\phi$ under
doubling with respect to a partition of the circle, which is defined by the
diameter joining $\theta/2$ and $(\theta+1)/2$\,. See \cite{hss, ser, bks}.
The following result from \cite[Section~9]{bks} motivates the definition
of combinatorial biaccessibility by itineraries:

\begin{lem}[Rays landing together, Bruin--Schleicher]\label{Lland}
Assume that $\K_c$ is locally connected and for some $\theta\in S^1$,
$\r_c(\theta)$ lands at the critical value $c$, or $\theta$ is a
characteristic angle of $z_1$ at the Fatou basin around $c$.\\
Suppose the rays $\r_c(\phi_1)$ and $\r_c(\phi_2)$ do not land at
preimages of $c$ or $z_1$\,. Then they land together, if and only if
the itineraries of $\phi_1$ and $\phi_2$ with respect to $\theta$ are equal.
\end{lem}

\textbf{Proof:} A separation line is formed by the rays $\r_c(\theta/2)$ and
$\r_c((\theta+1)/2)$ together with the common landing point $z=0$; in the
parabolic or hyperbolic case, there are two landing points connected with
an interior arc. The itinerary of an angle with respect to
$\theta$ is equal to the itinerary of the landing point with respect to the
separation line, since preimages of $c$ or $z_1$ are excluded. So if the
landing point is the same, the itineraries coincide. Now suppose the
landing points $w_i$ of $\r_c(\phi_i)$ are different: we need to show that
$w_1$ and $w_2$ are separated by a preimage of $z=0$. In the case of empty
interior, the arc $[w_1\,,\,w_2]$ is containing two biaccessible points
$w_1'\,,\,w_2'$\,, which satisfy
$f_c^n([w_1'\,,\,w_2'])\subset[-\beta_c\,,\,\beta_c]$. The subarc covers
$0$ in finitely many steps, since otherwise the corresponding intervals of
angles could be doubled indefinitely. If the interior of $\K_c$ is not empty,
the arc is not defined uniquely in interior components, but an image of the
subarc must cross the component around $z=0$ for the same reason. \mybox

So the biaccessibility dimension will be the same
for all external angles of $c$. If $\K_c$ is not locally connected, there may
be pairs of rays approximated by pairs landing together, but not doing so
themselves; this may look like the parameter space picture b) or d) in
Figure~\ref{FFB}, or be more complicated. With countably many exceptions,
the criterion in terms of itineraries shows, which rays should land
together according to the combinatorics. When $\theta$ is not periodic, for
any word $w$ angles with the itineraries $w*\nu$, $wA\nu$, and $wB\nu$ have
rays landing together as well. Here $\nu$ denotes the kneading sequence,
i.e., the itinerary of $\theta$ with respect to itself. If $\theta$ is
periodic, you cannot see from the itinerary alone whether two angles belong
to the same preimage of $z_1$\,. Again, the ambiguity concerns a countable
number of angles only; it does not matter for the Hausdorff dimension:

\begin{dfn}[Biaccessibility dimension]\label{Dbiacdim}
$1$. With respect to an angle $\theta\in S^1$, $\phi_1$ is
\emph{combinatorially biaccessing}, if there is an angle $\phi_2\neq\phi_1$
such that $\phi_1$ and $\phi_2$ are connected by a leaf of the lamination
associated to $\theta$ $\cite{tgdr}$. With countably many exceptions, this
means that $\phi_1$ and $\phi_2$ have the same itinerary, which does not
contain an $\ast$. The Hausdorff dimension of these angles is the
\emph{combinatorial biaccessibility dimension} $B_\comb(\theta)$.

$2$. For $c\in\M$, an angle $\phi_1\in S^1$ is
\emph{topologically biaccessing}, if there is an angle $\phi_2\neq\phi_1$
such that $\r_c(\phi_1)$ and $\r_c(\phi_2)$ land at the same
pinching point $z\in\K_c$\,. The Hausdorff dimension of these angles is
the \emph{topological biaccessibility dimension} $B_\top(c)$.
\end{dfn}

The following result is basically due to Bruin and Schleicher; see the
partial results in \cite[Section~9]{bks} and the announcement in \cite{bshb}.
The proof below will discuss the dynamics of various cases; some arguments
could be replaced by the geometric observation, that the angles of non-landing
rays have Hausdorff dimension 0, see the references in \cite{bshb}.

\begin{thm}[Combinatorial and topological biaccessibility]\label{Tcombtop}
Suppose that $\theta\in S^1$ and $c\in\partial\M$ belongs to the impression
of $\r_\sM(\theta)$, or $c\in\M$ is hyperbolic and the ray lands at the
corresponding root. Then $B_\comb(\theta)=B_\top(c)$.
\end{thm}

Recall the notions of maximal-primitive renormalization and of pure
satellite components from Section~\ref{3ren}. When $f_c(z)$ is simply
$p$-renormalizable, there is a small Julia set $\K_c^p\subset\K_c$ and
$f_c^p(z)$ is conjugate to $f_{\hat c}(\hat z)$ in a neighborhood; the
hybrid-equivalence $\psi_c(z)$ is mapping $\K_c^p\to\K_{\hat c}$\,. Now
$\K_c\setminus\K_c^p$ consists of a countable family of decorations, which
are attached at preimages of the small $\beta$-fixed point. Each decoration
corresponds to an open interval of angles, and their recursive removal shows
that the angles of rays accumulating or landing at $\partial\K_c^p$ form a
Cantor set of Hausdorff dimension $1/p$ and finite positive
Hausdorff measure \cite{mlr, wje}.
The parameter $c$ belongs to a small Mandelbrot set $\M_p\subset\M$, which
has corresponding decorations and the same Cantor set of angles, such that
the parameter ray lands or accumulates at $\partial\M_p$\,. The Douady
substitution of binary digits \cite{daa, mlr} relates the external angles
$\phi$ and $\hat\phi$: if $\r_c(\phi)$ lands at $z\in\partial\K_c$\,,
then $\psi_c(\r_c(\phi))$ lands at
$\hat z=\psi_c(z)\in\partial\K_{\hat c}$ and according to Liouville,
$\r_{\hat c}(\hat\phi)$ lands at $\hat z$ as well. The converse statement
is not obvious, because $\psi_c^{-1}(\r_{\hat c}(\hat\phi))$
may cut through the decorations. But when the ``tubing'' is chosen
appropriately, $\r_{\hat c}(\hat\phi)$ can be deformed such that it
is avoiding the images of decorations, still landing at $\hat z$.
The same arguments apply to parameter rays \cite{wje}.

\textbf{Proof of Theorem~\ref{Tcombtop}:}
The problem is that $\K_c$ may be non-locally connected; then two rays with
the same kneading sequence may land at different points of a
non-trivial fiber, or one may accumulate there without landing.
Cf.~Figure~\ref{FFB}. 
By the Yoccoz Theorem~\cite{mlcj, hy3}, non-local connectivity can
happen only for infinitely simply renormalizable parameters $c$, or for
$f_c(z)$ with a neutral cycle. Consider the following cases:

1. If $f_c(z)$ is hyperbolic or parabolic, $\K_c$ is locally connected and
$\theta$ is an external angle of the characteristic point \cite{mer, ser}.
If $f_c(z)$ is non-renormalizable or finitely renormalizable with all cycles
repelling, $\K_c$ is locally connected and $\r_c(\theta)$ lands at the
critical value: both for the dynamic plane and the parameter plane, $\theta$
is approximated by corresponding rational angles \cite{sf2, sf3, ser}. The
statement follows from Lemma~\ref{Lland}.

2. If $\K_c$ contains a Siegel disk or a Cremer point of period 1, we
have $B_\top(c)=0$ \cite{zb2}. If $c$ is a Siegel or Cremer parameter
on the boundary of a pure satellite component, there are countably
many biaccessing angles from the satellite bifurcations. Each copy of the
small Julia set has biaccessing angles of Hausdorff dimension 0, since the
renormalized rays would land at $\K_{\hat c}$\,. On the other hand,
$B_\comb(\theta)=0$ according to \cite{bsoos, bshb}. For an irrationally
neutral parameter in a primitive Mandelbrot set, the biaccessibility
dimension will be positive; these parameters are included in case~4.

3. Suppose the parameter $c$ is infinitely pure satellite renormalizable.
Choose a pure satellite component of period $m$ before $c$ and consider
the angles of the rays accumulating or landing at the small Julia set and
its preimages. Their Hausdorff dimension is $\frac1m$\,; so
$B_\top(c)\le\frac1m$ and $B_\comb(\theta)\le\frac1m$\,,
and letting $m\to\infty$ gives 0. All other infinitely renormalizable
parameters are contained in a primitive Mandelbrot set.

4. Suppose $c=c_p\ast\hat c$ for a primitive center $c_p$ of period $p$, and
denote the angles of the corresponding root by $\theta_\pm$\,. Since
$\r_\sM(\theta)$ lands or accumulates at $\M_p=c_p\ast\M$, $\theta$
belongs to the Cantor set of angles constructed above and $\r_c(\theta)$
accumulates or lands at $\K_c^p$\,. We will not need to check, whether
it accumulates at the critical value $z=c$. By item~1,
Proposition~\ref{Prenpcf}.2,
and the relation to core entropy from Proposition~\ref{Pbhpcf}.1,
we have $B_\top(c_p)=B_\comb(\theta_\pm)>\frac1p$\,. But $\frac1p$ is the
Hausdorff dimension of those angles, whose rays land or accumulate
at the small Julia sets, so these are negligible. Every biaccessible point
$z\in\K_{c_p}$ is iterated to the spine $[-\beta_{c_p}\,,\,\beta_{c_p}]$,
and the angles of the spine alone have Hausdorff dimension $B_\top(c_p)$
already. This arc consists of a Cantor set of biaccessing points
(except for $\pm\beta_{c_p}$) and a countable family of interior arcs.

Each Fatou component on the spine defines a strip bounded by four dynamic
rays, which are moving holomorphically for all parameters in the wake of
$\M_p$ by the composition of B\"ottcher conjugations. This holomorphic motion
is extended according to the S{\l}odkowsky Theorem \cite{shm}. Neglecting
branch points, the biaccessing rays of the Cantor set are approximated by
strips from both sides, so their motion is determined uniquely. The moved
ray is still a dynamic ray of the new polynomial, and it is still
biaccessing. On the other hand, a pair of topologically biaccessing rays
for the parameter $c\in\M_p$ is iterated to a pair separating $-\beta_c$
from $\beta_c$\,. Neglecting pairs separating small Julia sets, this pair
belongs to the moving Cantor set. Thus $B_\top(c)=B_\top(c_p)$, and
$B_\comb(\theta)=B_\comb(\theta_\pm)$ is shown analogously, observing
that $\r_c(\theta/2)$ and $\r_c((\theta+1)/2)$ belong to the strip
around the small Julia set at $z=0$. \mybox

\subsection{Core entropy revisited} \label{4ce}
According to \cite{mesh, gaotl}, the following relation is due to Thurston.
The proof follows Bruin--Schleicher \cite{bshb} and Tiozzo \cite{tiob}.
See also the partial results by Douady \cite{dte}.

\begin{thm}[Thurston, Bruin--Schleicher, Tiozzo]\label{Tbhlc}
Suppose $\K_c$ is locally connected with empty interior, or $f_c$ is
parabolic or hyperbolic with a real multiplier.
Using regulated arcs, define the tree $T_c$ as the path-connected hull of
the critical orbit. If $T_c$ is compact, define the core entropy $h(c)$ as
the topological entropy of $f_c$ on $T_c$\,. Then it is related to
the biaccessibility dimension by $h(c)=B_\top(c)\cdot\log2$.
\end{thm}

\textbf{Proof:} First, assume that $\K_c$ has empty interior.
According to Douady \cite[p.~449]{dcc}, every biaccessible 
point $z\in\K_c$ has an image separating $-\beta_c$ from $\beta_c$\,:
when two of its external angles differ in the $n$-th binary digit,
$f_c^{n-1}(z)$ has external angles in $(0,\,1/2)$ and in $(1/2\,,\,1)$.
It will be iterated to the arc of $T_c$ through $z=0$
in finitely many steps. Conversely, every external angle of $T_c\subset\K_c$
is biaccessing, with at most countably many exceptions. Setting
$X_c=\gamma_c^{-1}(T_c)$ and considering a countable union of linear
preimages shows $B_\top(c)=\dim X_c$\,. Now $X_c$ is forward invariant under
the doubling map $F(\phi)=2\phi$, which has constant slope $2$.
So $\log2\cdot\dim X_c=h_\top(F,\,X_c)$ by \cite[Proposition~III.1]{fu}.
But $\gamma_c:X_c\to T_c$ is a semi-conjugation from $F(\phi)$ to
$f_c(z)$, so $h_\top(F,\,X_c)=h_\top(f_c\,,\,T_c)=h(c)$. Note that by
the Bowen Theorem \cite{bow1, mevst}, 
it is sufficient that every point $z\in T_c$ has finitely many preimages
under $\gamma_c$\,, but this number need not be bounded globally.
Here every branch point with more than four branches is periodic or
preperiodic by the No-wandering-triangles Theorem \cite{tgdr}, 
and periodic points in $\partial\K_c$ are repelling. So this
condition will be satisfied, even if there are infinitely many branch points.

If $f_c(z)$ is hyperbolic or parabolic, $\gamma_c(X_c)=T_c\cap\partial\K_c$
will be a disconnected subset of $T_c$\,. To see that
$h_\top(f_c\,,\,T_c\cap\partial\K_c)=h_\top(f_c\,,\,T_c)$ in the
maximal-primitive case of period $p$, the estimate $h(c)>\frac{\log2}p$ from
Proposition~\ref{Prenpcf}.2 shows again that the interior does not contribute
more to the entropy. In the pure satellite case, the number of preimages under
$f_c^n(z)$ is not growing exponentially on $T_c$\,, so $h(c)=0$. \mybox

Similar techniques show that the set of biaccessing pairs
$(\phi_1\,,\,\phi_2)\in S^1\times S^1$ has the same Hausdorff dimension
\cite{gao}. 
In the postcritically finite case, results on Hausdorff measure are found
in Proposition~\ref{Pbhpcf}. Thurston \cite{tabstr} has suggested to define
a core $T_c\subset\K_c$ for every $c\in\M$ as ``the minimal closed and
connected forward-invariant subset containing [the critical point], to which
the Julia set retracts.'' Consider the following cases:
\begin{itemize}
\item If $c$ is parabolic or hyperbolic of period $p$, the critical orbit is
connected with arcs from $z$ to $f_c^p(z)$, but $T_c$ will not be unique. If
$\K_c$ contains a locally connected Siegel disk $\mathcal{D}$ of period $1$,
the core is $\overline{\mathcal{D}}$.
For a $p$-cycle of locally connected Siegel disks, $T_c$ contains arcs
through preimages of these disks in addition.
\item In the fixed Cremer case and in the infinitely pure satellite
renormalizable case, the following situation may arise according to Douady
and S{\o}rensen \cite{soeacc, soeren}: the rays $\r_c(\theta/2)$ and
$\r_c((\theta+1)/2)$ each accumulate at a continuum containing both
$\alpha_c$ and $-\alpha_c$\,.
$\K_c$ is neither locally connected nor arcwise connected. The core
according to the above definition may be too large, even all of $\K_c$\,.
The same problem applies to Cremer cycles of pure satellite type, 
and possibly to the non-locally connected Siegel case.
\item In the primitive renormalizable case, we have
$h_\top(f_c\,,\,T_c)=B_\top(c)\cdot\log2$ for a fairly arbitrary choice of
$T_c$\,. Even if the small Julia set $\K_c^p$ is an indecomposable continuum,
we may include complete preimages of it in $T_c$\,, because $B_\top(c)>1/p$.
\item Suppose $\K_c$ is locally connected with empty interior, but the
arcwise connected hull of the critical orbit is not compact. Then it is
not clear, whether taking the closure of $T_c$ may produce an entropy
$h_\top(\overline{T_c})>B_\top(c)\cdot\log2$.
\end{itemize}

\begin{rmk}[Locally connected model and non-compact core] \label{Rlcm}
1. For any parameter angle $\theta$, Bruin--Schleicher \cite[version~1]{bshb}
consider the set $X_\theta$ of angles $\phi$ that are either postcritical or
combinatorially biaccessing and before a postcritical angle. In the cases
where $\K_c$ has empty interior and is not of Cremer type, this set
corresponds to a tree $T_\theta$ within a locally connected model of
$\K_c$\,, which is a quotient of $S^1$ or of a space of itineraries. Now
$h_\top(f,\,T_\theta)=h_\top(F,\,X_\theta)=B_\comb(\theta)\cdot\log2$ may be
shown using a generalized notion of topological entropy when $T_\theta$ is
not compact. It would be interesting to know whether the entropy on the
closure can be larger, and whether $T_\theta$ can be embedded into $\K_c$ in
the infinitely renormalizable case; probably this will not work in the
situation of \cite{soeren} described above.

2. Suppose the binary expansion of $\theta$ is a concatenation of all finite
words. Then the orbit of $\theta$ under doubling is dense in $S^1$. The
parameter $c$ is not renormalizable and not in the closed main cardioid,
so $\K_c$ is locally connected with empty interior.
Since the critical orbit is dense in $\K_c$\,, the tree $T_\theta$
corresponds to a non-compact union $T_c$ of a countable family of arcs in
$\K_c$\,. Now
$B_\top(c)=\dim\gamma_c^{-1}(T_c)<1=h_\top(\overline{T_c})/\log2$
shows that the generalized topological entropy satisfies
$\dim\gamma_c^{-1}(T_c)\neq h_\top(T_c)/\log2$ or
$h_\top(T_c)<h_\top(\overline{T_c})$.
\end{rmk}

\subsection{Monotonicity, level sets, and renormalization} \label{4ctmlr}
Monotonicity of the biaccessibility dimension (or core entropy) has been
shown by Penrose \cite{chprk} for abstract kneading sequences and by
Tao Li \cite{taoli} for postcritically finite polynomials,
cf.~Proposition~\ref{Pmonopcf}. In terms of angles it reads as follows:

\begin{prop}[Monotonicity of laminations, folklore]\label{Pmonocomb}
Suppose $\theta_-<\theta_+$ are such that for some parameter $c$, the
corresponding dynamic rays land together at the critical value $c$,
or $\theta_\pm$ are characteristic angles at $z_1$\,. If an angle $\phi$
is combinatorially biaccessing with respect to $\theta_\pm$\,, it stays
biaccessing for all angles $\theta\in(\theta_-\,,\,\theta_+)$.
In particular we have $B_\comb(\theta)\ge B_\comb(\theta_\pm)$.
\end{prop}

The assumption may be restated in terms of laminations of the disk \cite{tgdr}
as follows: $\theta_\pm$ belong to a leaf or a polygonal gap of the quadratic
minor lamination QML. If $c$ is a branch point, we may assume that
$\theta_\pm$ bound the whole wake. Their preimages divide the circle into four
arcs; the preimages of $\theta$ are contained in
$I_0=(\theta_-/2\,,\,\theta_+/2)$ and
$I_1=((\theta_-+1)/2\,,\,(\theta_++1)/2)$. Consider the set $B\subset S^1$ of
angles combinatorially biaccessing with respect to $\theta_\pm$\,, and the
subset $U\subset B$ of biaccessing angles $\phi$ that are never iterated to
$I_0\cup I_1$ under doubling $F(\phi)=2\phi$. For each $\phi\in U$ the
itinerary with respect to $\theta_-$ is the same as the itinerary with respect
to $\theta_+$\,. Now $B$ is the union of
iterated preimages of $U$ and $B_\comb(\theta_\pm):=\dim B=\dim U$. To see
this, we may assume that $c$ belongs to a vein according to
Proposition~\ref{Pbetatree}; otherwise $c$ would have a non-trivial fiber
intersecting a vein and it should be moved there. Then $\K_c$ is locally
connected, the beginning of the critical orbit is defining a finite tree
$T_c\subset\K_c$\,, and every biaccessible point is iterated to
$T_c\cap\partial\K_c$\,. Moreover, $T_c$ does not have external angles in
$I_0\cup I_1$\,, since no $z\in T_c\cap\partial\K_c$ is iterated behind
$c$ or $z_1$\,. So both $U$ and $B$ are iterated to
$\gamma_c^{-1}(T_c\cap\partial\K_c)$. Variations of the following argument
where used by Douady \cite{dte} in the real case, by Tao Li \cite{taoli}
in the postcritically finite case, and by Tan Lei for veins
\cite[Proposition~14.6]{tiob}.

\begin{figure}[h!t!b!]
\unitlength 0.001\textwidth 
\begin{picture}(990, 280)
\put(10, 0){\includegraphics[width=0.28\textwidth]{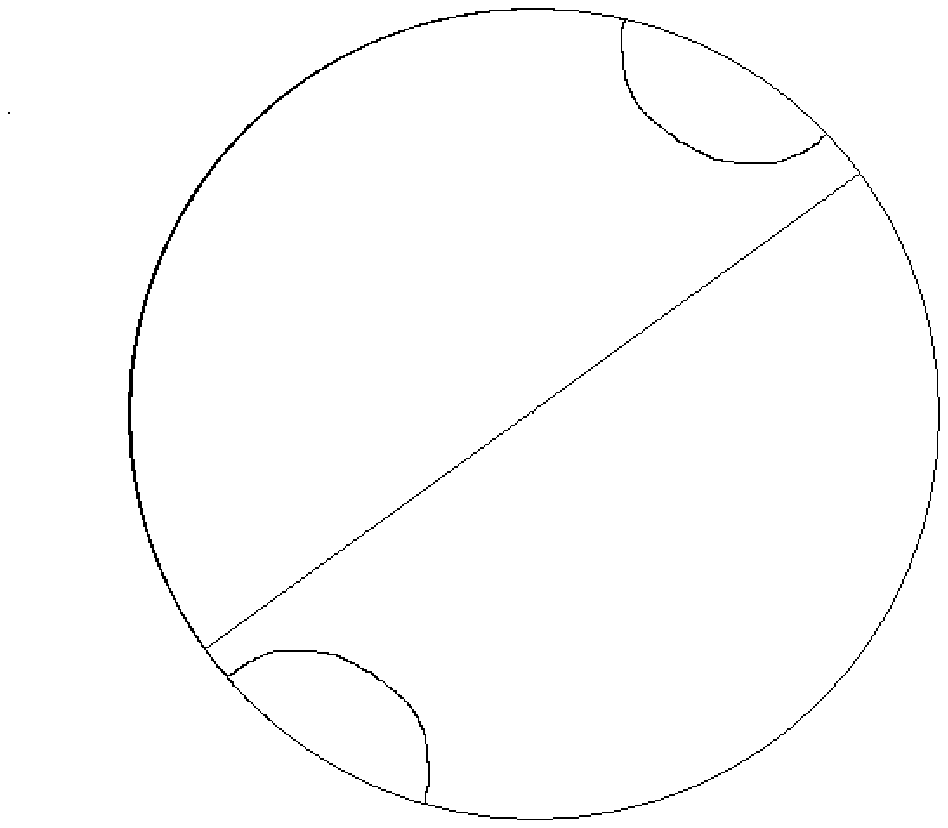}}
\put(360, 0){\includegraphics[width=0.28\textwidth]{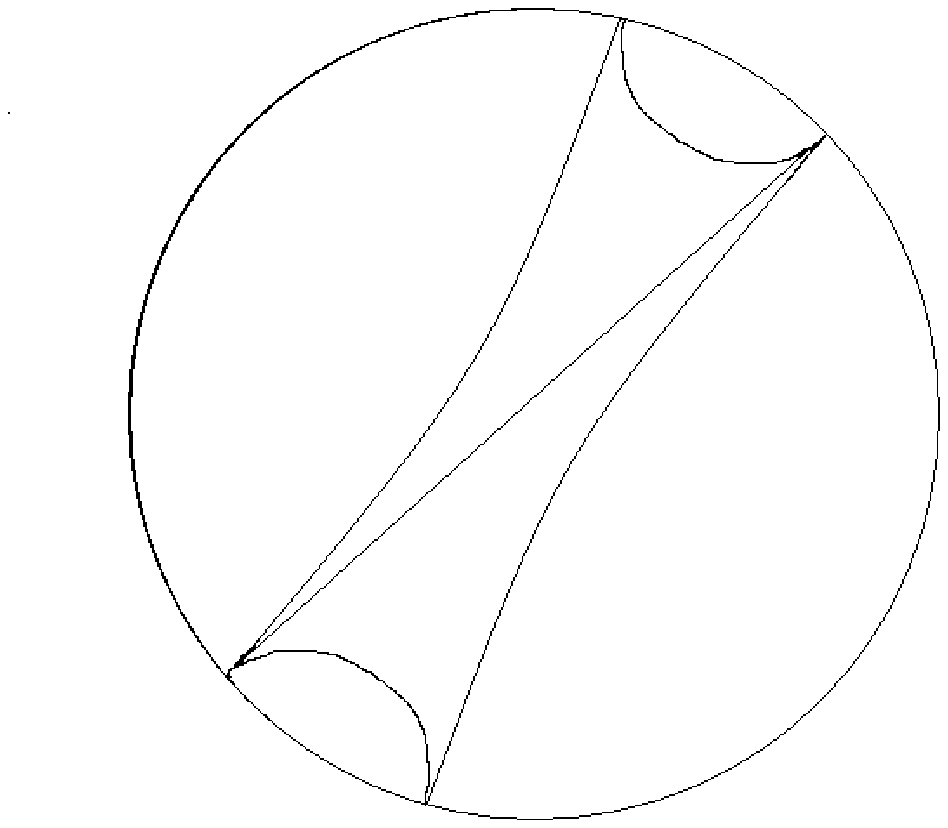}}
\put(710, 0){\includegraphics[width=0.28\textwidth]{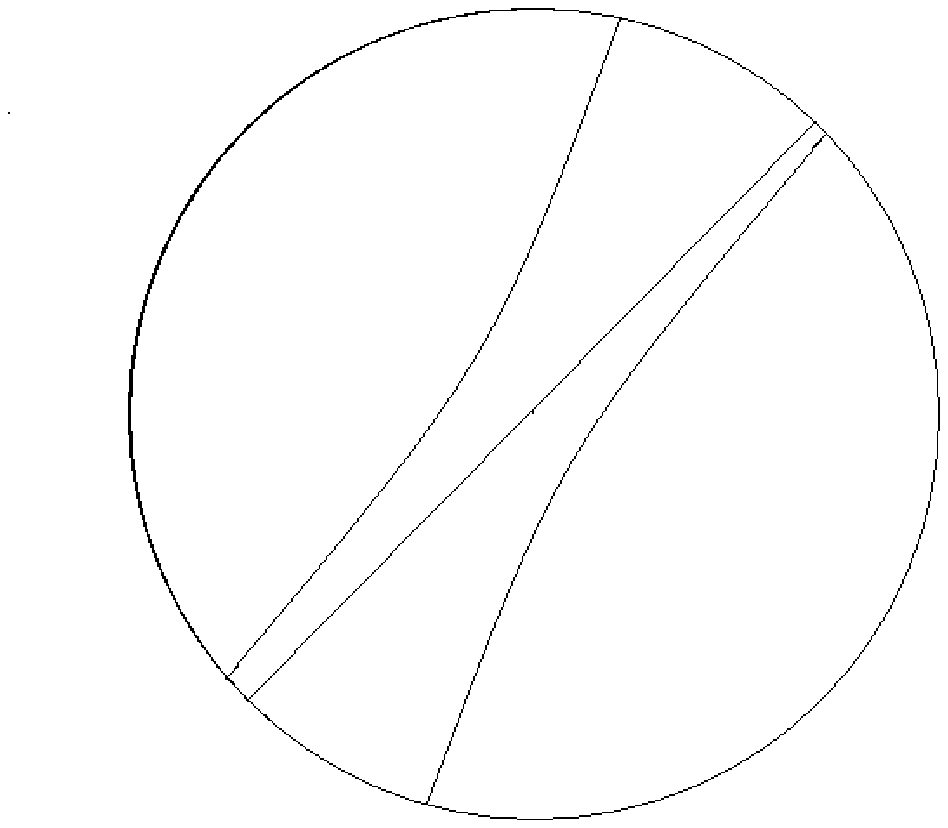}}
\put(11, 279){\makebox(0, 0)[lt]{a)}}
\put(361, 279){\makebox(0, 0)[lt]{b)}}
\put(711, 279){\makebox(0, 0)[lt]{c)}}
\end{picture} \caption[]{\label{Fcircles}
The pattern of rays landing together can be described by a lamination of the
disk \protect\cite{tgdr}. As $\theta$ is varied, the diameter defined by
$\theta/2$ and $(\theta+1)/2$ is moving and disconnecting or reconnecting
chords. The closed lamination is describing combinatorially biaccessing pairs
of angles; the corresponding rays need not land together if $\K_c$ is
non-locally connected.}
\end{figure}

\textbf{Proof of Proposition~\ref{Pmonocomb}:} The itinerary of any pair
$\phi_1\,,\phi_2\in U$ remains the same for all angles
$\theta\in(\theta_-\,,\,\theta_+)$, so $\phi_1$ and $\phi_2$ remain
biaccessing. Their preimages in $B$ remain biaccessing
as well. As an alternative argument, consider the possible change of
partners for $\phi\in B$ according to Figure~\ref{Fcircles}: since $\phi$
is iterated through $I_0\cup I_1$ and $(\theta_-\,,\,\theta_+)$ a finite
number of times, it has a finite number of possible itineraries with respect
to $\theta\in(\theta_-\,,\,\theta_+)$. Now $\phi$ is changing the partner when
$\theta$ is an iterate of $\phi$ or conjugate to an iterate of $\phi$.
The partners remain faithful on the finite number of remaining intervals
for $\theta$, since their itineraries do not change there. \mybox

So, what kind of bifurcation is changing the biaccessibility dimension?
The bifurcation of rational rays at a root is known explicitly, but it
does not affect the dimension. The periodic pairs stay the same for
parameters in the wake, while preperiodic and irrational pairs may be
regrouping. Now the dimension is changed only when an uncountable family
of irrational angles is gaining or losing partners. The bifurcations should
be studied for real parameters first: with $\|\phi\|:=\min\{\phi,\,1-\phi\}$,
and real $c=\gamma_\sM(\theta)$, the external angles $\phi$ of the Hubbard
tree $[c,\,f_c(c)]$ are characterized by $\|F^n(\phi)\|\le\|\theta\|$ for
all $n\ge0$. All biaccessible angles satisfy $\|F^n(\phi)\|\le\|\theta\|$ for
$n\ge N_\phi$\,. E.g., $\phi=0.010010001\dots$ is biaccessible for
$c=-2=\gamma_\sM(1/2)$ but not for any other real parameter.

Suppose $c_n\nearrow c_0$ with $B_\top(c_n)$ strictly increasing. Denote
the sets of biaccessing angles by $B_n$ and $B_0$\,. Consider the union
$B_-=\bigcup B_n$ and the newly biaccessible angles $B_+=B_0\setminus B_-$\,.
Then it will be hard to prove $\dim B_n\to\dim B_0$ by estimating
$\dim B_+$\,, since we have $\dim B_-=\lim\dim B_n\le\dim B_0=\dim B_+$ in
general. The last equality is obtained by considering the Hausdorff measure
for $b=\dim B_0$\,: now $\dim B_n<b$ gives $\mu_b(B_-)=0$, but
$\mu_b(B_0)>0$ according to Proposition~\ref{Pbhpcf}.

For pure satellite renormalization, the relation
$B_\top(c)=\frac1p\,B_\top(\hat c)$ from Proposition~\ref{Prenpcf}.1 extends
to all parameters $\hat c\in\M$ by the Douady substitution of binary digits
\cite{daa, wje}, since the covering intervals of length $2^{-np}$
correspond to intervals of length $2^{-n}$. The scaling implies that
$B_\top(c)\to0=B_\top(c_0)$ for a pure satellite Feigenbaum point $c_0$ and
$c\searrow c_0$ on the vein behind it. Consider primitive
renormalization:

\begin{thm}[Strict monotonicity and primitive Mandelbrot sets]
\label{Tmonoren}
$1$. The topological biaccessibility dimension is $B_\top(-2)=1$ and
$B_\top(c)=0$ for parameters $c$ in the closed main cardioid, in closed pure
satellite components, and for infinitely pure satellite renormalizable
parameters. We have $0<B_\top(c)<1$ for all other $c\in\M$.

$2$. For $c\prec c'$ we have $B_\top(c)\le B_\top(c')$. When
$B_\top(c')>0$, the inequality is strict unless $c$ and $c'$ belong to the
same primitive Mandelbrot set.

$3.$ The biaccessibility dimension $B_\top(c)$ is constant on every primitive
Mandelbrot set $\M_p=c_p\ast\M$. Every connected component of a level set
$\{c\in\M\,|\,B_\top(c)=b>0\}$ is either a maximal-primitive Mandelbrot set
or a single point.

$4$. For the primitive $\M_4\subset\M_{1/3}$ there is a sequence of parameters
$c_n\in\M\setminus\M_4$ converging to the root of $\M_4$\,, such that
$B_\top(c_n)=B_\top(\M_4)$.\\
There is a sequence of parameters $c_n\in\M\setminus\{\i\}$ converging to
the non-renormalizable Misiurewicz point $c=\i$, such that
$B_\top(c_n)=B_\top(\i)$.
\end{thm}

Analogous statements hold for $B_\comb(\theta)$. When $B_\top(c)$ is
restricted to a vein, items~1--3 imply that it is constant before a Feigenbaum
point, and on any arc corresponding to a primitive small Mandelbrot set. For
principal veins this was shown independently by Tiozzo \cite{tiob}.
See also the historical remarks after Theorem~\ref{Tbev}.

\textbf{Proof of Theorem~\ref{Tmonoren}:}\\
1. $B_\top(-2)=1$ is obtained from $\K_{-2}=[-2,\,2]$ or from the
$1\times1$-Markov matrix $A=(2)$. The pure satellite case was discussed in
the proof of Theorem~\ref{Tcombtop}.
In the remaining cases, $B_\top(c)>0$ follows from
Proposition~\ref{Prenpcf}.2 and monotonicity, or from the estimate $L_2$ in
\cite{bshb}. See \cite{zdun, ssms, bks, mesh, bshb} for $B_\top(c)<1$.

2. There are angles $\theta_-\le\theta\le\theta_+$\,, such that $\theta_\pm$
are external angles of $c$ and $\r_\sM(\theta)$ accumulates or lands at
$c'$, or at the corresponding root, or at a small Mandelbrot set around the
fiber of $c'$. The monotonicity statement of Proposition~\ref{Pmonocomb} is
transferred to the parameter plane. Under the additional assumptions, the
parameters are separated by a maximal-primitive Mandelbrot set \cite{sf2}, and
Proposition~\ref{Prenpcf}.3 gives strict monotonicity.

3. On a primitive Mandelbrot set, the biaccessibility dimension is constant
by the proof of Theorem~\ref{Tcombtop}. The components of $\M\setminus\M_p$
are attached to $\M_p$ at the root and at Misiurewicz points. If $\M_p$ is
maximal in the family of primitive Mandelbrot sets, these points are
approached by smaller maximal-primitive Mandelbrot sets, so $B_\top(c)$ is
strictly larger behind the Misiurewicz points. For the same reason, it is
strictly smaller on the vein before the root, but not necessarily on branches
at the vein, see item~4. If $B_\top(c)>0$ and $c$ is not primitive
renormalizable, then it is at most finitely renormalizable with all cycles
repelling, so its fiber is trivial by the Yoccoz Theorem. It can be separated
from any other parameter by a maximal-primitive Mandelbrot set again.

4. Denote the root by $c_\ast$ in the first case, and $c_\ast=\i$ in the
second case. Define $a_n$ as the sequence of $\alpha$-type Misiurewicz
points of lowest orders before $c_\ast$\,, then $B_\top(a_n)<B_\top(c_\ast)$
by item~2. In the other branch of $a_n$\,, the $\beta$-type Misiurewicz point
of lowest order satisfies $B_\top(b_{n\pm1})>B_\top(c_\ast)$ by
Examples~\ref{X315} and~\ref{X16}, respectively. Continuity on the vein to
$b_{n\pm1}$ according to Theorem~\ref{Tbev} provides a parameter $c_n$ on the
arc from $a_n$ to $b_{n\pm1}$ with $B_\top(c_n)=B_\top(c_\ast)$. \mybox

\begin{figure}[h!t!b!]
\unitlength 0.001\textwidth 
\begin{picture}(990, 420)
\put(10, 0){\includegraphics[width=0.42\textwidth]{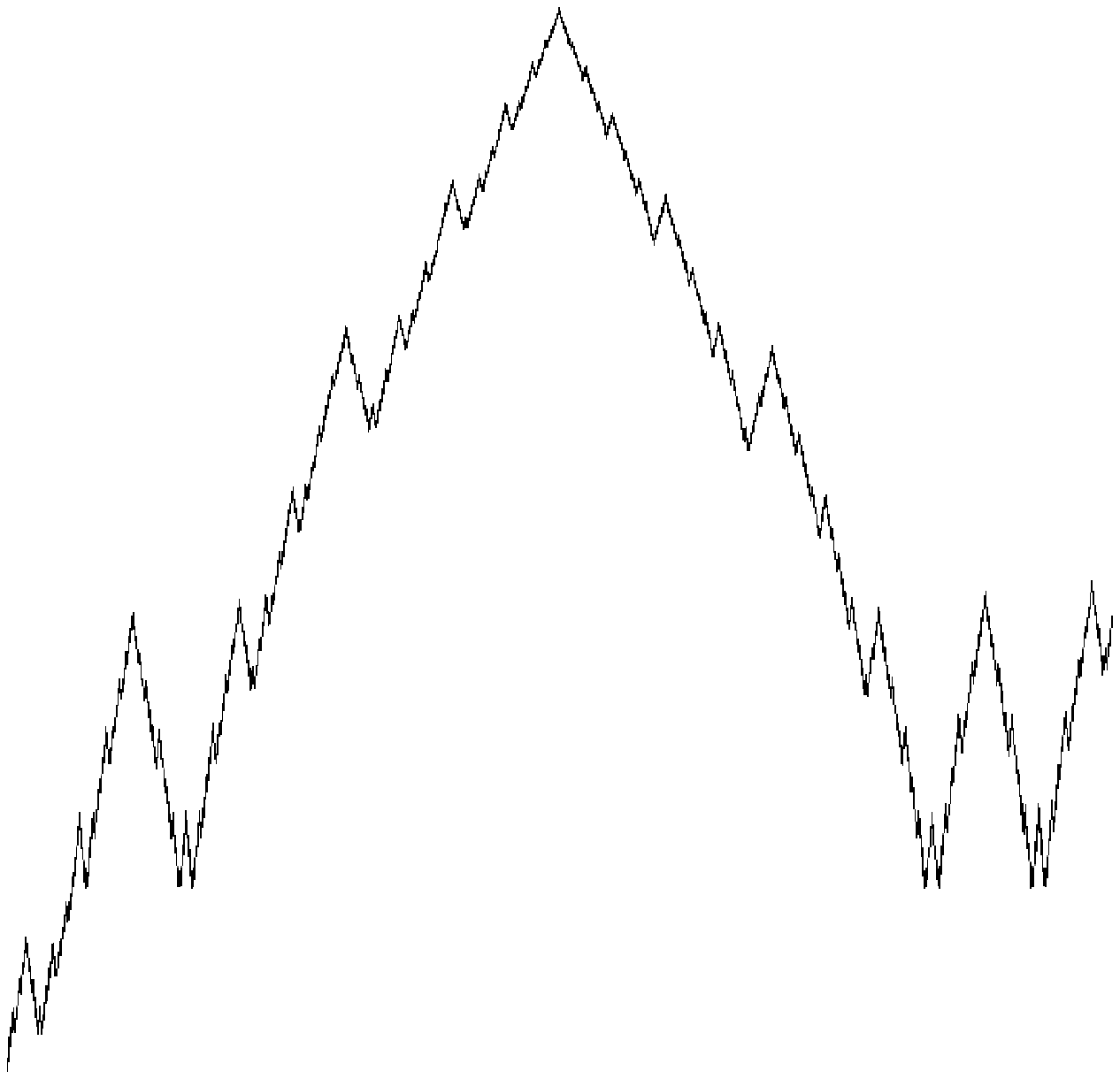}}
\put(570, 0){\includegraphics[width=0.42\textwidth]{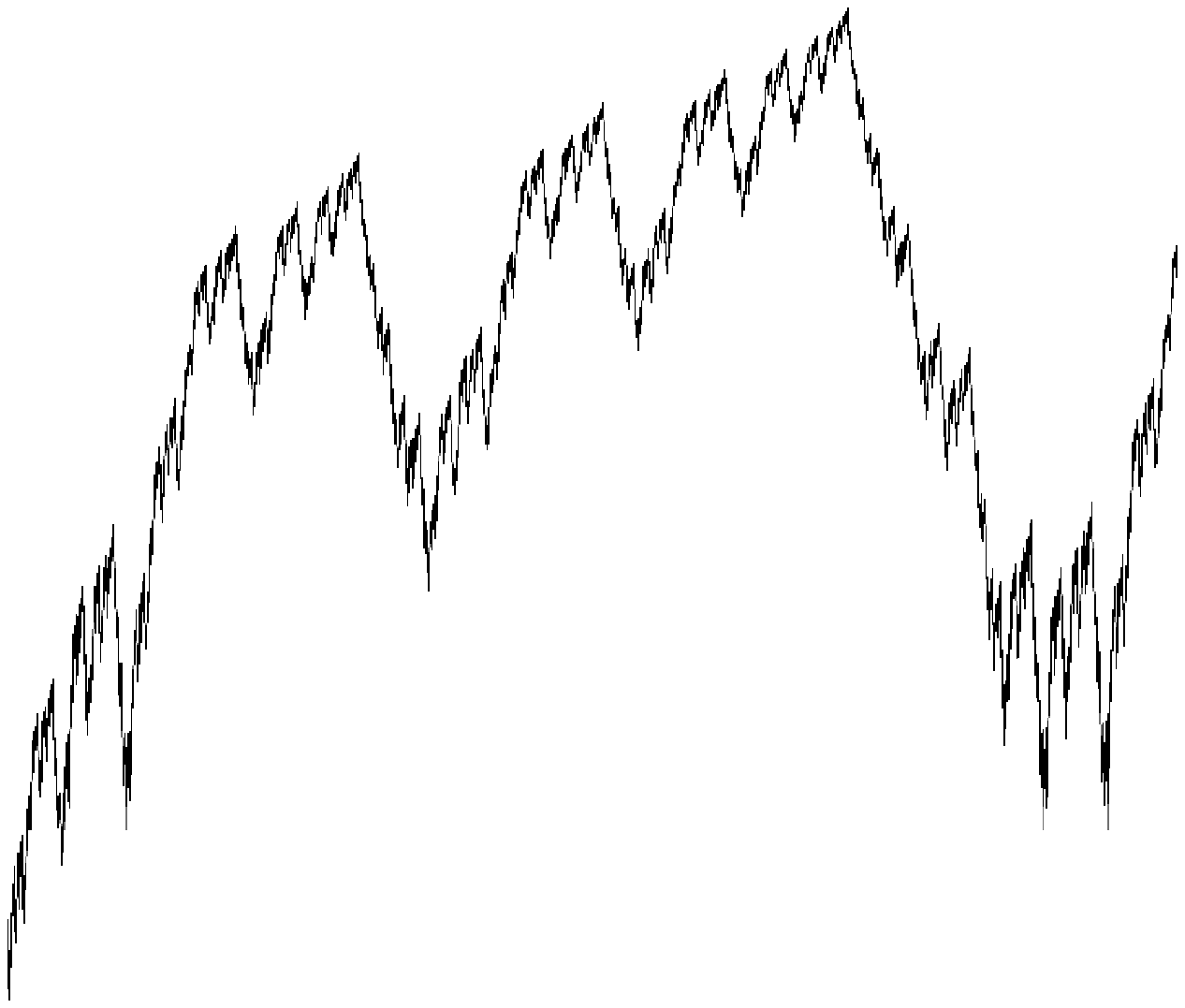}}
\thinlines
\multiput(10, 0)(560, 0){2}{\line(1, 0){420}}
\multiput(10, 420)(560, 0){2}{\line(1, 0){420}}
\multiput(10, 0)(560, 0){2}{\line(0, 1){420}}
\multiput(430, 0)(560, 0){2}{\line(0, 1){420}}
\end{picture} \caption[]{\label{Fm34}
The graph of $\lambda(\theta)$ representing $B_\comb(\theta)$ on intervals
containing $[3/7\,,\,25/56]$ and $[3/15\,,\,49/240]$, which are corresponding
to one decorated side of the primitive Mandelbrot sets $\M_3$ and $\M_4$\,.
Note that in both cases, $p$-renormalizable angles are on the same level and
most of the interval represents the $1/3$-limb.
There seems to be a left-sided maximum at $3/7$ but not at $3/15$.
In both cases there is a similarity between the various sublimbs; see
Figure~\protect\ref{Fm4l} for $\M_4$\,. In the case of $\M_3$\,, the symmetry
of the $1/2$-limb seems to be transferred to other sublimbs.}
\end{figure}

\subsection{Entropy and biaccessibility on veins} \label{4ebv}
Summarizing some classical results:

\begin{prop}[The tree of veins, folklore]\label{Pbetatree}
$1$. Every $\beta$-type Misiurewicz point is connected to 0 by a unique
regulated arc within $\M$. These arcs form an infinite tree; its
branch points are Misiurewicz points and centers. The tree is the disjoint
union of veins, such that veins of higher preperiods branch off from
a vein of lower preperiod.

$2$. Every biaccessible parameter belongs to a vein, and the
corresponding dynamic rays land at the critical value or at the
characteristic point. Hyperbolicity is dense on every vein.

$3$. For every parameter $c$ on a vein, $\K_c$ is locally connected and the
regulated path-connected hull $T_c$ of the critical orbit is a finite tree.
We have $T_c\subset\tilde T_c$\,, where the tree $\tilde T_c$ is
homeomorphic to the Hubbard tree of the $\beta$-type endpoint of the vein.
\end{prop}

Note that it is not known, whether $\M$ is locally connected on a vein;
this is open for Feigenbaum points in particular. The vein has a linear order
$\prec$ extending the partial order of $\M$ to interior parameters and to
possible non-accessible boundary points.

\textbf{Proof}: 1. The veins were described by Douady--Hubbard in
terms of the partial order $\prec$ in a locally connected Model of $\M$
\cite{dcc}. To see that they are topological arcs, note that any non-trivial
fiber would be infinitely renormalizable according to Yoccoz \cite{hy3};
but a vein meets a small Mandelbrot set in a quasi-conformal image of
a real interval, plus finitely many internal rays in satellite components.
This argument is due to Kahn \cite{dcc, sf3}, and homeomorphisms from
intervals onto complete veins are due to Branner--Douady and Riedl
\cite{bd, rt}. By the Branch Theorem \cite{sf2, tgdr}, the tree branches
at Misiurewicz points with finitely many branches and at centers with
countably many internal rays leading to sublimbs. There is a unique branch
of lowest preperiod, since any open interval contains a unique dyadic angle
of smallest denominator. By this order, the decomposition into veins is
well-defined: the branch of lowest order is continued before the branch point.

2. If $c$ is a biaccessible parameter, it is on the vein to the $\beta$-type
Misiurewicz point of lowest preperiod in its wake, which corresponds to the
dyadic angle of lowest preperiod in an interval. Hyperbolicity is dense on
the real axis according to Graczyk--Swi\c{a}tek and Lyubich \cite{grsw, l12},
so possible non-trivial fibers can meet it in a single point only, as
illustrated in Figure~\ref{FFB}. Triviality of fibers is preserved
under tuning \cite{sf3}, so all non-hyperbolic parameters on a vein are roots
or separated by roots, and approximated by hyperbolic arcs. The approximation
by roots shows in addition, that $c$ or $z_1$ has the same characteristic
angles as the parameter, employing the approximation by corresponding
characteristic points and local connectivity of $\K_c$ from item~3.
--- Suppose that $c\in\M$ and the critical value $c\in\K_c$ is biaccessible.
It is not known whether this implies that $c$ belongs to a vein of $\M$;
otherwise $c$ belongs to a non-trivial fiber crossing the candidate vein.

3. The Julia sets of real polynomials are locally connected according to
Levin--van~Strien \cite{lscr}. If $\K_c$ was not locally connected, it would
be infinitely simply renormalizable and obtained by tuning from a real
parameter; but local connectivity would be carried over according to
Schleicher \cite{sf3}. Consider the tree $\tilde T_c$ of regulated arcs in
$\K_c$ connecting the orbit of a preimage of $\beta_c$ corresponding to the
endpoint of the vein. The beginning of the critical orbit gives the endpoints
of $T_c$ by the same arguments as for postcritically finite Hubbard trees.
$\tilde T_c$  is forward invariant and contains the critical value $z=c$
according to item~2, thus the entire critical orbit, so
$T_c\subset\tilde T_c$.  The tree $\tilde T_c$ bifurcates only at parameters
on veins of lower order (and not at the origin of a non-principal vein).
The topological entropy on $T_c$ and $\tilde T_c$ is the same; for
postcritically finite dynamics this is Lemma~\ref{LHtree}.4, and the general
case is obtained from a recursion for the number of preimages of $z=0$ on
$\tilde T_c\setminus T_c$\,. \mybox

\begin{thm}[Biaccessibility and entropy on veins of $\M$,
generalizing Tiozzo]\label{Tbev}
On any vein $v$ of $\M$, we have:

$1$. The Hubbard tree $T_c$ is a finite tree for $c\in v$, thus compact,
and the core entropy satisfies $h_\top(f_c\,,\,T_c)=h(c)=B_\top(c)\cdot\log2$.

$2$. It is $0$ on an arc through pure satellite components, constant on
the arcs corresponding to maximal-primitive Mandelbrot sets, and strictly
increasing otherwise.

$3$. The core entropy $h(c)$ and biaccessibility dimension $B_\top(c)$
are continuous on $\overline{v}$.

$4$. The Hausdorff dimension of biaccessing parameter angles of $v$
before $c$ equals $B_\top(c)$.
\end{thm}

Item~1 is obtained from Proposition~\ref{Pbetatree}.3 and Theorem~\ref{Tbhlc},
and item~2 follows from Theorem~\ref{Tmonoren}. Continuity is proved below.
See Proposition~\ref{Parcpd} for item~4. Tiozzo \cite{tiob} has obtained
these results for principal veins with different methods, e.g., combinatorics
of dominant parameters, piecewise-linear maps on trees, and a kneading
determinant. Probably continuity can be proved with the approach of
Bruin--Schleicher \cite{bshb} as well.

The real interval $[-2,\,0]$ is the principal vein of the $1/2$-limb.
Here Milnor--Thurston \cite{mtimi} constructed a semi-conjugacy from
$f_c(z)$  to a tent map and showed that the topological entropy depends
continuously and monotonically on the kneading invariant. Later the Douady map
$\Phi_\sM(c):=\Phi_c(c)$ implied monotonicity with respect to the parameter
$c$, by relating the symbolic description to the external angle of $c$
\cite{dte}. An alternative proof is due to Tsujii \cite{tsu, mtm}.

The following proofs are based on continuity results for the entropy of tree
endomorphisms. Now the Hubbard tree $T_c$ or $\tilde T_c$ is homeomorphic
to a linear tree and $f_c(z)$ is conjugate to a tree map, but it would be
difficult to control parameter dependence of the homeomorphisms. Instead only
one polynomial is transferred to a linear tree and perturbed continuously
there; combinatorics relates postcritically finite maps back to polynomials.
The second proof is motivated by the idea of \cite{dte} for real polynomials,
but it needs deeper results on tree endomorphisms than the first one.

\textbf{First proof of continuity on veins:} The core entropy is constant on
closed arcs corresponding to small Mandelbrot sets of maximal-primitive
type. If the origin of the vein $v$ is a pure satellite center, then $h(c)$
is 0 on a closed arc meeting pure satellite components, and it is continuous
at the corresponding Feigenbaum point by tuning. 
Since $h(c)$ is monotonic, we need to show that it does not jump.
Suppose $c_0$ is behind the Feigenbaum point and not primitive renormalizable.
Then $c_0$ has trivial fiber \cite{sf2} and is not a satellite root,
so it is approximated monotonically from below by Misiurewicz points $c_n$\,.
In the Hubbard tree $T_{c_0}\subset\K_{c_0}$ denote the characteristic point
corresponding to $c_n$ by $z_n$\,, then $z_n\to z_0$ from below. Now the tree
$T_{c_0}$ is homeomorphic to a tree consisting of line segments, with a map
$g_0(x)$ conjugate to $f_{c_0}(z)$ and points $x_n$ corresponding to $z_n$\,.
Choose $x_n'$ with $x_{n-1}\preceq x_n'\prec x_n$\,, such that the orbit of
$x_n$ does not return to $(x_n'\,,\,x_0]$, and define a sequence of continous
maps $\eta_n(x)$ with $\eta_n:[x_n'\,,\,x_0]\to[x_n'\,,\,x_n]$ monotonically
and $\eta_n(x)=x$ for $x$ before $x_n'$\,, such that
$g_n=g_0\circ\eta_n\to g_0$ uniformly. Then
$\liminf h_\top(g_n)\ge h_\top(g_0)=h(c_0)$ according to \cite{lmg}. But
$g_n(x)$ is postcritically finite with the same order of the critical orbit
as $f_{c_n}(z)$, which implies $h_\top(g_n)=h(c_n)$ since lap numbers of
iterates are determined by itineraries on the tree. 
By monotonicity we have $\limsup h(c_n)\le h(c_0)$. So there is no jump of
$h(c)$ for $c\to c_0$ from below. The construction of $g_n(x)$ needs to be
modified when $c_0$ is a maximal-primitive root,  since $z_0$ is strictly
before the critical value $c_0$\,: all arcs through Fatou basins are squeezed
to points first. (By Proposition~\ref{Prenpcf}.2, this does not change the
core entropy.)

Now suppose $c_0$ is behind the Feigenbaum point and either a
maximal-primitive tip or not primitive renormalizable, so
$h(c_0)=\log\lambda_0>0$. According to \cite{alm}, $\lambda_0$ is the
growth rate of lap numbers,
$\lim\frac1n\log L(f_{c_0}^n)=h(c_0)$. For $\lambda>\lambda_0$ there is
an $N$ with $L(f_{c_0}^n)\le\lambda^n$ when $n\ge N$. Lap numbers
are changing at postcritically finite parameters in $v$ only, since the
existence and the linear order of endpoints of laps is determined from
their itineraries. So there is a one-sided neighborhood
$U=[c_0\,,\,c')\subset v$ with $c'\succ c_0$\,, such that the lap number
$L(f_c^n)$ is constant for $c\in U$ when $n\le N$. Since each $f_c(z)$ on
$T_c$ is conjugate to a tree endomorphism, its entropy is related to the
growth of lap numbers. Now these are sub-multiplicative, so
$L(f_c^{kN})\le [L(f_c^N)]^k\le\lambda^{kN}$ and letting $k\to\infty$ gives
$\log\lambda_0=h(c_0)\le h(c)\le\log\lambda$. This shows continuity of $h(c)$
for $c\to c_0$ from above by taking $\lambda\searrow\lambda_0$\,. \mybox

On the real vein, continuity was proved by different methods in
\cite{mtimi, miszepm, dte}. The typical counterexample to continuity involves
a change of the lap number $L(g^1)$ or a modification of a 0-entropy map
with a periodic kink. According to \cite{miszepm, alm}, a continuous
piecewise-monotonic interval map $g_0(x)$ has the following properties
under a $C^0$-perturbation $g(x)$ respecting the lap number: the topological
entropy is lower semi-continuous for $g\to g_0$ and it can jump at most to
$\frac q p\log2$, if $g_0(x)$ has a $p$-cycle containing $q$ critical points.
Continuity is obtained when $h_\top(g_0)\ge\frac q p\log2$ already, or from
$C^2$-convergence.

These results were generalized to non-flat piecewise monotonic-continuous
maps in \cite{mishpc}. The highest possible jump is bounded by
$\limsup h_\top(g)\le\max(h_\top(g_0)\,,\,\log\lambda)$, where $\lambda$ is
the highest eigenvalue of the transition matrix describing the orbits of
critical points and points of discontinuity, as long as these are mapped to
each other in a finite number of iterations. By concatenating the edges, any
continuous tree map can be described by a piecewise continuous interval map;
if the topological entropy is discontinuous, the bound $\lambda$ is related to
the orbits of critical points and branch points. For quadratic polynomials,
the critical point is never a periodic branch point, and a cycle of
branch points gives $\lambda=1$. By Proposition~\ref{Prenpcf}.2, a center has
$h(c_p)>\frac1p\log2$ unless it is of pure satellite type. So:

\textit{Suppose $g_0(x)$ is an endomorphism of a finite linear tree $T$ with
$h_\top(g_0)>0$, which is topologically conjugate to a quadratic polynomial
$f_c(z)$ on $T_c$\,. Then $h_\top(g)\to h_\top(g_0)$ when $g\to g_0$
uniformly, where the maps $g^1$ have the same lap number as $g_0^1$ on $T$.}

\textbf{Second proof of continuity on veins:} Assume that the monotonic core
entropy $h(c)$ has a jump discontinuity on the vein $v$, omitting an interval
$(h_0\,,\,h_1)$. Again by tuning, an initial pure satellite Feigenbaum point
has $h(c_{\mbox{\tiny F}})<h_0$. Choose
Misiurewicz points $c_0\prec c_1$ on $v$ with $0<h(c_0)\le h_0<h_1\le h(c_1)$.
Now $f_{c_1}(z)$ on $T_{c_1}$ is topologically conjugate to an endomorphism
$g_1(x)$ on a tree consisting of line segments. It contains the orbit of
a preperiodic characteristic point corresponding to $c_0$\,, so there
is a continuous family $g_t(x)$, $0\le t\le1$, on the linear tree with
$g_0(x)$ representing the combinatorics of $c_0$\,. $g_t(x)$ shall have the 
same critical point 
as $g_1(x)$. By the adaptation of \cite{mishpc} above, the topological entropy
of $g_t(x)$ could be discontinuous only when the critical point is
$p$-periodic and $h_\top(g_t)<\frac1p\log2$; there would be a center $c\in v$
with the same order of the critical orbit \cite{bks} and $h(c)>\frac1p\log2$.
So there are $0\le t_0<t_1\le 1$ with $h_\top(g_{t_i})=h_i$\,.
Since the topological entropy on the tree is determined by the growth rate
of lap numbers according to \cite{alm}, and lap numbers are changing only
at postcritically finite parameters, there is a $t$ with $t_0<t<t_1$ and
$h_0<h_\top(g_t)<h_1$\,, such that $g_t(x)$ is postcritically finite. Again,
the critical orbit is defining an oriented Hubbard tree and there is a unique
quadratic polynomial $f_c(z)$ with the same dynamics. The parameter
$c$ with $c_0\prec c\prec c_1$ 
and $h_0<h(c)<h_1$ contradicts the assumed jump discontinuity. \mybox

The proofs did not use the semi-conjugation from $f_c(z)$ to
a piecewise-linear map of slope $\lambda$ \cite{bcplm}. An alternative
proof might be given recursively by showing that the set of admissible
$\lambda$ is connected for a particular shape of the tree. For principal
veins, the piecewise-linear model was described explicitly in terms
of external angles \cite{tiob}.

\subsection{Continuity with respect to the angle and the parameter}\label{4cc}
Monotonicity and continuity of the core entropy on veins implies:

\begin{prop}[Partial results on continuity] \label{Ppcont}
$1$. The biaccessibility dimension $B_\top(c)$ is continuous on the union
of regulated arcs connecting a finite number of postcritically finite
parameters within $\M$.

$2$. Suppose that $\theta_0\in S^1$ and $c_0\in\partial\M$ belongs to the
impression of $\r_\sM(\theta_0)$. Assume that $c_0$ is biaccessible,
postcritically finite, simply renormalizable, or in the closed main cardioid.
Then the biaccessibility dimension is lower semi-continuous as
$\theta\to\theta_0$ or as $c\to c_0$\,: we have
$\liminf B_\comb(\theta)\ge B_\comb(\theta_0)$ and
$\liminf B_\top(c)\ge B_\top(c_0)$.
\end{prop}

\textbf{Proof:} 1. Continuity on veins was shown in Section~\ref{4ebv}. If
a postcritically finite parameter is not on a vein, it is a Misiurewicz point
of period $>1$ with one external angle, and $\M$ is openly locally connected
at this point. So an arc consisting of a countable family of subsets of veins
ends at this parameter. Items~2 and~3 of Proposition~\ref{Pbetatree}
remain true for this arc, and the proofs of continuity are transferred.

2. If $c_0$ is renormalizable and not a primitive root, there is a
biaccessible $c_1\prec c_0$ with $B_\top(c_0)=B_\top(c_1)$. Monotonicity
gives $B_\top(c)\ge B_\top(c_0)$ for $c\succeq c_1$\,, which is a
neighborhood of $c_0$\,. If $c_0$ is a primitive root, a postcritically finite
endpoint, or a non-renormalizable biaccessible parameter with $B_\top(c_0)>0$,
there is a sequence of parameters $a_n\nearrow c_0$ with
$B_\top(a_n)\nearrow B_\top(c_0)$. Now $c\to c_0$ implies that $c\succ a_n$
with $n\to\infty$. \mybox

According to \cite{bks} the biaccessibility of an angle $\phi_1$ can be
seen by comparing its itinerary to the kneading sequence, without knowing
the second angle $\phi_2$\,: every branch of $\K_c$ at the landing point is
characterized by an eventually unique sequence of closest precritical points.
This approach was used in version~1 of \cite{bshb} to obtain H\"older
continuity with respect to the kneading sequence $\nu$ by estimating the
number of words appearing in the itineraries of biaccessing angles.
However, the proof is currently under revision:

\begin{conj}[Continuity on $S^1$, Thurston and Bruin--Schleicher]%
\label{Ccombcont}
The combinatorial biaccessibility dimension $B_\comb(\theta)$ is
continuous on $S^1$.
Moreover, when $\theta_0\in\Q$ and $B_\comb(\theta_0)>0$, it is H\"older
continuous at $\theta_0$ with exponent $B_\comb(\theta_0)$\,.
\end{conj}

To transfer H\"older continuity from $\nu$ to $\theta$, consider the
intervals of angles $\theta$ where the first $N$ entries of $\nu(\theta)$
are constant: when $\theta_0$ is rational, they do not shrink to $\theta_0$
faster than $\sim2^{-N}$. The H\"older exponent may become smaller for
irrational $\theta_0$\,. And $B_\comb(\theta)$ is not H\"older continuous at,
e.g., the root at $\theta_0=0$: denote by $\theta_n=2^{-n}$ the dyadic angle
of lowest preperiod in the limb of rotation number $1/(n+1)$. According to
Example~\ref{Xprinc}, we have $B_\comb(\theta_n)=\log\lambda_n/\log2$ with
$\lambda_n\to1$ and $(\lambda_n)^n=2/(\lambda_n-1)\to\infty$.

Consider the Feigenbaum point $c_{\mbox{\tiny F}}= -1.401155189$ of pure
satellite period doubling. It is approximated by real Misiurewicz points
or centers $c_k$ of preperiod and period $\asymp2^k$ and
$B_\comb(\theta_k)=B_\top(c_k)\asymp2^{-k}$ \cite{sts}, since the tuning
map $\hat c\mapsto (-1)\ast\hat c$ is halving the dimension. Repeated Douady
substitution gives $\log(\theta_k-\theta_{\mbox{\tiny F}})^{-1}\asymp2^k$,
so $B_\comb(\theta)$ is not H\"older continuous at $\theta_{\mbox{\tiny F}}$
\cite{tiob}. Since $c_{\mbox{\tiny F}}-c_k\asymp\delta_{\mbox{\tiny F}}^{-k}$
with $\delta_{\mbox{\tiny F}}=4.669201609$, $B_\top(c)$ is H\"older
continuous on the real axis at $c_{\mbox{\tiny F}}$ with exponent
$\log2/\log\delta_{\mbox{\tiny F}}=0.449806966$. Note that the
biaccessibility dimension on a sequence of angles or parameters gives
upper and lower bounds on a particular vein due to monotonicity, but only the
lower bound applies to branches of the vein as well. The bounds
at the 0-entropy locus in \cite{bshb} are not restricted to veins.

\begin{thm}[Continuity on $\M$]\label{TCtopcont}
\textsc{Assuming Conjecture~\ref{Ccombcont},
that $B_\comb(\theta)$ is continuous on $S^1$\,:}\\
The topological biaccessibility dimension $B_\top(c)$ is continuous on $\M$.
\end{thm}

This implication was obtained independently by Bruin--Schleicher.
Two different proofs shall be discussed here: the first one considers similar
cases as in the proof of Theorem~\ref{Tcombtop}. The second one will be
shorter but more abstract. See version~2 of \cite{bshb} for an alternative
argument based on compactness of fibers \cite{sf2, sf3}.

\textbf{First proof of Theorem~\ref{TCtopcont}:}
According to Theorem~\ref{Tcombtop}, $B_\top(c)$ is determined from
$B_\comb(\theta)$ when $c\in\partial\M$ is in the impression of
$\r_\sM(\theta)$. By the proof of that Theorem, $B_\top(c)$ is constant
on the closures of interior components, since they belong to a primitive
Mandelbrot set or they are pure satellite renormalizable 
(non-hyperbolic components belong to non-trivial fibers \cite{dcc, sf2},
so they are infinitely renormalizable). Now $B_\comb(\theta)$ is
continuous on a compact domain by assumption,
so it is uniformly continuous.

1. Suppose $c\in\partial\M$ is non-renormalizable or finitely renormalizable.
According to Yoccoz \cite{hy3}, the fiber of $c$ is trivial; there is a
sequence of open neighborhoods shrinking to $c$, which are bounded by parts
of periodic rays, equipotential lines, and arcs within hyperbolic components.
More specifically:\\
a) If $c$ is not on the boundary of a hyperbolic component,
$\M\setminus\{c\}$ has finitely many branches, and $c$ is approximated by
separating roots on every branch; their external rays are connected with
equipotential lines.\\
b) If $c$ is a Siegel or Cremer parameter, or the root of the main
cardioid, use two external rays landing at different satellite roots,
an equipotential line, and an arc within the component.\\
c) If $c$ is a satellite bifurcation point, apply this construction
symmetrically for both components.\\
d) If $c$ is a primitive root, combine b) for the component and a) for
the vein before it.\\
The rays and equipotential lines are used to visualize a bounded open
neighborhood in $\C$; the relatively open neighborhood in $\M$ is
determined by pinching points and interior arcs.
When $c_n\to c$, the neighborhoods can be chosen such that the angles
of the rays converge to external angles of $c$,
so $B_\top(c_n)\to B_\top(c)$.

2. Suppose $c\in\partial\M$ is infinitely pure satellite renormalizable.
It is characterized by a sequence of pure satellite components, each one
bifurcating from the previous one, with periods $p_k\to\infty$ and
external angles $\theta_k^-\nearrow\theta_-$ and
$\theta_k^+\searrow\theta_+$\,. There are two cases \cite{wje}:\\
a) If $p_{k+1}=2p_k$ for $k\ge K$, $c$ belongs to the fiber of a tuned
Feigenbaum point: we have $\theta_-<\theta_+$ and a sequence of roots
with external angles $\tilde\theta_k^-\searrow\theta_-$ and
$\tilde\theta_k^+\nearrow\theta_+$\,. When $c_n\to c$, there are angles
$\theta_n$ with $B_\top(c_n)=B_\comb(\theta_n)$ and
$\theta_n\in[\theta_{k_n}^-\,,\,\tilde\theta_{k_n}^-]$ or
$\theta_n\in[\tilde\theta_{k_n}^+\,,\,\theta_{k_n}^+]$. Now
$k_n\to\infty$, so $B_\top(c_n)\to0=B_\top(c)$.\\
b) If $p_{k+1}\ge3p_k$ for infinitely many values of $k$, we have
$\theta_-=\theta_+$ and $\M$ does not continue beyond the fiber of $c$.
When $c_n\to c$, there are angles
$\theta_n$ with $B_\top(c_n)=B_\comb(\theta_n)$ and
$\theta_n\in[\theta_{k_n}^-\,,\,\theta_{k_n}^+]$ with
$k_n\to\infty$, so $B_\top(c_n)\to0=B_\top(c)$ again.

3. Suppose $c\in\partial\M$ is infinitely renormalizable and belongs to
a primitive Mandelbrot set $\M_p$\,. When $c_n\to c$,
$c_n\in\M\setminus\M_p$\,, these points belong to decorations of order
tending to $\infty$, since there are only finitely many decorations
of any bounded order, and these have a positive distance to $c$. By
the Douady substitution, or the recursive construction of decorations, the
corresponding intervals of angles tend to length $0$. Since
$B_\comb(\theta)$ is uniformly continuous, and constant on the angles
of rays bounding the decorations, we have $B_\top(c_n)\to B_\top(\M_p)$.
Note that we do not need to show that the connecting Misiurewicz points
tend to $c$. \mybox

There is a closed equivalence relation on $S^1$, such that
$\theta_1\sim\theta_2$ if the corresponding parameter rays land
together or their impressions form the same fiber. (The latter
situation would be a tuned image of the real case sketched in
Figure~\ref{FFB}, so both rays land or none does.) Now $S^1\!/\!\sim$
is a locally connected Hausdorff space, which is the boundary of the
abstract Mandelbrot set \cite{dcc, sf2}. The natural projection
$\pi_{\scriptscriptstyle S}:S^1\to S^1\!/\!\sim$ is continuous. By collapsing
non-trivial fibers, a continuous projection
$\pi_\sM:\partial\M\to S^1\!/\!\sim$ is obtained. (In Figure~\ref{FFB},
it would map any of the latter configurations to the first one.)

\textbf{Second proof of Theorem~\ref{TCtopcont}:}
Now there is a factorization $B_\comb=B\circ\pi_{\scriptscriptstyle S}$ with
a continuous map $B:S^1\!/\!\sim\,\to\,[0,\,1]$, since $\theta_1\sim\theta_2$
implies $B_\comb(\theta_1)=B_\comb(\theta_2)$: Siegel and Cremer parameters
are not biaccessible; parameter rays with non-trivial impressions belong to
infinitely renormalizable fibers, and $B_\comb(\theta)$ is constant on the
corresponding angles; and $\M$ is locally connected at the remaining
parameters $c$, as is $\K_c$ at $c$.
By Theorem~\ref{Tcombtop} and assuming Conjecture~\ref{Ccombcont} we have
$B_\top=B\circ\pi_\sM$ on $\partial\M$ as a composition of continuous maps.
The interior of $\M$ is treated as in the first proof above. \mybox

\section{Biaccessibility of the Mandelbrot set} \label{5}
For the Mandelbrot set, topological biaccessibility is defined by
parameter rays landing together. Combinatorial biaccessibility is obtained
from the equivalence relation defining the abstract Mandelbrot set
$S^1\!/\!\sim$ \cite{dcc, tgdr}; a pair of angles approximated by pairs
of equivalent periodic angles is equivalent as well, and the angles are
combinatorially biaccessing. The rays may fail to land together as
illustrated in Figure~\ref{FFB}. By the Yoccoz Theorem~\cite{hy3}, this can
happen only for infinitely renormalizable angles, which are negligible in
terms of Hausdorff dimension \cite{mlr}. Zakeri \cite{zerm, zbpm} has shown
that the biaccessing angles of $\M$ or $\M\cap\R$ have Hausdorff dimension
1 and Lebesgue measure 0. According to Bruin--Schleicher \cite{bks, bshb},
the dimension is $<1$ when a neighborhood of $c=-2$ or $\theta=1/2$ is
excluded.

\begin{prop}[Biaccessibility dimension of arcs of $\M$,
generalizing Tiozzo]\label{Parcpd}
Suppose $c'\prec c''$ do not belong to the same primitive Mandelbrot set;
$c''$ is biaccessible or of $\beta$-type. Consider the external angles of
the regulated arc $[c',\,c'']\subset\M$. Their Hausdorff dimension is related
to the dynamic one by $\dim\gamma_\sM^{-1}[c',\,c'']=B_\top(c'')$.
\end{prop}

\textbf{Proof:} 1. Suppose that $c''$ is a non-renormalizable Misiurewicz
point. To obtain $\dim\gamma_\sM^{-1}[c',\,c'']\le B_\top(c'')$, note that
$\gamma_\sM^{-1}[c',\,c'']\subset\gamma_{c''}^{-1}[\alpha_{c''}\,,\,c'']$\,:
any biaccessible parameter $c\in[c',\,c'']$ is either a satellite root or
approximated by roots, and there are corresponding characteristic points
in $[\alpha_{c''}\,,\,c'']$\,; their limit is biaccessible since $\K_{c''}$
is locally connected.

To show $\dim\gamma_\sM^{-1}[c',\,c'']\ge B_\top(c'')$, choose
non-renormalizable Misiurewicz points $c_n$ with $c'\preceq c_n\prec c''$
such that $c_n\nearrow c''$, and such that the preperiodic point $z_c''$
corresponding to $c''$ moves holomorphically for $c$ in the wake of $c_1$\,.
Due to non-renormalizability, the Hubbard tree for $c_n$ is covered by
iterates of the arc $[c_n\,,\,z_{c_n}'']\subset\K_{c_n}$\,, and the related
maps of angles are piecewise-linear, so
$\dim\gamma_{c_n}^{-1}[c_n\,,\,z_{c_n}'']=B_\top(c_n)$. We have
$\dim\gamma_\sM^{-1}[c',\,c'']\ge\dim\gamma_\sM^{-1}[c_n,\,c'']\ge
B_\top(c_n)$, since the angles of $[c_n\,,\,z_{c_n}'']\subset\K_{c_n}$
remain biaccessing for parameters $c\succeq c_n$ according to
Proposition~\ref{Pmonocomb}; assuming that the corresponding parameter rays
do not land at $[c_n\,,\,c'']\subset\M$ gives
a contradiction at least for angles that are not infinitely renormalizable.
The proof is completed by $\lim B_\top(c_n)=B_\top(c'')$.

2. Suppose that $c''$ is a non-renormalizable parameter and approximate it
with non-renormalizable Misiurewicz points $c'\prec c_n''\prec c''$. Then
we have $[c',\,c'')=\bigcup\,[c',\,c_n'']$ and
$\dim\gamma_\sM^{-1}[c',\,c'']=\lim\dim\gamma_\sM^{-1}[c',\,c_n'']
=\lim B_\top(c_n'')=B_\top(c'')$ by monotonicity and continuity according
to Theorem~\ref{Tbev}.

3. Suppose $c''$ is primitive renormalizable but not pure satellite
renormalizable, so it belongs to a maximal primitive Mandelbrot set of
period $p$. The argument of case~2 applies to the corresponding root
$c_*''$ as well. We have
$\dim\gamma_\sM^{-1}[c',\,c'']=\dim\gamma_\sM^{-1}[c',\,c_*'']
=B_\top(c_*'')=B_\top(c'')$ since
$\dim\gamma_\sM^{-1}[c_*'',\,c'']\le\frac1p<B_\top(c_*'')$
by Proposition~\ref{Prenpcf}.2.

4. Suppose $c''$ belongs to a pure satellite Mandelbrot set of maximal
period $p$. By neglecting finitely many angles, we may assume that $c'$
is $p$-renormalizable as well. The proof of cases~1--3 can be copied,
noting that now the iterates of $[c_n\,,\,z_{c_n}'']\subset\K_{c_n}$
cover the small Hubbard tree; in the pure satellite case the
biaccessibility dimension is dominated by the small Julia set according
to $B_\top(c)=\frac1p\,B_\top(\hat c)$.
On the other hand, the statement will be false when $c',\,c''$ belong to
the same primitive Mandelbrot set of period $p$\,: then
$\dim\gamma_\sM^{-1}[c',\,c'']\le\frac1p<B_\top(c'')$.
\mybox

More generally, $c''$ may be an endpoint with trivial fiber, by approximating
it with biaccessible parameters \cite{sf3}.
Tiozzo \cite{tiob} has obtained Proposition~\ref{Parcpd} with his proof of
continuity on principal veins. The following result relates the local
biaccessibility dimension of $\M$ to the dynamic one. Common upper and lower
estimates for both quantities have been determined combinatorially for
$c\to-2$ in \cite{bks, bshb} and at the 0-entropy locus in \cite{bshb}.

\begin{thm}[Biaccessibility dimension of pieces of $\M$]\label{Tpiecepd}
\textsc{Assuming Conjecture~\ref{Ccombcont},
that $B_\comb(\theta)$ is continuous on $S^1$\,:}

$1$. Define a closed piece $\P\subset\M$ by disconnecting $\M$ at finitely
many pinching points. Then the biaccessible parameters $\P'$ have
$\dim\gamma_\sM^{-1}(\P')=\max\{B_\top(c)\,|\,c\in\P\}$.

$2$. Suppose $c\in\M$ is a parameter with trivial fiber, not belonging to the
closure of a hyperbolic component. For any sequence of nested pieces $\P_n$
with $\bigcap\,\P_n=\{c\}$ we have $\lim\dim\gamma_\sM^{-1}(\P_n')=B_\top(c)$.
\end{thm}

\textbf{Proof:} Intersecting $\P$ with the tree of veins according to
Proposition~\ref{Pbetatree} gives a countable family of full veins and
finitely many truncated veins; these arcs $\mathcal{A}$ contain all
biaccessible parameters, and truncated veins within a primitive Mandelbrot
set are negligible. By Proposition~\ref{Parcpd} we have
$\dim\gamma_\sM^{-1}(\P')=\sup\dim\gamma_\sM^{-1}(\mathcal{A})
\le\max\{B_\top(c)\,|\,c\in\P\}$. On the other hand, this maximum is attained
at a parameter $c_0\in\partial\M\cap\P$ approximated by the $\beta$-type
endpoints of veins $\mathcal{A}_n\subset\P$, so
$\dim\gamma_\sM^{-1}(\P')\ge\limsup\dim\gamma_\sM^{-1}(\mathcal{A}_n)
=B_\top(c_0)=\max\{B_\top(c)\,|\,c\in\P\}$.
Item~2 is immediate from item~1 and continuity.
\mybox

\section[Asymptotic self-similarity and local maxima]%
{Asymptotic self-similarity and local maxima\\
   of the biaccessibility dimension} \label{6}
The biaccessibility dimension $B_\comb(\theta)$ shows H\"older asymptotics
at rational angles $\theta_0$ for specific sequences $\theta_n\to\theta_0$\,.
The same techniques give partial results towards the Tiozzo Conjecture
\cite{tiob}. A self-similarity of $B_\comb(\theta)$ for $\theta\to\theta_0$
was considered by Tan Lei and Thurston \cite{gaotl}; the geometric sequences 
suggest a possible scaling factor.

\subsection{H\"older asymptotics} \label{6h}
The following example describes sequences converging to a $\beta$-type
Misiurewicz point. See Figure~\ref{Fssb} for related zooms of
$B_\comb(\theta)$.

\begin{xmp}[Asymptotics at $\theta_0=1/4$]\label{X14}
$a=\gamma_\sM(1/4)$ is the principal $\beta$-type Misiurewicz point in
the $1/3$-limb. It is approached by the sequences $c_n$ and $a_n$ on the vein
according to Example~\ref{Xpriser}. $B_\top(c_n)$ or $B_\top(a_n)$ is given by
$\log\lambda_n/\log2$, where $\lambda_n$ is the largest root of the polynomial
below:\\
Center $c_n$ of lowest period $n\ge 4$\,:
$\,x^{n-2}\cdot(x^3-x^2-2)+(x+1)=0$\\
$\alpha$-type Misiurewicz point $a_n$ of preperiod $n\ge3$\,:
$\,x^{n-2}\cdot(x^3-x^2-2)+2=0$\\
In the other branch at $a_{n-1}$\,, there is a $\beta$-type Misiurewicz point
$b_n$ of preperiod $n\ge4$\,:
$\,x^{n-2}\cdot(x^3-x^2-2)+2(x-1)=0$\\
These polynomials are obtained in Appendix~\ref{Asm}. They imply
$\lambda_n<\lambda_0$ and give geometric asymptotics
$\lambda_n\sim\lambda_0-K\cdot\lambda_0^{-n}$ with $K>0$ in the three cases.
Here $B_\top(a)=B_\comb(1/4)$ is determined from $\lambda_0$\,, which
satisfies $x^3-x^2-2=0$.
\end{xmp}

For any $\beta$-type Misiurewicz point $a$, $f_a(z)$ maps the arc
$[-\alpha_a\,,\beta_a]\to[\alpha_a\,,\beta_a]$ and this defines sequences
of preimages of $0$ and $\alpha_a$\,, respectively, approaching $\beta_a$
monotonically. There are corresponding sequences of centers $c_n$ and
$\alpha$-type Misiurewicz points $a_n$ approaching $a$. Their critical
orbit is described as follows \cite{wjm}: it stays close to the orbit of
the distinguished preimage of $\beta_c$ until it is close to $\beta_c$\,,
and then it moves monotonically on $[\alpha_c\,,\beta_c]$ until it meets
$0$ or $\alpha_c$\,. The latter steps are increased with $n$, while the first
part of the orbit is combinatorially independent of $n$.

\begin{prop}[Asymptotics at angles of Misiurewicz points]\label{Pasympbeta}
A $\beta$-type Misiurewicz point $a$ is approached monotonically
by a sequence of centers $c_n$ with period $n$ on the vein before
it, such that there is no lower period behind $c_n$\,. The core
entropy converges geometrically:
$h(c_n)\sim h(a) - K\cdot\lambda_0^{-n}$ with $K>0$ and
$h(a)=\log\lambda_0$\,. The same result holds for a sequence of
$\alpha$-type Misiurewicz points $a_n$ of preperiod $n$.
\end{prop}

\textbf{Proof:} For large $n$, choose $n'=n-n_0$ edges on
$[\alpha_c\,,\beta_c]$ mapped monotonically, such that they are covered by the
image of a unique edge. Label them such that the former edges are $1\dots n'$
and the latter edge is the last one. The Markov matrix $A_n$ has the
characteristic matrix $A_n-xI$ in (\ref{eqmatseqb}), where the off-diagonal
blocks are $0$ except for the last column each. The top left block has
dimension $n'$ and the lower right block is independent of $n$, as is the last
column of the lower left block. (These blocks are different for $c_n$ and
$a_n$\,.)
\be\label{eqmatseqb} A_n-xI \:=\:
\left(\begin{array}{ccccc|ccc}
-x & & & & & & & 1 \\
1 & \ddots & & & & & & \vdots \\
 & \ddots & \ddots& & & & & \vdots \\
 & & \ddots & \ddots & & & & \vdots \\
 & & & 1 & -x & & & 1 \\ \hline
 & & & & \ast & & &  \\
 & & & & \vdots & & B-xI &  \\
 & & & & \ast & & & 
\end{array}\right)
\ee
The characteristic polynomial is obtained by Laplace expansion of the
determinant with upper $n'\times n'$ minors and complementary lower minors.
There are only two non-zero contributions; taking the first $n'$ columns
for the upper minor gives $(-x)^{n'}$ and taking the first $n'-1$ columns
plus the last one gives $\pm(1+x+\dots+x^{n'-1})$, so
\be\pm\det(A_n-xI)=
 x^{n'}\cdot\tilde P(x)+\frac{x^{n'}-1}{x-1}\,\tilde Q(x) \ , \ee
taking care of the relative signs. This polynomial is multiplied with
$x^{n_0}(x-1)$ and regrouped, and the largest common factor is determined. So
we have $h(c_n)=\log\lambda_n$ or $h(a_n)=\log\lambda_n$\,, respectively,
where $\lambda_n$ is the largest root of an equation
\be\label{eqprfasb}  R(x)\cdot\Big(x^n\cdot P(x)-Q(x)\Big)=0 \ .\ee
Here $P(x)$ is monic and $P(x)$ and $Q(x)$ do not have common roots. By
monotonicity and continuity on the vein according to Theorem~\ref{Tbev}
we have $\lambda_n\nearrow\lambda_0$\,, so the roots of $R(x)$ have modulus
$<\lambda_0$ and they are negligible at least for large $n$. Then
$\lambda_n^n\to\infty$ gives $P(\lambda_0)=0\neq Q(\lambda_0)$. Now
$\lambda_0$ is a simple root of $P(x)$ and the largest in modulus; otherwise
for large $n$, (\ref{eqprfasb}) would have a root larger than $\lambda_0$ or a
small circle of such roots around $\lambda_0$\,. Rewriting this equation for
$\lambda_n$ and performing a fixed point iteration according to the
Banach Contraction Mapping Principle gives
\be x=\lambda_0+Q(x)\,\frac{x-\lambda_0}{P(x)}\,x^{-n} \qquad\mbox{and}\qquad
 \lambda_n=\lambda_0+\frac{Q(\lambda_0)}{P'(\lambda_0)}\cdot\lambda_0^{-n}
 +\O(n\lambda_0^{-2n}) \ . \ee
The logarithm provides a corresponding asymptotic formula for
$\log\lambda_n$\,. \mybox

\begin{figure}[h!t!b!]
\unitlength 0.001\textwidth 
\begin{picture}(990, 1055)
\put(10, 740){\includegraphics[width=0.42\textwidth]{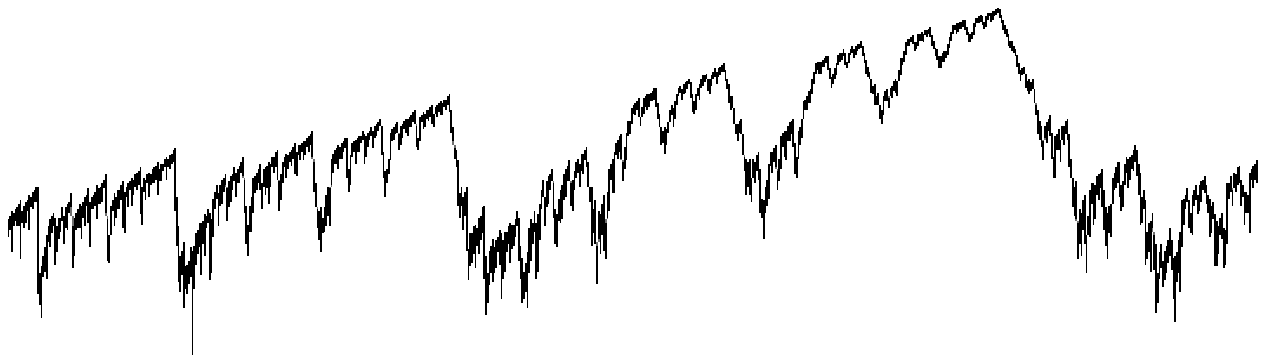}}
\put(570, 740){\includegraphics[width=0.42\textwidth]{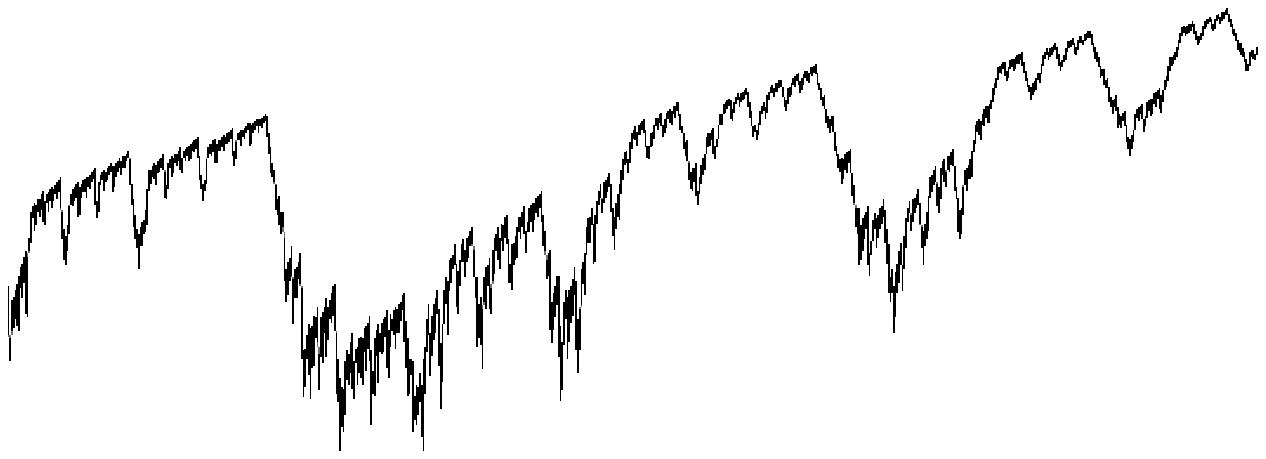}}
\put(10, 370){\includegraphics[width=0.42\textwidth]{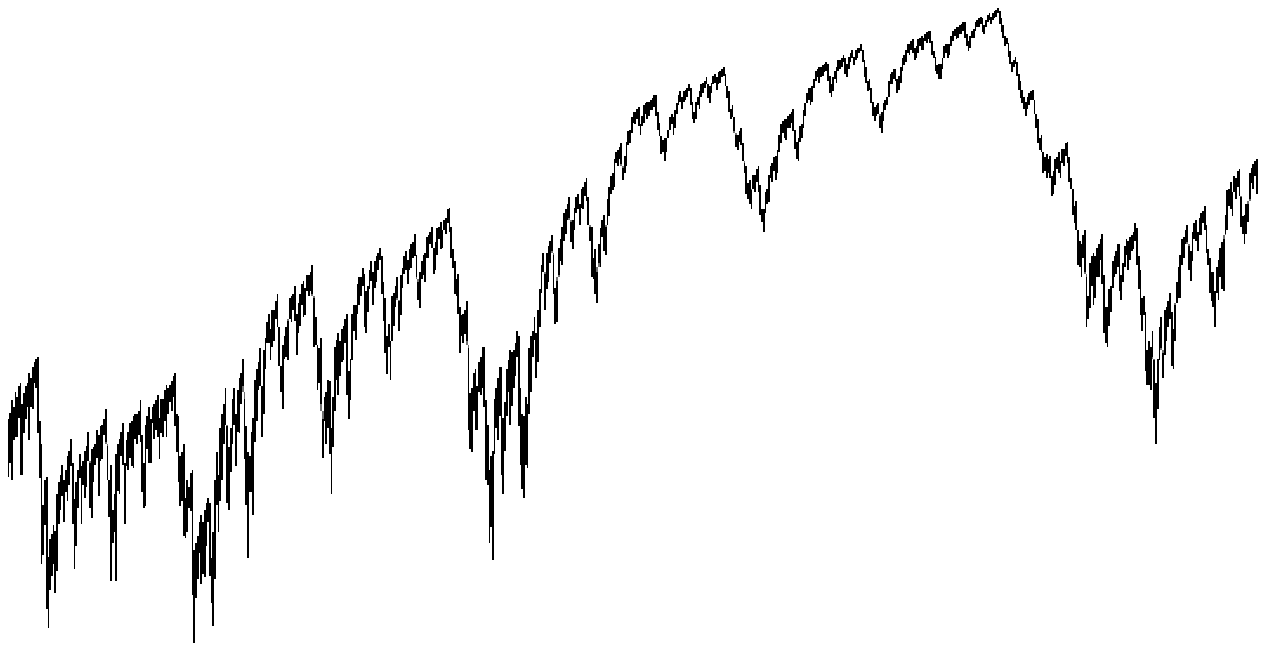}}
\put(570, 370){\includegraphics[width=0.42\textwidth]{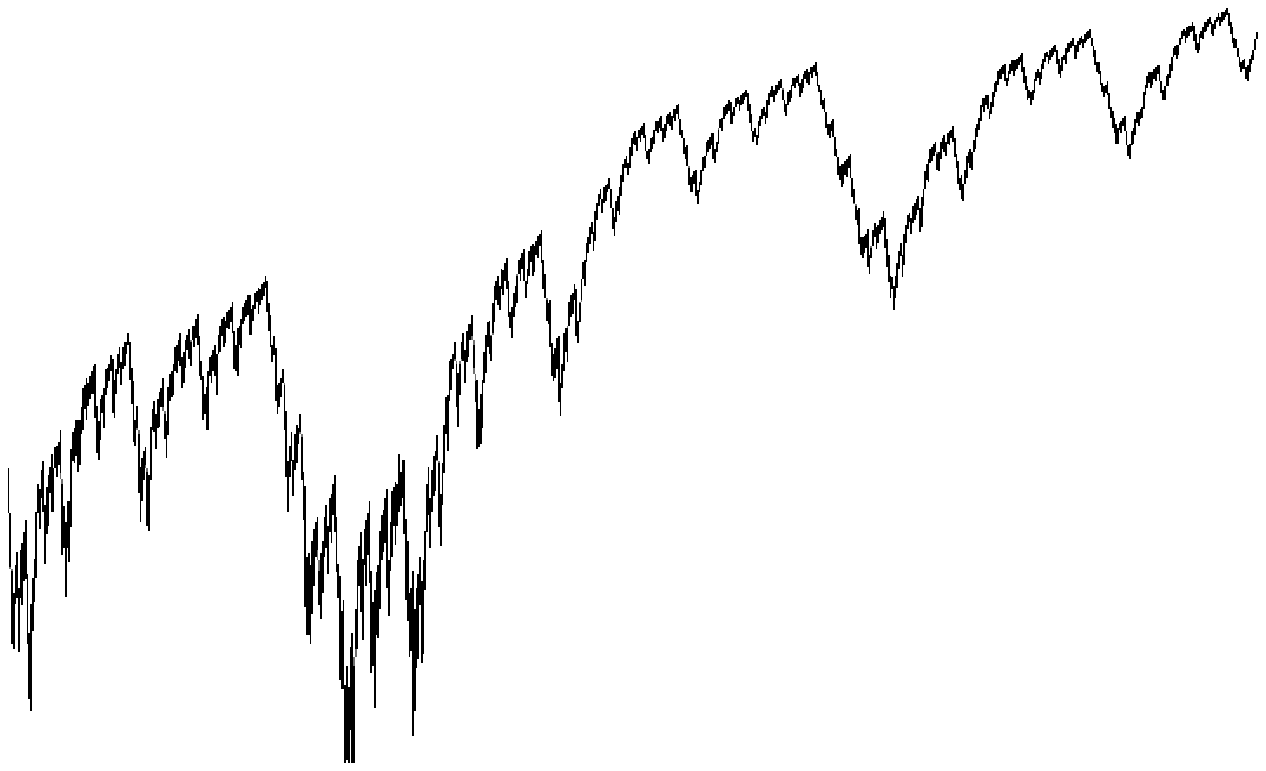}}
\put(10, 0){\includegraphics[width=0.42\textwidth]{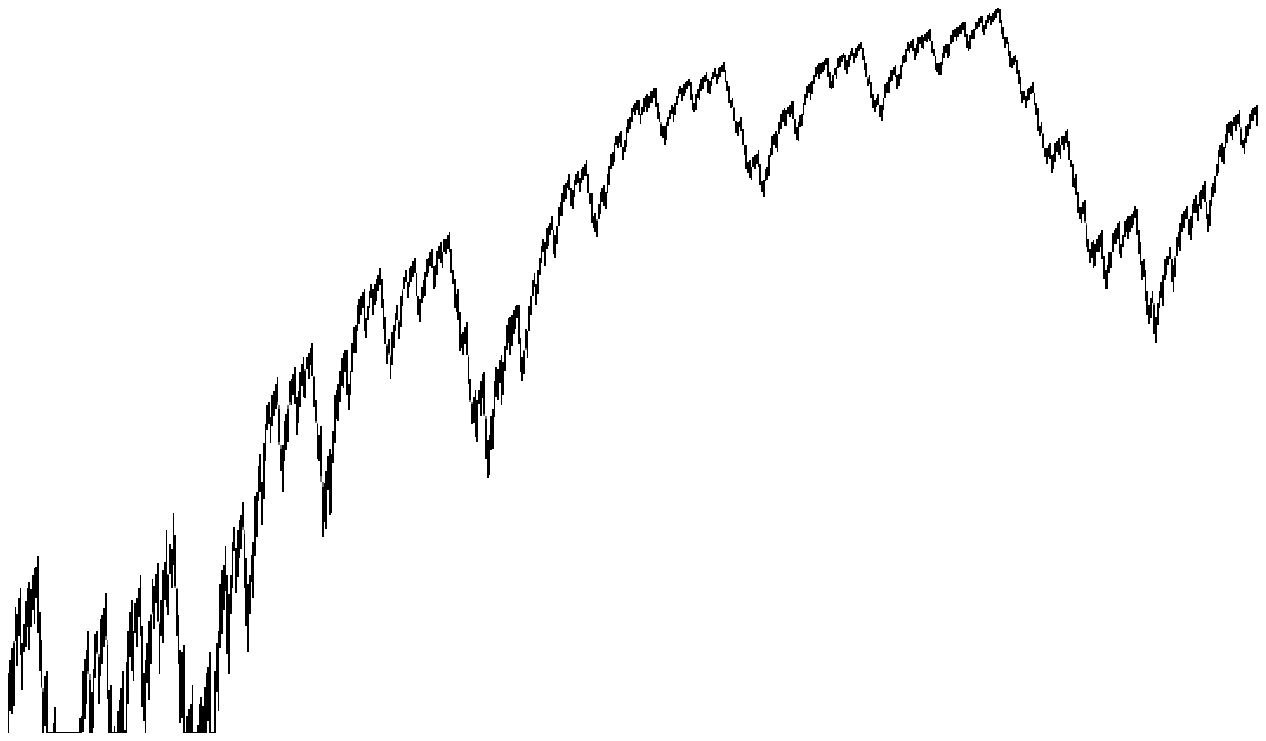}}
\put(570, 0){\includegraphics[width=0.42\textwidth]{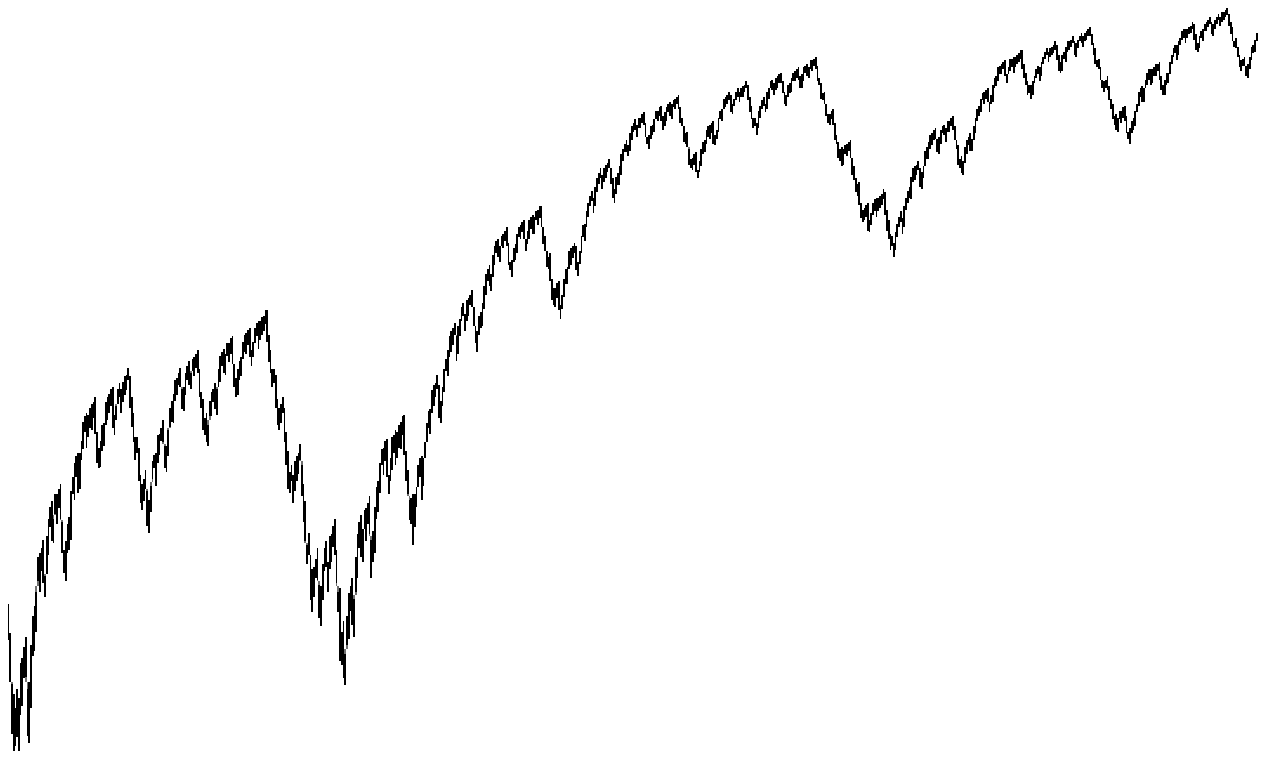}}
\thinlines
\multiput(10, 740)(560, 0){2}{\line(1, 0){420}}
\multiput(10, 1055)(560, 0){2}{\line(1, 0){420}}
\multiput(10, 740)(420, 0){2}{\line(0, 1){315}}
\multiput(570, 740)(420, 0){2}{\line(0, 1){315}}
\multiput(10, 370)(560, 0){2}{\line(1, 0){420}}
\multiput(10, 685)(560, 0){2}{\line(1, 0){420}}
\multiput(10, 370)(420, 0){2}{\line(0, 1){315}}
\multiput(570, 370)(420, 0){2}{\line(0, 1){315}}
\multiput(10, 0)(560, 0){2}{\line(1, 0){420}}
\multiput(10, 315)(560, 0){2}{\line(1, 0){420}}
\multiput(10, 0)(420, 0){2}{\line(0, 1){315}}
\multiput(570, 0)(420, 0){2}{\line(0, 1){315}}
\multiput(500, 157)(0, 370){3}{\makebox(0, 0)[cc]{$\rightarrow$}}
\multiput(500, 342)(0, 370){2}{\makebox(0, 0)[cc]{$\swarrow$}}
\put(35, 1035){\makebox(0, 0)[lt]{$n=0$}}
\put(595, 1035){\makebox(0, 0)[lt]{$n=1$}}
\put(35, 665){\makebox(0, 0)[lt]{$n=2$}}
\put(595, 665){\makebox(0, 0)[lt]{$n=3$}}
\put(35, 295){\makebox(0, 0)[lt]{$n=4$}}
\put(595, 295){\makebox(0, 0)[lt]{$n=5$}}
\end{picture} \caption[]{\label{Fssi}
Consider zooms of $\lambda(\theta)$ centered at $\theta_0=1/6$ with
$\lambda_0=1.521379707$. The width is $0.284\times2^{-n}$ and the height
is $2.185\times\lambda_0^{-n}$\,. A left-sided maximum at $\theta_0$ was
observed by Tan Lei and Thurston \protect\cite{gaotl}. There seems to be
a kind of self-similarity with respect to the combined scaling by $2^2$
and by $\lambda_0^2$, not by $2^1$ and $\lambda_0^1$\,. See
Example~\protect\ref{X16} for the asymptotics of specific sequences; one of
these shows that there is not a right-sided maximum at $\theta_0=1/6$\,.}
\end{figure}

The following example shows that the phenomenon is not limited to $\beta$-type
Misiurewicz points, see Figure~\ref{Fssi}. The higher period of $a$ means that
the period of $c_n$ and the preperiod of $a_n$ grows in steps of $2$:

\begin{xmp}[Asymptotics at $\theta_0=1/6$]\label{X16}
$a=\i=\gamma_\sM(1/4)$ is a Misiurewicz point of preperiod $1$ and period $2$
in the $1/3$-limb of $\M$. It is approached by the sequences $c_n$ and $a_n$
of lowest periods or preperiods on the arc before $a$. $\lambda_n$ is the
largest root of the following polynomial:\\
Center $c_n$ of lowest period $n=5,\,7,\,9,\,\dots$\,:
$\,x^{n-1}\cdot(x^3-x-2)+(x^2+1)=0$\\
$\alpha$-type Misiurewicz point $a_n$ of preperiod $n=3,\,5,\,7,\,\dots$\,:
$\,x^{n-1}\cdot(x^3-x-2)+2=0$\\
In the other branch at $a_{n+1}$\,, consider the $\beta$-type Misiurewicz
point $b_n$ with the preperiod $n=2,\,4,\,6,\,\dots$\,:
$\,x^{n-1}\cdot(x^3-x-2)-2=0$\\
These polynomials give geometric asymptotics
$\lambda_n\sim\lambda_0\pm K\cdot\lambda_0^{-n}$ with $\lambda_n<\lambda_0$ in
the first two cases. Here $B_\top(\i)=B_\comb(1/6)$ is determined from
$\lambda_0$\,, which satisfies $x^3-x-2=0$. The polynomial for $b_n$  is
obtained in Appendix~\ref{Asm}. 
It shows $\lambda_n>\lambda_0$\,, so $B_\comb(\theta)$ does not have a local
maximum at $\theta_0=1/6$ and $B_\top(c)$ does not have a local maximum at
$c=\i$. See Theorem~\ref{Tmonoren}.4 for another application of $b_n$\,.
\end{xmp}

\begin{rmk}[H\"older asymptotics]\label{Rash}
1. On every branch at a Misiurewicz point $a$ with ray period $rp$, there
is a sequence of centers $c_n\to a$, such that the period $n$ is increasing
by $rp$ and there is no lower period between $c_n$ and $a$.
For a suitable choice of angles we have
$\theta_n\sim\theta_0\pm\tilde K\cdot2^{-n}$. In the examples given here and
in Section~\ref{3cce}, and for all $\beta$-type Misiurewicz points according
to Proposition~\ref{Pasympbeta}, we have $B_\top(c_n)=B_\comb(\theta_n)
\sim B_\comb(\theta_0)\pm K'\cdot\lambda_0^{-n}$ with
$B_\top(a)=B_\comb(\theta_0)=\log\lambda_0/\log2$. This confirms the
H\"older exponent $B_\comb(\theta_0)$ for $B_\comb(\theta)$ at
$\theta=\theta_0$ given by Bruin--Schleicher \cite{bshb}, see
Conjecture~\ref{Ccombcont}. Similar statements apply to Misiurewicz points
$a_n\to a$ and to periodic angles $\theta_0$ with $B_\comb(\theta_0)>0$,
see below.

2. Suppose the center parameters $c_n$ are on an arc before or behind the
Misiurewicz point $a$. They converge geometrically as well,
$c_n\sim a+\hat K\cdot\rho_a^{-rj}$ with $n=n_0+rpj$ and the multiplier
$\rho_a=(f_a^p)'(z_a)$ at the periodic point $z_a$ \cite{dhp}.
By the corresponding estimate for fundamental domains
\cite{tls, wjm}, $h(c)$ is H\"older continuous on the arc at $a$
with the optimal exponent $h(a)\cdot p/\log|\rho_a|$\,.
(If this was $>1$, then $h$ would be differentiable with $h'(a)=0$.)
\end{rmk}

\begin{xmp}[Asymptotics at $\theta_0=3/15$]\label{X315}
The primitive Mandelbrot set $\M_4$ in the $1/3$-limb of $\M$ has the external
angles $\theta_0=3/15$ and $4/15$ at its root $a$. It is approached from below
by sequences $c_n$ and $a_n$ of lowest period or preperiod, such that the
growth factor $\lambda_n$ is the largest root of the following polynomial:\\
Center $c_n$ of period $n=7,\,11,\,15,\,\dots$\,:
$\,x^n\cdot(x^4-2x-1)+(x^4+1)=0$\\
$\alpha$-type Misiurewicz point $a_n$ of preperiod $n=3,\,7,\,11,\,\dots$\,:
$\,x^n\cdot(x^4-2x-1)+2=0$\\
In the other branch at $a_{n-1}$\,, consider the $\beta$-type Misiurewicz
point $b_n$ with the preperiod $n=4,\,8,\,12,\,\dots$\,:
$\,x^n\cdot(x^4-2x-1)-2(x^3+x^2+1)=0$\\
These polynomials give geometric asymptotics
$\lambda_n\sim\lambda_0\pm K\cdot\lambda_0^{-n}$ with $\lambda_n<\lambda_0$ in
the first two cases. Here $\lambda_0$ satisfies $x^4-2x-1=0$ and
$B_\top(\M_4)=B_\comb(3/15)=\log\lambda_0/\log2$. The polynomial for $b_n$ is
obtained in Appendix~\ref{Asm}. It shows $\lambda_n>\lambda_0$\,, so
$B_\comb(\theta)$ does not have a left-sided local maximum at $\theta_0=3/15$.
See Theorem~\ref{Tmonoren}.4 for another application of $b_n$\,.
In the $p/q$-sublimb of $\M_4$\,, the $\beta$-type Misiurewicz point $b_n'$ of
lowest preperiod $n=4q-6=2,\,6,\,10,\,\dots$ has a $\lambda_n$ according to
$\,x^n\cdot(x-1)(x^4-2x-1)-2(x^2+1)=0$\,.
\end{xmp}

The sequence of $b_n'$ is generalized to all hyperbolic components as follows:

\begin{prop}[Comparing sublimbs]\label{Pmaxlimbs}
$1$. Consider a hyperbolic component $\Omega$ of period $m$ and the
$\beta$-type Misiurewicz point $b_{p/q}$ of lowest preperiod
$qm-m_0$ in the sublimb with rotation number $p/q$\,. Then the core entropy
$h(b_{p/q})=\log\lambda_q$ is strictly decreasing with $q$
$($and independent of $p)$.

$2$. If $\Omega$ is not of pure satellite type, so
$h(\Omega)=\log\lambda_0>0$, the core entropy converges geometrically:
$h(b_{p/q})\sim h(\Omega) + K\cdot\lambda_0^{-mq}$ with $K>0$.
\end{prop}

\textbf{Proof:} We may assume $m>1$, since the limbs of the main cardioid
are treated explicitly according to Example~\ref{Xprinc}. There is an
$m$-cycle of small Julia sets with $q$ branches at the small $\alpha$-fixed
points. For $q>2$, label the edges such that the first one connects the 
critical value to a small $\alpha$-fixed point, \dots,
edge number $q'=mq-q_0=1+m(q-2)$
is an image at the same small $\alpha$, and the last edge contains the
critical point. When analogous edges are used for $q=2$ as well, the cycle
of small $\alpha$ points is marked in addition, but this does not change the
highest eigenvalue of the Markov matrix $A_q$\,: the $n$-th order preimages of
any point in the Hubbard tree are growing as $\lambda_q^n$\,, as does the
sum of any row of $A_q^n$ for any subdivision into edges.
The characteristic matrix $A_q-xI$ is shown in (\ref{eqmatseqs}), where the
off-diagonal blocks are $0$ except for the last column each. The top left
block has dimension $q'$ and the lower right block is independent of $q$,
as is the last column of the lower left block.
\be\label{eqmatseqs} A_q-xI \:=\:
\left(\begin{array}{ccccc|ccc}
-x & & & & & & & 2 \\
1 & \ddots & & & & & & 0 \\
 & \ddots & \ddots& & & & & \vdots \\
 & & \ddots & \ddots & & & & \vdots \\
 & & & 1 & -x & & & 0 \\ \hline
 & & & & 1 & & &  \\
 & & & & 0 & & B-xI &  \\
 & & & & \vdots & & & 
\end{array}\right)
\ee
The characteristic polynomial is obtained from Laplace cofactor expansion
along the first row and multiplied with $x^{q_0}$ afterwards. We have
$h(b_{p/q})=\log\lambda_q$\,, where $\lambda_q$ is the largest root of an 
equation
\be\label{eqpolylqs}  x^{mq}\cdot P(x)-Q(x)=0 \ee
with $P(x)$ monic. Now $P(x)=x^k\cdot(x-1)\cdot P_0(x)$, where $P_0(x)$
is the characteristic polynomial for the center $c$ of $\Omega$\,: when the
entry $2$ in the first row of $A_2$ is omitted, this removes $Q(x)$ from
the characteristic polynomial. On the other hand, the same Markov matrix will
be obtained for $c$, when edges to preimages of $\beta_c$ are added to the
Hubbard tree according to Lemma~\ref{LHtree}.4. So $P(\lambda_0)=0$ and
$P(x)>0$ for $x>\lambda_0$\,. If $\lambda_0>1$ then $\lambda_0$ is a simple
root of $P(x)$ by Lemma~\ref{Lrenpcf}.1 and Proposition~\ref{Prenpcf}.2.
Now $\lambda_q>\lambda_0$ shows $P(\lambda_q)>0$ and $Q(\lambda_q)>0$.
Setting $x=\lambda_{q+1}$ in the equation (\ref{eqpolylqs}) for
$\lambda_q$ gives
\be \lambda_{q+1}^{mq}\cdot P(\lambda_{q+1})-Q(\lambda_{q+1})
 =\Big(\lambda_{q+1}^{mq}-\lambda_{q+1}^{m(q+1)}\Big)\cdot P(\lambda_{q+1})
 <0 \ , \ee
which implies monotonicity $\lambda_{q+1}<\lambda_q$\,. Set
$\lambda_\ast:=\lim\lambda_q\ge\lambda_0$\,. Assuming $\lambda_0>1$, we have
$\lambda_\ast>1$ and taking $x=\lambda_q$\,, $q\to\infty$ in (\ref{eqpolylqs})
gives $Q(\lambda_\ast)/P(\lambda_\ast)=\infty$, so $P(\lambda_\ast)=0$,
$\lambda_\ast=\lambda_0$\,, and $Q(\lambda_0)\neq0$. Geometric asymptotics
are obtained as in Proposition~\ref{Pasympbeta} since $P'(\lambda_0)>0$.
The asymptotics will be different when $\lambda_0=1$ but we still have
$\lambda_q\searrow\lambda_0$\,: assuming $\lambda_\ast>\lambda_0$ and taking
$q\to\infty$ in (\ref{eqpolylqs}) again gives the contradiction
$P(\lambda_\ast)=0$. \mybox

\subsection{Local maxima} \label{6max}
The wake of a pinching point $c$ is bounded by two parameter rays, such that
it contains all parameters $c'\succ c$. If $c$ is a root or a Misiurewicz
point, these parameter rays have rational angles. When $c$ is a branch point,
the wake consists of subwakes corresponding to the branches behind $c$. In the
wake of a hyperbolic component, the wakes of satellite components may be
called subwakes of the original component, since they correspond to sublimbs.
Motivated by an analogous result for $\alpha$-continued fractions,
Tiozzo has conjectured that on any wake or subwake, or union of neighboring
subwakes, the maximal entropy is attained at the $\beta$-type Misiurewicz
point of lowest preperiod \cite{tiob}. According to
Proposition~\ref{Pbetatree}.1, these points are organized in an infinite tree
of veins, whose branch points are centers and Misiurewicz points. This
suggests to address the Tiozzo Conjecture by considering the two cases
separately. The first case is Proposition~\ref{Pmaxlimbs}.1, but the second
case is still open:

\begin{conj}[Comparing branches]\label{Cmaxbranch}
Consider a branch point $a$ and for each branch behind $a$, the
$\beta$-type Misiurewicz point of lowest preperiod. Of these,
the point with the lowest preperiod has strictly maximal entropy
and biaccessibility dimension.
\end{conj}

This conjecture is equivalent to the Tiozzo Conjecture (except for neighboring
subwakes of a branch point). Figure~\ref{Fssb} suggests a local maximum of
$B_\comb(\theta)$ at $\theta_0=1/4$, which implies a local maximum of
$B_\top(c)$ at $c_0=\gamma_\sM(\theta_0)$. On the other hand,
Figure~\ref{Fssi} implies only a left-sided maximum \cite{gaotl} at the
external angle $\theta_0=1/6$ of $c_0=\i$. In Example~\ref{X16} a sequence
of $\beta$-type Misiurewicz points $b_n$ with angles $\searrow1/6$ and 
$h(b_n)>h(\i)$ is constructed, proving that there is no right-sided
maximum at $1/6$\,. By Theorem~\ref{Tmonoren}.2
or Proposition~\ref{Prenpcf}.3, $B_\comb(\theta)$ cannot be constant on an
interval, since $\beta$-type Misiurewicz points are dense in $\partial\M$
and approximated by maximal-primitive Mandelbrot sets.

\begin{thm}[Maximal entropy] \label{Tmax}
$1$. $B_\top(c)$ has a local minimum at $c_0\in\M$, if and only if
$B_\top(c_0)=0$ or $c_0$ is primitive renormalizable 
and not a maximal-primitive root.\\
If there is a local maximum at $c_0\in\partial\M$, then $c_0$ is an endpoint
and neither simply renormalizable nor on the boundary of the main cardioid.\\
Analogous statements hold for local minima and maxima of
$B_\comb(\theta)$, but there are strict local minima at the inner angles of
branch points in addition,

$2$. Conjecture~$\ref{Cmaxbranch}$ implies: on the interval of any wake or
subwake, $B_\comb(\theta)$ restricted to dyadic angles has a strict absolute
maximum at the dyadic angle of lowest preperiod. The restriction of
$B_\comb(\theta)$ to dyadic angles will be continuous.

$3$. \textsc{Assuming Conjecture~\ref{Ccombcont},
that $B_\comb(\theta)$ is continuous on $S^1$\,:}\\
Conjecture~$\ref{Cmaxbranch}$ implies: on any wake or subwake, 
$B_\top(c)$ has a strict absolute maximum at the $\beta$-type Misiurewicz
point $c_0$ of lowest preperiod. There will be a local maximum at
$c\in\partial\M$, if and only if $c$ is of $\beta$-type.
\end{thm}

\textbf{Proof:} 1. Suppose $B_\top(c_0)>0$. If $c_0$ is primitive
renormalizable and not a maximal-primitive root, then
$B_\top(c)\ge B_\top(c_0)$ for $c$ behind the corresponding
primitive-maximal root $c_1$\,. Otherwise $c_0$ is only finitely
renormalizable, so its fiber is trivial, and for $c$ on the regulated arc
before $c_0$ we have $B_\top(c)<B_\top(c_0)$ by Theorem~\ref{Tmonoren}.2.\\
If $c_0$ is simply renormalizable or on the boundary of the main cardioid,
it is approximated by $\beta$-type Misiurewicz points $c$ with
$B_\top(c)>B_\top(c_0)$. Otherwise it has trivial fiber, so there is a
finite number of rays landing, and this number must be one when there is a
maximum at $c_0$\,: otherwise $B_\top(c)>B_\top(c_0)$ for $c\succ c_0$\,.

2. Fix an interval of angles corresponding to a wake and denote the dyadic
angle of smallest denominator by $\theta_0$\,. Choose another dyadic angle
$\theta$ in the interval. Intersecting the regulated arc from
$c_0=\gamma_\sM(\theta_0)$ to $c=\gamma_\sM(\theta)$ with the tree of veins
according to Proposition~\ref{Pbetatree} gives a finite number of branch
points on a finite number of veins. The first vein is ending at $c_0$ and the
last vein at $c$. For each vein, the endpoint has lowest preperiod in the
subwake of the origin of the vein. The finite collection of endpoints is
compared successively by applying Proposition~\ref{Pmaxlimbs} or
Conjecture~\ref{Cmaxbranch} to the corresponding branch point of the tree.\\
A sequence of dyadic angles $\theta_n\to\theta_0$ has higher preperiods
eventually. There are parameters $a_n\nearrow c_0=\gamma_\sM(\theta_0)$ with
$c_n=\gamma_\sM(\theta_n)$ behind $a_n$\,, so
$B_\top(a_n)<B_\top(c_n)<B_\top(c_0)$. Continuity on the vein according to
Theorem~\ref{Tbev} gives $B_\top(a_n)\to B_\top(c_0)$, which implies
$B_\top(c_n)\to B_\top(c_0)$ and $B_\comb(\theta_n)\to B_\comb(\theta_0)$.

3. Choose a non-simply-renormalizable parameter $c\in\partial\M$ within the
given wake or subwake. $c$ has trivial fiber. If it is not of $\beta$-type,
the regulated arc from $c_0$ to $c$ is meeting a countable family of veins
with endpoints $c_n\to c$. The recursive application of item~2 shows that
$B_\top(c_n)$ is strictly decreasing; we have $B_\top(c_n)\to B_\top(c)$
by assuming continuity of $B_\comb(\theta)$ and Theorem~\ref{TCtopcont}.
So $B_\top(c)<B_\top(c_0)$, there is no local maximum at $c$, and there
is a strict absolute maximum at $c_0$ for the given wake or subwake. \mybox

\subsection{Self-similarity} \label{6sim}
Figures~\ref{Fssb} and~\ref{Fssi} suggest that the graph of the
biaccessibility dimension $B_\comb(\theta)$ may be self-similar; see also
Tan Lei--Thurston \cite{gaotl}. According to the examples and propositions in
Section~\ref{6h}, there are periodic and preperiodic angles $\theta_0$ with
the following property: $B_\comb(\theta_0)=\log\lambda_0/\log2>0$ and there
is a sequence of rational angles $\theta_n\to\theta_0$ with $n$ growing
by $rp$, $\theta_n\sim\theta_0+\tilde K\cdot2^{-n}$, and
$B_\comb(\theta_n)\sim B_\comb(\theta_0)+K'\cdot\lambda_0^{-n}$\,. So we
shall zoom into the graph by scaling with $2^{rpj}$ in the horizontal and by
$\lambda_0^{rpj}$ in the vertical direction:
\begin{itemize}
\item Is there a limit set for $j\to\infty$ in local Hausdorff topology
\cite{tls}?
\item Is it the graph of a function $S(x)=\lim\,\lambda_0^{rpj}\cdot
 \Big(B_\comb(\theta_0+2^{-rpj}x)-B_\comb(\theta_0)\Big)$,
which would be self-similar under combined scaling by $2^{rp}$ and
$\lambda_0^{rp}$? 
\end{itemize}
The latter property can hold only when $\theta_n'-\theta_n=o(2^{-n})$
implies $B_\comb(\theta_n')-B_\comb(\theta_n)=o(\lambda_0^{-n})$, which would
follow from a suitable uniform H\"older estimate for $B_\comb(\theta)$ in
a neighborhood of $\theta_0$ \cite{bshb}.

\begin{figure}[h!t!b!]
\unitlength 0.001\textwidth 
\begin{picture}(990, 420)
\put(10, 0){\includegraphics[width=0.42\textwidth]{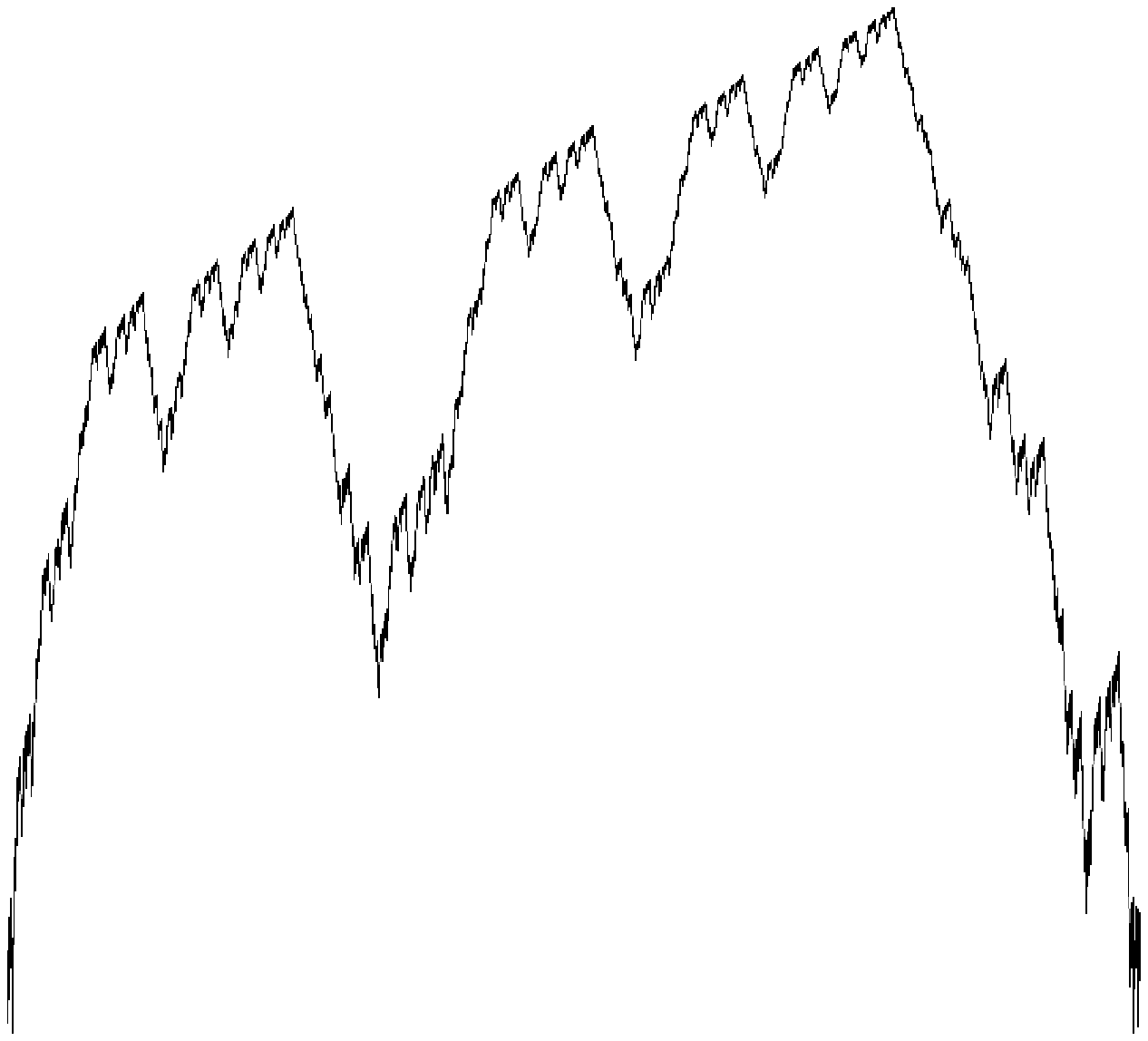}}
\put(570, 0){\includegraphics[width=0.42\textwidth]{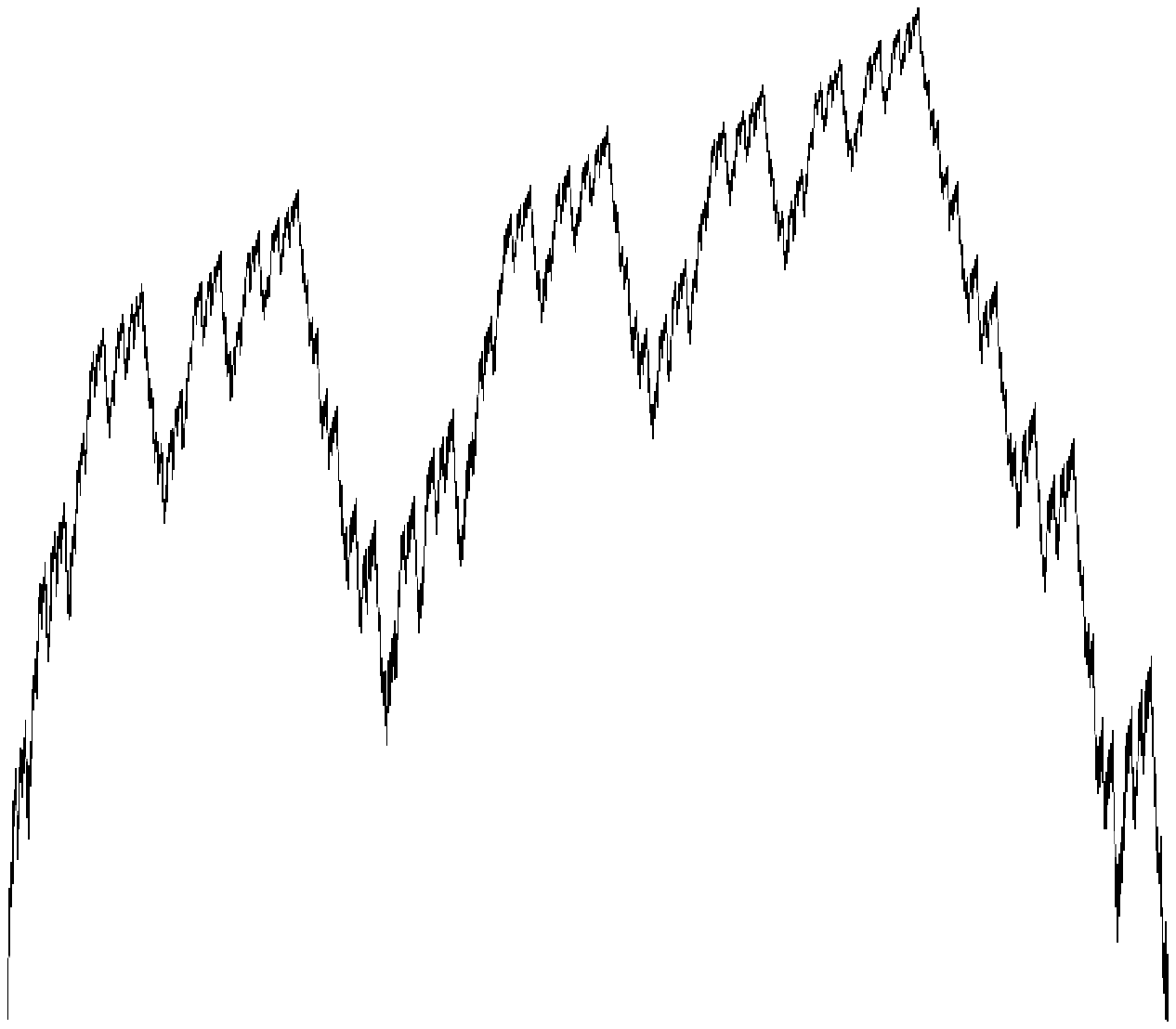}}
\thinlines
\multiput(10, 0)(560, 0){2}{\line(1, 0){420}}
\multiput(10, 420)(560, 0){2}{\line(1, 0){420}}
\multiput(10, 0)(560, 0){2}{\line(0, 1){420}}
\multiput(430, 0)(560, 0){2}{\line(0, 1){420}}
\end{picture} \caption[]{\label{Fm4l}
The graph of $\lambda(\theta)$ representing
$B_\comb(\theta)$ on the intervals $[52/255\,,\,67/255]$ and
$[820/4095\,,\,835/4095]$, which are corresponding to the limbs $1/2$ and
$1/3$ of the primitive Mandelbrot set $\M_4$\,,
see also Figure~\protect\ref{Fm34}.
}
\end{figure}

A different similarity phenomenon is suggested by Figure~\ref{Fm4l}.
It has a qualitative explanation by the linear map between angles of
$\alpha$-type and $\beta$-type Misiurewicz points in different sublimbs.
But this correspondence cannot describe the graph of $B_\comb(\theta)$ in
detail, because the angles of a branch point in one limb may correspond to
endpoints or to several branch points in the other limb.

\appendix
\section{Markov matrices and characteristic polynomials} \label{Asm}
Let us start with the \textbf{proof of Lemma~\ref{LHtree}:}\\
1. When $c$ is preperiodic and the two edges at $z=0$ are numbered first
and second, the Markov matrix $\tilde A$ has the block form given in
(\ref{eqlht1}). Note that the first two rows are identical, because no
marked point is mapped to $0$, so the image of any edge either covers
both edges at $0$ or neither. $u=1$ or $u=0$ indicates whether the edge
intersecting $(0,\,\beta_c]$ is mapped over itself.
\be\label{eqlht1} \tilde A \,=\,
\left(\begin{array}{c|c|ccc}
0 & u &  & C & \\ \hline
0 & u &  & C & \\ \hline
 &  &  &  & \\
B_1 & B_2 &  & D & \\
 &  &  &  &
\end{array}\right)
\quad\to\quad
\left(\begin{array}{c|c|ccc}
0 & 0 &  & 0 & \\ \hline
0 & u &  & C & \\ \hline
 &  &  &  & \\
B_1 & B_3 &  & D & \\
 &  &  &  &
\end{array}\right)
\quad=\quad
\left(\begin{array}{c|ccc}
0 &  & 0 & \\ \hline
 &  &  & \\
B &  & A & \\
 &  &  &
\end{array}\right) \ee
A conjugate matrix is obtained by adding the first column to the second
column and subtracting the second row from the first, denoting $B_1+B_2=B_3$.
In the new block form, the second row and column may be interpreted as the
transitions to and from the joint edge containing $z=0$, so a block equal
to $A$ is identified; now $\tilde A$ has the same eigenvalues as $A$ and an
additional eigenvalue $0$. Note that $B_3$ and $A$ have a unique entry of
$2$. 
While $\tilde A$ corresponds to the usual definition of the Hubbard tree,
$A$ is used for the examples in the present paper; $\tilde A$ should be
used as the adjacency matrix for a subshift of finite type.
--- Now $A$ is irreducible if and only if $\tilde A$ is irreducible, since
total connectivity of the Markov graph is transferred. Consider eigenvalues
of modulus $\lambda$ to show that primitivity is preserved as well.\\
2. The fixed point $\alpha_c$ is not marked when the parameter $c$ in the
$1/2$-limb is not an $\alpha$-type Misiurewicz point. Splitting the edge
containing $\alpha_c$ gives the following transition matrix $\tilde A$, since
the new edges at $\alpha_c$ are mapped over each other. The proof proceeds as
in item~1 with the additional eigenvalue $-1$.
\be \tilde A \,=\,
\left(\begin{array}{c|c|ccc}
0 & 1 &  & C & \\ \hline
1 & 0 &  & C & \\ \hline
 &  &  &  & \\
B_1 & B_2 &  & D & \\
 &  &  &  &
\end{array}\right)
\quad\to\quad
\left(\begin{array}{c|c|ccc}
-1 & 0 &  & 0 & \\ \hline
1 & 1 &  & C & \\ \hline
 &  &  &  & \\
B_1 & B_3 &  & D & \\
 &  &  &  &
\end{array}\right)
\quad=\quad
\left(\begin{array}{c|ccc}
-1 &  & 0 & \\ \hline
 &  &  & \\
B &  & A & \\
 &  &  &
\end{array}\right) \ee
3. When an edge is split by marking a preimage of a marked point, $\tilde A$
will have an additional eigenvalue $0$ by the same proof as for item~1.\\
4. When, e.g., edges towards $\beta_c$ and to a first and second preimage are
attached to the Hubbard tree, the new matrix has the following block
form. $\tilde A$ is reducible, since the original edges are not mapped to
the new ones. The additional eigenvalues are $1$ and $0$.
\be \tilde A \,=\,
\left(\begin{array}{ccc|ccc}
0 & 0 & 0 &  &  & \\
1 & 0 & 0 &  & 0 & \\
0 & 1 & 1 &  &  & \\ \hline
&  &  &  &  & \\
 & B &  &  & A & \\
 &  &  &  &  &
\end{array}\right)
\qquad\qquad\begin{array}{c} \\ \\ \\ \\ \\ \mybox\end{array}\ee

To obtain the characteristic polynomials for a sequence of matrices, the
characteristic matrix $A-xI$ may be transformed to a companion matrix with
$x$-dependent coefficients. When every orbit of edges is passing through a
small subfamily, the rome method can be used \cite{alm}. We shall employ
matching conditions for a piecewise-linear model with expansion rate
$\lambda>1$ again; as the characteristic polynomial, the resulting equation
has the growth factor as its largest positive solution, since all edges have
expressions for the length that are positive for $\lambda>1$, and the
topological entropy is determined uniquely by the combinatorics.
Using the normalization of length $1$ for $[0,\,\pm\alpha_c]$ for parameters
$c$ behind $\gamma_\sM(9/56)$, we have $2/\lambda$ for
$[\alpha_c\,,\,\gamma_c(9/28)]$, $2/\lambda^2$ for
$[\alpha_c\,,\,\gamma_c(9/56)]$, $2/(\lambda-1)$ for
$[\mp\alpha_c\,,\,\pm\beta_c]$, and $2/(\lambda(\lambda-1))$ for
$[-\alpha_c\,,\,\gamma_c(3/4)]$.

\textbf{Computation for Example~\ref{Xpriser} with $q=3$ and for
Example~\ref{X14}:} Pulling back edges towards $\beta_c$ gives the
following matching conditions:
\begin{eqnarray*}
c_n &:& \lambda^3=2+\frac2\lambda+\frac2{\lambda^2}+\dots
 +\frac2{\lambda^{n-4}}+\frac1{\lambda^{n-3}} \\
a_n &:& \lambda=\frac2{\lambda^2}+\dots
 +\frac2{\lambda^{n-1}} \\
b_n &:& \lambda=\frac2{\lambda^2}+\dots
 +\frac2{\lambda^{n-2}}+\frac2{\lambda^{n-1}(\lambda-1)}
\end{eqnarray*}

\textbf{Computation for Example~\ref{X16}:} For the $\beta$-type
Misiurewicz parameter $c=b_n$ consider the edges between $\alpha_c$
and $f_c(c)$. Each edge is mapped to the adjacent one by $f_c^2(z)$,
as is the corresponding branch. The matching condition is
\[ \lambda^2=\frac2\lambda+\frac2{\lambda^3}+\dots
 +\frac2{\lambda^{n-1}}+\frac2{\lambda^{n-1}(\lambda-1)} \ . \]

\textbf{Computation for Example~\ref{X315}:} The Hubbard tree for
$c_n$\,, $a_n$\,, or $b_n$ contains a sequence of small edges scaled by
$\lambda^4$, since they are mapped to the next one by $f_c^4(z)$. The
matching condition for $b_n$ is
\[ \lambda=\frac2{\lambda^2}+\frac2{\lambda^6}+\dots
 +\frac2{\lambda^{n-2}}+\frac2{\lambda^{n-1}(\lambda-1)} \ . \]

\section{Piecewise-linear models and Galois conjugates} \label{Amt}
Suppose $c\in\M\cap\R$ and $h(c)=\log\lambda$ with $\lambda>1$. Then
$f_c(z)=z^2+c$ is semi-conjugate to a piecewise-linear map with slope
$\pm\lambda$, which can be normalized to $g_\lambda(x)=\lambda|x|-1$
\cite{mtimi, dte}. If $c$ is primitive renormalizable, the small Julia sets
are squeezed to points by the semi-conjugation. Given $\lambda>1$, we may
iterate $g_\lambda^n(0)$ to obtain a kneading sequence and an external angle;
in the postcritically finite case, the Hubbard tree is obtained as well, and
the parameter $c$ is found from the real or complex Spider Algorithm
\cite{hss, bct}.

In our real case, the external angle $\theta$ and the kneading sequence $\nu$
are easily converted: an entry $A$ in $\nu$ means that a binary digit of
$\theta$ is changing. From the kneading sequence, $g_\lambda^n(0)$ is obtained
as a polynomial of degree $n-1$ in $\lambda$, which has coefficients $\pm1$.
Here $\lambda|x|-1$ is replaced with $\pm\lambda x-1$ according to the
corresponding entry in $\nu$. If $\theta$ is rational, or $\nu$ is preperiodic
or $\ast$-periodic, $\lambda$ is obtained from a polynomial equation. The
polynomial has coefficients $\pm1$ in the periodic case, and the lower
coefficients are $\pm2,\,0$ in the preperiodic case. For rational and
irrational angles $\theta$, the kneading determinant
$D_c(t)=\sum_{n=0}^\infty\pm t^n$
is holomorphic for $t\in\disk$; its coefficients $\pm1$ correspond to the
binary digits $0$ and $1$ of the external angle $\theta\le1/2$ associated to
the real parameter $c$. This function is related to the generating
function of lap numbers on $[-\beta_c\,,\,\beta_c]$ as follows \cite{mtimi}:
\be\label{eqlapdet} L_c(t)=\sum_{n=1}^\infty L(f_c^n)\,t^{n-1}
 =\frac1{1-t}\,+\,\frac1{(1-t)^2\,D_c(t)} \ . \ee
Now $L_c(t)$ is meromorphic in $\disk$ and $\log L(f_c^n)$ grows as
$n\log\lambda$\,. Thus the smallest pole of $L_c(t)$ and the smallest root of
$D_c(t)$ are located at $t=1/\lambda$. These functions are rational, if and
only if the dynamics is postcritically finite, hyperbolic, or parabolic. In
the $p$-periodic case we have $D_c(t)=P(t)\,/\,(1-t^p)$ with $P(t)$ of degree
$p-1$, and the Douady substitution for any $\hat c\in\M\cap\R$, and for
period doubling in particular, reads
\be\label{eqdsmr} D_c(t)=\frac{P(t)}{1-t^p} \quad\Rightarrow\quad
 D_{c\ast\hat c}(t)=P(t)\,D_{\hat c}(t^p) \quad\mbox{and}\quad
 D_{c\ast(-1)}(t)
 =D_c(t)\,\frac{1-t^p}{1+t^p} \ . \ee
Note that this gives the maximum relation of Lemma~\ref{Lrenpcf}.2 for
every real parameter $\hat c$. Now $D_c(t)$ is constant for parameters $c$
between a center and each of the neighboring parabolic parameters, and the
satellite bifurcation changes it discontinuously at the center. But the
roots of $D_c(t)$ in $\disk$ depend continuously on $c$, and $h(c)$ is
continuous in particular \cite{mtimi}: if a parameter is not hyperbolic or
parabolic, it is approximated by parameters from both sides such that the
first $N$ coefficients of $D_c(t)$ are constant and $N\to\infty$. The same
applies to a primitive parabolic parameter
approximated from below, and the explicit change at a center does not affect
roots in $\disk$. Alternatively, Douady \cite{dte} shows that $h(c)$
cannot have a jump discontinuity: for any $1<\lambda\le2$, the kneading
sequence of the tent map $g_\lambda(x)$ is realized in any full family of
unimodal maps. Thus the monotonic map $h(c)$ is surjective, hence continuous.

For a postcritically finite real polynomial with the topological entropy
$h=\log\lambda$, $\lambda$ is the highest eigenvalue of a non-negative integer
matrix, so it is an algebraic integer and its Galois conjugates are bounded by
$\lambda$ in modulus. Conversely, given $\lambda>1$ with this property,
Thurston \cite{tjack, tedo} constructs a postcritically finite real polynomial
with the topological entropy $h=\log\lambda$. This polynomial will not be
quadratic in general. 
For a postcritically finite real quadratic polynomial,
the Galois conjugates of $\lambda$ must belong to the set $\M_2$ defined
below, which is related to an iterated function system: for complex $\lambda$
with $|\lambda|>1$ consider the holomorphic affine maps
\be\label{eqifs} g_{\lambda,\pm}(z)=\pm\lambda z-1 \qquad
g_{\lambda,\pm}^{-1}(z)=\pm\frac1\lambda(z+1) \qquad
g_\lambda^{-1}(K)=g_{\lambda,+}^{-1}(K) \cup g_{\lambda,-}^{-1}(K) \ . \ee
The corresponding attractor $\K_\lambda$ is non-empty and compact; according
to Hutchinson \cite{hifs} it is the unique compact set with
$g_\lambda^{-1}(K)=K$. It is obtained as the intersection of iterated
preimages of a large disk as well, so $0\in\K_\lambda$ implies that the two
preimages according to (\ref{eqifs}) are never disjoint, and $\K_\lambda$
will be connected. But $\K_\lambda$ may be connected also for parameters
$\lambda$ with $0\notin\K_\lambda$; then it cannot be full.

\begin{figure}[h!t!b!]
\unitlength 0.001\textwidth 
\begin{picture}(990, 785)
\put(10, 0){\includegraphics[width=0.98\textwidth]{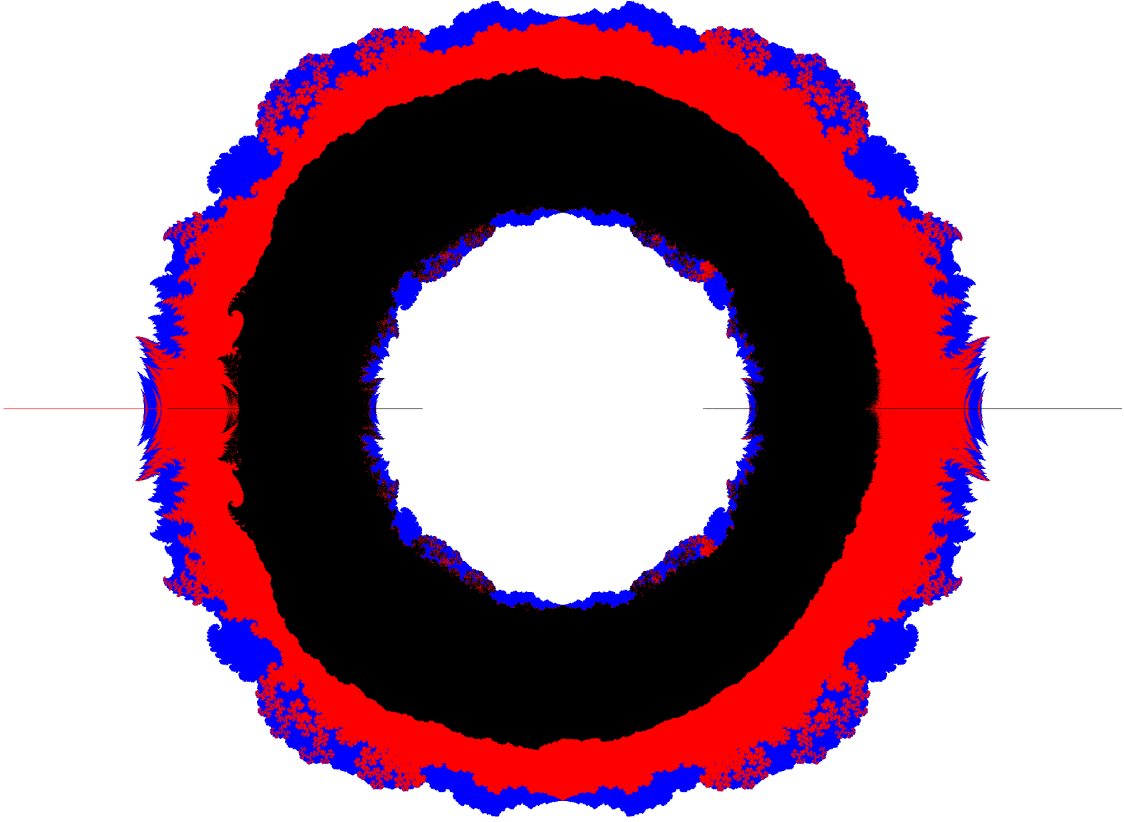}}
\end{picture} \caption[]{\label{Fbbt}
The sets $\M_2\subset\M_1\subset\M_0$ of Proposition~\protect\ref{Pifs} are
shown in black, red (light gray), and blue (dark gray), respectively
\cite{tjack}. The computation is based on polynomials of degree $\le32$.
Due to this restriction, $\M_2$ looks thinner than $\M_1$ for $|\lambda|<1$,
but they are equal there in fact.}
\end{figure}

\begin{prop}[Barnsley--Bousch IFS and Thurston set]\label{Pifs}
Consider the compact sets $\M_2\subset\M_1\subset\M_0$ from
Figure~$\ref{Fbbt}$, which are defined by taking the closure of roots of
families of polynomials,
$\M_i=\overline{\{P^{-1}(0)\,|\,P(\lambda)\mbox{ as follows}\,\}}$\,:\\
For the Barnsley connectedness locus $\M_0$\,, $P(\lambda)$ has coefficients
$\pm1,\,0$ with $P(0)\neq0$.\\
For the Bousch set $\M_1$\,, $P(\lambda)$ has coefficients $\pm1$
from any composition $g_{\lambda,\pm}^n(0)$.\\
For the Thurston set $\M_2$\,, $P(\lambda)$ has coefficients $\pm1$
from a composition $g_{\lambda,\pm}^n(0)$ corresponding to a real
periodic kneading sequence.

$1$. For $|\lambda|>1$ consider the Hutchinson set $\K_\lambda$ of the IFS
$(\ref{eqifs})$. It is connected if and only if $\lambda\in\M_0$\,. We have
$0\in\K_\lambda$ if and only if $\lambda\in\M_1$\,.

$2$. $\M_0$ and $\M_1$ are invariant under inversion and under taking the
$n$-th root. 
Denoting a symmetric annulus by
$\mathcal{A}_r=\overline{\disk_r}\setminus\disk_{1/r}$\,,
we have $\mathcal{A}_{\sqrt{2}}\subset\M_0\subset\mathcal{A}_{2}$ and
$\mathcal{A}_{\sqrt[4]{2}}\subset\M_1$\,.

$3$. $\M_0$ and $\M_1$ are locally connected and connected by
H\"older paths. 

$4$. $\M_2$ is locally connected and for $|\lambda|<1$, the sets $\M_2$ and
$\M_1$ agree. Moreover, $\M_2$ is the closure of the set of Galois conjugates
for $\e{h(c)}$, when all real centers $c$ are considered.
\end{prop}

See \cite{bnplm, tbmksc, tiom} for the \textbf{proof}. Note that Bousch is
using different dynamics, which are conjugate to $\frac1\lambda(z\pm1)$
instead of $\pm\frac1\lambda(z+1)$, but the attractor $\K_\lambda$ will be
the same due to its symmetry. For item~1, the points in $\K_\lambda$ are
parametrized by sequences of signs as $z=\sum\pm\lambda^{-n}$, and the
intersection $g_{\lambda,+}^{-1}(\K_\lambda)
 \cap g_{\lambda,-}^{-1}(\K_\lambda)\subset\K_\lambda$ is considered:
$\K_\lambda$ is connected if the intersection is not empty. If $P(\lambda)=0$
for a polynomial with coefficients $\pm1$, $z=0$ will be periodic under
$g_\lambda^{-1}$\,, so $0\in\K_\lambda$. Moreover, there are
parameters $\lambda\in\partial\M_0\cap\partial\M_1$\,, such that the
intersection is $\{0\}$ and the dynamics on $\K_\lambda$ is
quasi-conformally equivalent to a quadratic polynomial on a dendrite Julia
set \cite{bmplm, ersqs}, cf.~Remark~\ref{Rsedht}.3. See \cite{bbr} for more
detailed pictures of $\M_1$ and for the similarity between $\M_1$ and
$\K_\lambda$ for suitable $\lambda$. 

Due to renormalization, $\M_2$ is invariant under an $n$-th root as well.
In Figure~\ref{Fbbt}, $\M_2$ is restricted strongly for $|\lambda|>1$, and
this is interpreted as follows: for $\lambda>\sqrt2$ and a postcritically
finite $f_c(z)$ with core entropy $h(c)=\log\lambda$, the Galois conjugates
of $\lambda$ are considerably smaller than $\lambda$ \cite{tjack, tedo}.
Now $\M_2$ is pathwise connected as well: for $|\lambda|>1$
consider the parametrization
$\M_2\setminus\overline\disk=\{\lambda\,|\,D_c(1/\lambda)=0,\,-2\le c\le0\}$
by the kneading determinant. All roots of $D_c(t)$ in $\disk$ depend
continuously on the parameter $c$ by the arguments sketched above. In the
case of $|\lambda|<1$, Tiozzo \cite{tiom} shows that $\M_2$ and $\M_1$ agree
by constructing a suitable dense set of polynomials. Irreducible polynomials
yield the statement on Galois conjugates. See \cite{tiob, tiom} for images of
sets, which are analogous to $\M_2$ for other principal veins. Thurston has
started a description of quadratic and higher parameter spaces in terms of 
critical portraits as well, see \cite{gaotl, gao}.

The program Mandel is available from
\href{http://www.mndynamics.com}{www.mndynamics.com}\,.
For dyadic angles you can compute $B_\comb(\theta)$ and zoom into its graph.
Or display the fractals of Appendix~\ref{Amt}.

\begin{thebibliography}{99}\small
\bibitem{alfa} L.~Alsed\`a, N.~Fagella, Dynamics on Hubbard trees,
   Fund.~Math.~\textbf{164}, 115--141 (2000). 
\bibitem{alm} L.~Alsed\`a, J.~Llibre, M.~Misiurewicz,
   \emph{Combinatorial dynamics and entropy in dimension one}, 
   Advanced Series in Nonlinear Dynamics~\textbf{5}, World Scientific 2000.
\bibitem{bbr} J.~C.~Baez, The beauty of roots,
   \href{http://math.ucr.edu/home/baez/roots/beauty.pdf}{slides} and
   blog entry (2011).\\
\href{https://johncarlosbaez.wordpress.com/2011/12/11/the-beauty-of-roots/}%
   {https://johncarlosbaez.wordpress.com/2011/12/11/the-beauty-of-roots/}
\bibitem{bcplm} M.~Baillif, A.~de Carvalho, Piecewise linear model for tree
   maps, Internat.~J.~Bifur.\ Chaos Appl.~Sci.\ Engineering.~\textbf{11},
   3163--3169 (2001).
\bibitem{bmplm} C.~Bandt, On the Mandelbrot set for pairs of linear maps,
   Nonlinearity~\textbf{15}, 1127--1147 (2002). 
\bibitem{bnplm} M.~F.~Barnsley, A.~N.~Harrington, A Mandelbrot set for a pair
   of linear maps, Physica~\textbf{15~D}, 421--432 (1985).
\bibitem{tbmksc} T.~Bousch, Connexit\'e locale et par chemins h\"olderiens
   pour les syst\`emes it\'er\'es de fonctions, unpublished manuscript (1993).
   Available from \href{http://topo.math.u-psud.fr/~bousch/preprints/}%
   {topo.math.u-psud.fr/$\sim$bousch/preprints/} . 
\bibitem{bow1} R.~Bowen, Entropy for group endomorphisms and homogeneous
   spaces, Trans.~AMS~\textbf{153}, 401--414 (1971). 
\bibitem{bd} B.~Branner, A.~Douady, Surgery on complex polynomials, in:
   \emph{Holomorphic dynamics}, X.~Gomez-Mont et al.~eds.,
   LNM~\textbf{1345}, Springer 1988, 11--72. 
\bibitem{bks} H.~Bruin, D.~Schleicher, \emph{Symbolic dynamics of
   quadratic polynomials}, 
   preprint of 2002. \href{https://www.mittag-leffler.se/preprints/}%
   {https://www.mittag-leffler.se/preprints/}
\bibitem{bshb} H.~Bruin, D.~Schleicher, Hausdorff dimension of biaccessible
   angles for quadratic polynomials, preprint (2012).
   \href{http://arXiv.org/abs/1205.2544}{arXiv:1205.2544}
\bibitem{bct} X.~Buff, Cui~G.-Zh., Tan~L., Teichm\"uller spaces and
   holomorphic dynamics, in: \emph{Handbook of Teichm\"uller theory,
   Vol.~IV}, A.~Papadopoulos ed., Soc.~math.~europ., to appear.
   \href{http://www.math.univ-angers.fr/~tanlei/papers/teich.pdf}{Preprint}
   of 2011.
\bibitem{bsoos} S.~Bullett, P.~Sentenac, Ordered orbits of the shift,
   square roots, and the devil's staircase,
   Math.~Proc.~Camb.~Phil.~Soc.~\textbf{115}, 451--481 (1994). 
\bibitem{dhp} A.~Douady, J.~H.~Hubbard, On the dynamics of polynomial-like
   mappings, Ann.~Sci.~\'Ecole Norm.~Sup.~\textbf{18}, 287--343 (1985).
\bibitem{daa} A.~Douady, Algorithms for computing angles in the Mandelbrot
   set, in: \emph{Chaotic Dynamics and Fractals}, M.~F.~Barnsley, S.~G.~Demko
   eds., Notes Rep.~Math.~Sci.~Eng.~\textbf{2}, 155--168 (1986).
\bibitem{dcc} A.~Douady, Descriptions of compact sets in $\C$, in:
   \emph{Topological methods in modern mathematics}, L.~R.~Goldberg ed.,
   Publish or Perish 1993, 429--465.
\bibitem{dte} A.~Douady, Topological entropy of unimodal maps: monotonicity
   for quadratic polynomials, in: \emph{Real and complex dynamical systems}
   (Hiller{\o}d 1993), NATO Adv.~Sci.~Inst.~Ser.~C
   Math.~Phys.~Sci.~\textbf{464}, Kluwer 1995, 65--87.
\bibitem{ersqs} K.~E.~Ero\v{g}lu, S.~Rohde, B.~Solomyak, Quasisymmetric
   conjugacy between quadratic dynamics and iterated function systems,
   Ergod.~Th.~Dynam.~Sys.~\textbf{30}, 1665--1684 (2010).
   \href{http://arXiv.org/abs/0806.3952}{arXiv:0806.3952}
\bibitem{fafr} K.~Falconer, \emph{Fractal geometry, mathematical foundations
   and applications}, Wiley 2003. 
\bibitem{fu} H.~Furstenberg, Disjointness in ergodic theory, minimal sets,
   and a problem in Diophantine approximation,
   Math.~Systems Theory~\textbf{1}, 1--49 (1967). 
\bibitem{gaotl} Gao~Y., Tan~L., On Bill Thurston's binding entropy theory,
   \href{http://milne.ruc.dk/~lunde/MLCQV/Slides/Tan_Lei_entropy.pdf}%
   {Presentation} in Holb{\ae}k 2012.
\bibitem{gao} Gao~Y., \emph{Dynatomic periodic curve and core entropy for
   polynomials}, Ph.D.~thesis Angers University and
   Chinese Academy of Science 2013.\\
   Available from \href{http://tel.archives-ouvertes.fr/tel-00874123}%
   {http://hal.archives-ouvertes.fr}.
\bibitem{grsw} J.~Graczyk, G.~Swi\c{a}tek, Generic hyperbolicity in the
   logistic family, Ann.~Math.~\textbf{146}, 1--56 (1997).
\bibitem{hojo} R.~A.~Horn, C.~R.~Johnson, \emph{Matrix analysis},
   Cambridge University Press 1985.
\bibitem{hy3} J.~H.~Hubbard, Local connectivity of Julia sets and bifurcation
   loci: Three theorems of J.-C.~Yoccoz, in: \emph{Topological methods in
   modern mathematics}, L.~R.~Goldberg~ed., Publish or Perish 1993,
   467--511, 375--378.
\bibitem{hss} J.~H.~Hubbard, D.~Schleicher, The Spider Algorithm, in:
   \emph{Complex Dynamical Systems}, R.~L.~Devaney ed.,
   Proc.\ Symp.\ Appl.\ Math.~\textbf{49}, AMS 1994, 155--180.
\bibitem{hifs} J.~E.~Hutchinson, Fractals and self similarity, Indiana
   Univ.~Math.~J.~\textbf{30}, 713--747 (1981). 
\bibitem{wje} W.~Jung, Renormalization and embedded Julia sets in
   the Mandelbrot set, preprint in preparation (2014).
\bibitem{wjm} W.~Jung, Self-similarity and homeomorphisms
   of the Mandelbrot set, preprint in preparation (2014).
\bibitem{lscr} G.~M.~Levin, S.~van~Strien, Local connectivity of the
   Julia set of real polynomials,\\
   Ann.~Math.~\textbf{147}, 471--541 (1998).
\bibitem{taoli} T.~Li, \emph{A monotonicity conjecture for the entropy of
   Hubbard trees}, Ph.D.~thesis SUNY Stony Brook 2007. Available from
   \href{http://www.math.sunysb.edu/dynamics/theses}%
   {www.math.sunysb.edu/dynamics/theses}.
\bibitem{lmg} J.~Llibre, M.~Misiurewicz, Horseshoes, entropy and periods
  for graph maps, Topology~\textbf{32}, 649--664 (1993). 
\bibitem{l12} M.~Lyubich, Dynamics of quadratic polynomials, I-II,
   Acta~Math.~\textbf{178}, 185--297 (1997).
\bibitem{mlr} A.~Manning, Logarithmic capacity and renormalizability
   for landing on the Mandelbrot set,
   Bull.~London Math.~Soc.~\textbf{28}, 521--526 (1996).
\bibitem{mcr} C.~T.~McMullen, \emph{Complex Dynamics and Renormalization},
   Annals of Mathematics Studies~\textbf{135},
   Princeton University Press 1995.
\bibitem{mesh} Ph.~Meerkamp, D.~Schleicher, Hausdorff dimension and
   biaccessibility for polynomial Julia sets,
   Proc.~Amer.~Math.~Soc.~\textbf{141}, 533--542 (2013).
   \href{http://arXiv.org/abs/1101.4702}{arXiv:1101.4702}
\bibitem{mevst} W.~de Melo, S.~van Strien, \emph{One-dimensional dynamics},
   Springer 1993. 
\bibitem{mtimi} J.~Milnor, W.~Thurston, On iterated maps of the
   interval, in: \emph{Dynamical Systems},\\ 
   J.~C.~Alexander ed., LNM~\textbf{1342}, Springer 1988, 465--563.
\bibitem{mlcj} J.~Milnor, Local connectivity of Julia sets:
   Expository lectures, in:
   \emph{The Mandelbrot Set, Theme and Variations}, Tan~L.~ed.,
   LMS Lecture Notes~\textbf{274}, Cambridge University Press 2000, 67--116.
   \href{http://arxiv.org/abs/math/9207220}{arxiv:math/9207220}
\bibitem{mer} J.~Milnor, Periodic Orbits, External Rays and the Mandelbrot
   Set: An Expository Account, Ast\'erisque \textbf{261}, 277--333 (2000).
\href{http://arxiv.org/abs/math/9905169}{arxiv:math/9905169}
\bibitem{mtm} J.~Milnor, Tsujii's Monotonicity Proof for Real Quadratic Maps,
   unpublished \href{http://www.math.sunysb.edu/~jack/}{manuscript} (2000).
\bibitem{miszepm} M.~Misiurewicz, W.~Szlenk, Entropy of piecewise monotone
   mappings, Studia Math.~\textbf{67}, 45--63 (1980). 
\bibitem{mishpc} M.~Misiurewicz, S.~V.~Shlyachkov, Entropy of piecewise
   continuous interval maps, in \emph{European Conference on Iteration
   Theory 89}, Ch.~Mira ed., World Scientific 1991, 239--245.
\bibitem{chprk} Ch.~Penrose, \emph{On quotients of shifts associated with
   dendrite Julia sets of quadratic polynomials}, Ph.D.~thesis
   University of Coventry 1994.\\
   A related manuscript of 2000 is available from
   \href{http://www.maths.qmul.ac.uk/~csp/}%
   {www.maths.qmul.ac.uk/$\sim$csp/} .
\bibitem{rscr} J.~Riedl, D.~Schleicher, On Crossed Renormalization
   of Quadratic Polynomials, in: \emph{Proceedings of the $1997$ conference
   on holomorphic dynamics}, RIMS Kokyuroku \textbf{1042}, 11--31, Kyoto 1998.
\bibitem{rt} J.~Riedl, \emph{Arcs in Multibrot Sets, Locally
   Connected Julia Sets and Their Construction by Quasiconformal Surgery},
   Ph.D.\ thesis TU~M\"unchen 2000.\\
   Available from \href{http://www.math.sunysb.edu/dynamics/theses}%
   {www.math.sunysb.edu/dynamics/theses} for 2001.
\bibitem{sf2} D.~Schleicher, On Fibers and Local Connectivity of
   Mandelbrot and Multibrot Sets, in:
   \emph{Fractal geometry and applications}, Lapidus, Frankenhuysen eds.,
   Proc.\ Symp.\ Appl.~Math.~\textbf{72}, AMS 2004, 477--512.
   \href{http://arxiv.org/abs/math/9902155}{arXiv:math/9902155}
\bibitem{sf3} D.~Schleicher, On Fibers and Renormalization of Julia Sets
   and Multibrot Sets, preprint (1998).
   \href{http://arxiv.org/abs/math/9902156}{arXiv:math/9902156}.
\bibitem{ser} D.~Schleicher, Rational Parameter Rays of the
   Mandelbrot Set, Ast\'erisque \textbf{261}, 405--443 (2000).
   \href{http://arxiv.org/abs/math/9711213}{arxiv:math/9711213}
\bibitem{shm} Z.~S{\l}odkowsky, Holomorphic motions and polynomial
   hulls, Proc. Am.~Math.~Soc. \textbf{111}, 347--355 (1991).
\bibitem{ssms} S.~Smirnov, On supports of dynamical laminations and
   biaccessible points in polynomial Julia sets, Colloq.~Math.~\textbf{87},
   287--295 (2001). 
\bibitem{soeacc} D.~E.~K.~S{\o}rensen, Accumulation theorems for quadratic
   polynomials, Ergod.~Th.\ Dynam.\ Sys.~\textbf{16}, 1--36 (1996).
\bibitem{soeren} D.~E.~K.~S{\o}rensen, Infinitely renormalizable quadratic
   polynomials, with non-locally connected Julia set,
   J.~Geometric Analysis~\textbf{10}, 169--206 (2000).
\bibitem{sts} P.~\v{S}tefan, A theorem of \v{S}harkovski\brevei~on the
   existence of periodic orbits of continuous endomorphisms of the real line,
   Comm.~Math.~Phys.~\textbf{54}, 237--248 (1977). 
\bibitem{tls} Tan~L., Similarity between the Mandelbrot set and Julia sets,
   Comm.\ Math.\ Phys.~\textbf{134}, 587--617 (1990).
\bibitem{tgdr} W.~Thurston, On the geometry and dynamics of iterated
   rational maps, in: \emph{Complex dynamics: families and friends},
   D.~Schleicher ed., AK Peters 2009, 1--137.
\bibitem{tabstr} W.~Thurston, Core entropy, talk at
   ``Holomorphic Dynamics around Thurston's Theorem'', Roskilde 2010.
   Citation according to the
   \href{http://milne.ruc.dk/~lunde/Holodyn/Abstracts.htm}{Abstract}.
\bibitem{tjack} W.~Thurston, Real polynomial entropy, talk at
   ``Frontiers in Complex Dynamics'', Banff 2011. Citation according to the
   \href{http://www.math.sunysb.edu/jackfest/Videos/}{slides and video.}
\bibitem{tedo} W.~Thurston, Entropy in dimension one, preprint 2011.
   Edited by S.~Koch, to appear.
\bibitem{tiob} G.~Tiozzo, Topological entropy of quadratic polynomials and
   dimension of sections of the Mandelbrot set, preprint (2013).
   \href{http://arXiv.org/abs/1305.3542}{arXiv:1305.3542}
\bibitem{tiom} G.~Tiozzo, Galois conjugates of entropies of real unimodal
   maps, preprint (2013).\\
   \href{http://arXiv.org/abs/1310.7647}{arXiv:1310.7647}
\bibitem{tsu} M.~Tsujii, A simple proof for monotonicity of entropy in the
   quadratic family, Ergodic Th.\ Dynam.~Sys.~\textbf{20}, 925--933 (2000).
\bibitem{zbqj} S.~Zakeri, Biaccessibility in quadratic Julia sets,
   Ergod.~Th.~Dynam.~Sys.~\textbf{20}, 1859--1883 (2000).
\bibitem{zb2} S.~Zakeri, Biaccessibility in Quadratic Julia Sets II:
   The Siegel and Cremer Cases, preprint (1998).
   \href{http://arXiv.org/abs/math/9801150}{arXiv:math/9801150}
\bibitem{zerm} S.~Zakeri, External rays and the real slice of the
   Mandelbrot set, Ergod.~Th.~Dynam.~Sys.~\textbf{23}, 637--660 (2003).
   \href{http://arXiv.org/abs/math/0210382}{arXiv:math/0210382}
\bibitem{zbpm} S.~Zakeri, On biaccessible points of the Mandelbrot set,
   Proc.~Amer.~Math.~Soc.~\textbf{134}, 2239--2250 (2006).
\bibitem{zdun} Anna Zdunik, On biaccessible points in Julia sets of
   polynomials, Fund.~Math.~\textbf{163}, 277--286 (2000). 
\end{thebibliography}
\end{document}